\documentclass[12pt,reqno,a4]{amsart} 

\usepackage[all]{xy}
\usepackage[latin1]{inputenc}  
\usepackage{rotating}

\def\theoremnumberstyle{changesection}
\newcommand{\remembersatz}{\arabic{section}.\arabic{satz}}
\usepackage{tom}

\newcommand{\tom}[1]{} 
\def\TOM{2} 
\newcommand{\problem}[2]{\def\AA{#1}\if\AA\TOM#2\fi} 
\newcommand{\oldversion}[1]{} 

\catcode`@=11
\renewenvironment{enumerate}{%
  \ifnum \@enumdepth >3 \@toodeep\else
      \advance\@enumdepth \@ne
      \edef\@enumctr{enum\romannumeral\the\@enumdepth}\list
      {\csname label\@enumctr\endcsname}{\usecounter
        {\@enumctr}\def\makelabel##1{\hss\llap{\upshape##1}}}\fi
      \itemsep1ex\partopsep1ex\labelsep1ex
}{%
  \endlist
}
\catcode`@=\active%

\renewcommand{\theequation}{\arabic{section}.\arabic{equation}}
\renewcommand{\thesubsection}{\arabic{section}.\alph{subsection}}

\newcommand{\theequationchange}[1]{\addtocounter{equation}{-1}\renewcommand{\theequation}{#1}}
\newcommand{\theequationchangeback}{\renewcommand{\theequation}{\arabic{section}.\arabic{equation}}}
\newcounter{sss}  

\newcommand{\ua}{\underline{a}}
\newcommand{\uai}{\underline{a}_{(0)}}
\newcommand{\uaii}{\underline{a}_{(1)}}
\newcommand{\uaiii}{\underline{a}_{(2)}}
\newcommand{\tua}{\tilde{\underline{a}}}
\newcommand{\tuai}{\tilde{\underline{a}}_{(0)}}
\newcommand{\tuaii}{\tilde{\underline{a}}_{(1)}}
\newcommand{\tuaiii}{\tilde{\underline{a}}_{(2)}}
\renewcommand{\i}{{i}}
\newcommand{\ii}[1]{{#1}}
\renewcommand{\j}{{j}}
\newcommand{\ifix}{{\i,fix}}
\newcommand{\jfix}{{\j,fix}}
\newcommand{\ta}{\tilde{a}_{k,l}^\i}
\newcommand{\tal}{\tilde{a}_{\kappa,\lambda}^\j}
\newcommand{\taj}{\tilde{a}_{k,l}^\j}

\SelectTips{cm}{12}  
\CompileMatrices
\begin{document} 

   \thispagestyle{empty}

   \title{Existence of Curves with Prescribed Topological
     Singularities}
   \author{Thomas Keilen}
   \address{Universit\"at Kaiserslautern\\
     Fachbereich Mathematik\\
     Erwin-Schr\"odinger-Stra\ss e\\
     D -- 67663 Kaiserslautern\\
     e-mail: keilen@mathematik.uni-kl.de}
   \author{Ilya Tyomkin}
   \address{Tel Aviv University\\
     School of Mathematical Sciences\\
     Ramat Aviv\\
     ISR -- Tel Aviv 69978\\
     e-mail: tyomkin@math.tau.ac.il
     }
   \begin{abstract}
     Throughout this paper we study the existence of irreducible
     curves $C$ on smooth 
     projective surfaces $\Sigma$ with singular points of prescribed
     topological types $\ks_1,\ldots,\ks_r$. There are \emph{necessary}
     conditions for the existence of the type 
     $\sum_{i=1}^r \mu(\ks_i)\leq \alpha C^2+\beta C.K+\gamma$ for
     some fixed divisor $K$ on $\Sigma$ and suitable coefficients
     $\alpha$, $\beta$ and $\gamma$, and
     the main \emph{sufficient} condition that we find is of the same
     type, saying it is asymptotically optimal. An important
     ingredient for the proof is a vanishing theorem for
     invertible sheaves on the blown up $\Sigma$ of the form
     $\ko_{\widetilde{\Sigma}}(\pi^*D-\sum_{i=1}^rm_iE_i)$, deduced
     from the Kawamata-Vieweg 
     Vanishing Theorem. A large part of the paper is devoted to the
     investigation of our conditions on ruled surfaces,
     products of elliptic curves, surfaces in $\PC^3$, and
     K3-surfaces.
   \end{abstract}

   \maketitle
   
   \tableofcontents

   \section{Introduction}
   \setcounter{equation}{0}

   \begin{varthm-roman-break}[General Assumptions and Notations]
     Throughout this paper $\Sigma$ will be a smooth projective surface
     over $\C$. 
 
     Given distinct points $z_1,\ldots,z_r\in\Sigma$, we denote by
     $\pi:\widetilde{\Sigma}=\Bl_{\underline{z}}(\Sigma)\rightarrow\Sigma$ 
     the blow up of $\Sigma$ 
     in $\underline{z}=(z_1,\ldots,z_r)$, and the exceptional divisors $\pi^*z_i$ will
     be denoted by $E_i$, $i=1,\ldots,r$. We shall write
     $\widetilde{C}=\Bl_{\underline{z}}(C)$ for the strict
     transform of a curve $C\subset\Sigma$. 

     For any smooth surface $S$ we will denote by $\Div(S)$ the 
     group of divisors on $S$ and by $K_S$ its canonical divisor. If $D$ is any
     divisor on $S$, $\ko_S(D)$ shall be a corresponding invertible
     sheaf.  $|D|_l=\P\big(H^0\big(S,\ko_S(D)\big)\big)$ denotes the
     system of curves linearly equivalent to $D$, while we use the notation
     $|D|_a$ for the system of curves algebraically equivalent to $D$
     (cf.~\cite{Har77} Ex.~V.1.7), that is the reduction of the
     connected component of $\Hilb_S$, 
     the Hilbert scheme of all curves on $S$, 
     containing any curve algebraically equivalent to $D$ (cf.~\cite{Mum66} Chapter 
     15).\tom{\footnote{Note that indeed the reduction of the Hilbert
         scheme gives the Hilbert scheme $\Hilb_S^{red}$ of
         curves on $S$ over reduced base spaces.}}
     We will use the notation
     $\Pic(S)$ for the Picard group of $S$, that is
     $\Div(S)$ modulo linear equivalence (denoted by $\sim_l$),
     $\NS(S)$ for the 
     N\'eron-Severi group, that is $\Div(S)$ modulo algebraic
     equivalence (denoted by $\sim_a$), and $\Num(S)$ for $\Div(S)$ modulo
     numerical equivalence (denoted by $\sim_n$). Note that for all
     examples of surfaces $\Sigma$ 
     which we consider in Section \ref{sec:examples}
     $\NS(\Sigma)$ and $\Num(\Sigma)$ coincide. 

     Given a curve $C\subset\Sigma$ we will write $p_a(C)$ for its
     arithmetical genus and $g(C)$ for the geometrical one.

     Let $Y$ be a Zariski topological space. 
     We say a subset $U\subseteq Y$ is \emph{very general} if it is an 
     at most countable intersection of open dense subsets of
     $Y$. 
     Some statement is said to hold for points
     $z_1,\ldots, z_r\in Y$ (or $\underline{z}\in Y^r$) 
     \emph{in very general position} if there is a suitable very
     general subset $U\subseteq Y^r$, contained in the complement of
     the closed subvariety $\bigcup_{i\not=j}\{\underline{z}\in
     Y^r\;|\; z_i=z_j\}$ of $Y^r$,  such that the statement holds for
     all $\underline{z}\in U$. The main results of this paper will
     only be valid for points in very general position.

     Given distinct points $z_1,\ldots,z_r\in\Sigma$ and non-negative integers
     $m_1,\ldots,m_r$ we denote by
     $X(\underline{m};\underline{z})=X(m_1,\ldots,m_r;z_1,\ldots,z_r)$ 
     the zero-dimensional subscheme of $\Sigma$ defined by the ideal
     sheaf $\kj_{X(\underline{m};\underline{z})/\Sigma}$ with stalks
     \begin{displaymath}
       \kj_{X(\underline{m};\underline{z})/\Sigma,z}=
       \left\{
         \begin{array}{ll}
           \m_{\Sigma,z_i}^{m_i}, & \mbox{ if } z=z_i, i=1,\ldots,r,\\
           \ko_{\Sigma,z}, & \mbox{ else.}
         \end{array}
       \right.
     \end{displaymath}
     We call a scheme of the type $X(\underline{m};\underline{z})$ a
     \emph{generic fat point scheme}.

     For a reduced curve $C\subset\Sigma$ we define the
     zero-dimensional subscheme $X^{es}(C)$ of $\Sigma$ via the ideal sheaf
     $\kj_{X^{es}(C)/\Sigma}$ with stalks
     \begin{displaymath}
       \kj_{X^{es}(C)/\Sigma,z}=I^{es}(C,z)=
       \{g\in\ko_{\Sigma,z}\;|\;f+\varepsilon g \mbox{ is equisingular
         over } \C[\varepsilon]/(\varepsilon^2)\},         
     \end{displaymath}
     where $f\in\ko_{\Sigma,z}$ is a local equation of $C$ at $z$.
     $I^{es}(C,z)$ is called the \emph{equisingularity ideal} of the 
     singularity $(C,z)$, and it is of course $\ko_{\Sigma,z}$
     whenever $z$ is a smooth point. If $x,y$ are local
     coordinates of $\Sigma$ at $z$, then
     $I^{es}(C,z)/\big(f,\frac{\partial f}{\partial x},\frac{\partial
       f}{\partial y}\big)$ can be identified with the tangent space
     of the equisingular stratum in the semiuniversal deformation of
     $(C,z)$. (cf.~\cite{Wah74}, \cite{DH88}, and Definition
     \ref{def:Tsmooth}) We call $X^{es}(C)$ the
     \emph{equisingularity scheme of $C$}.

     If $X\subset\Sigma$ is any zero-dimensional scheme with ideal
     sheaf $\kj_X$ and if $L\subset\Sigma$ is any curve with ideal
     sheaf $\kl$, we define the
     \emph{residue scheme} $X:L\subset\Sigma$ by the ideal sheaf $\kj_{X:L/\Sigma}=\kj_X:\kl$
     with stalks
     \begin{displaymath}
       \kj_{X:L/\Sigma,z}= \kj_{X,z}:l_z,
     \end{displaymath}
     where $l_z\in\ko_{\Sigma,z}$ is a local equation for
     $L$ and ``$:$'' denotes the ideal quotient. This 
     naturally leads to the definition of the \emph{trace scheme} 
     $X\cap L\subset L$ via the ideal sheaf $\kj_{X\cap L/L}$ given by
     the exact sequence
     \begin{displaymath}
       \xymatrix@C0.6cm{
         0\ar[r] & {\kj_{X:L/\Sigma}(-L)}\ar[r]^(0.6){\cdot L} & {\kj_{X/\Sigma}}\ar[r]
         &{\kj_{X\cap L/L}}\ar[r] &0.
         }
     \end{displaymath}

     Given topological singularity types\footnote{For the definition
       of a singularity type and more information see \cite{Los98}
       1.2.} $\ks_1,\ldots,\ks_r$ and a divisor 
     $D\in\Div(\Sigma)$, we denote by $V_{|D|}(\ks_1,\ldots,\ks_r)$
     the locally closed subspace of $|D|_l$ of reduced curves in the
     linear system $|D|_l$ having 
     precisely $r$ singular points of types
     $\ks_1,\ldots,\ks_r$. Analogously,
     $V_{|D|}(m_1,\ldots,m_r)=V_{|D|}(\underline{m})$ denotes the
     locally closed subspace of $|D|_l$ of reduced curves having
     precisely $r$ ordinary singular points of multiplicities
     $m_1,\ldots,m_r$. (cf.~\cite{Los98} 1.3.2) 

     The spaces $V=V_{|D|}(\ks_1,\ldots,\ks_r)$ respectively
     $V=V_{|D|}(\underline{m})$ are the main objects of interest of this 
     paper. We say $V$ is \emph{T-smooth} at $C\in V$ if the germ
     $(V,C)$ is smooth of the (expected) dimension
     $\dim|D|_l-\deg(X)$, where $X=X^{es}(C)$ respectively 
     $X=X(\underline{m};\underline{z})$ with
     $\Sing(C)=\{z_1,\ldots,z_r\}$. By \cite{Los98} Proposition 2.1
     T-smoothness of $V$ at $C$ is implied by the vanishing of
     $H^1\big(\Sigma,\kj_{X/\Sigma}(C)\big)$. 
   \end{varthm-roman-break}

   It is the aim of this paper to give sufficient conditions for the 
   non-emptiness of $V$ in terms of the linear system $|D|_l$ and
   invariants of the imposed singularities. The results are
   generalisations of known results for $\PC^2$, and for an
   overview on these we refer to \cite{Los98}
   Chapter 4.
   
   We basically follow the ideas described in \cite{Los98} 4.1.2. 
   The case of ordinary singularities (Corollary
   \ref{cor:existence-II}) is treated by applying a vanishing 
   theorem for generic fat point schemes (Theorem
   \ref{thm:vanishing}), and the more interesting case of prescribed
   topological types $\ks_1,\ldots,\ks_r$ is then dealt with by
   gluing local equations into a curve with ordinary singularities.
   Upper bounds for the minimal possible degrees of these local
   equations can be taken from the $\PC^2$-case (cf.~\cite{Los98}
   Theorem 4.2).
   
   Thus the main results of this paper are the following theorems and
   their corollaries Corollary \ref{cor:existence-II} and Corollary
   \ref{cor:existence-IV}.

   \renewcommand{\thesatz}{\ref{thm:vanishing}}
   \begin{theorem}
     Let $m_1\geq\ldots\geq m_r\geq 0$ be non-negative integers,
     $\alpha\in\R$ with $\alpha>1$,
     $k_\alpha=\max\big\{n\in\N\;\big|\;n<\tfrac{\alpha}{\alpha-1}\big\}$ and let
     $D\in\Div(\Sigma)$ be a divisor satisfying the following three 
     conditions\tom{\footnote{The proof uses
         Kawamata--Viehweg vanishing which needs characteristic zero
         for the ground field.}}
     \begin{equationlist}
        \theequationchange{\ref{eq:vanishing:1+}}
        \item $(D-K_\Sigma)^2 \geq
          \max\left\{\alpha\cdot\sum\limits_{i=1}^r(m_i+1)^2, 
          (k_\alpha\cdot m_1+k_\alpha)^2\right\}$, 
        \theequationchange{\ref{eq:vanishing:2+}}
        \item $(D-K_\Sigma).B\geq k_\alpha\cdot(m_1+1)$\; for any
          irreducible curve $B$ with $B^2=0$ and $\dim|B|_a>0$, and
        \theequationchange{\ref{eq:vanishing:3+}}
        \item $D-K_\Sigma$ is nef.
        \theequationchangeback
     \end{equationlist}
     Then for $z_1,\ldots,z_r\in\Sigma$ in very
     general position and $\nu>0$ 
     \begin{displaymath}
       H^\nu\left(\Bl_{\underline{z}}(\Sigma),\pi^*D-\sum\limits_{i=1}^r m_iE_i\right)=0.
     \end{displaymath}
     In particular,
     \begin{displaymath}
       H^\nu\big(\Sigma,\kj_{X(\underline{m};\underline{z})/\Sigma}(D)\big)=0.
     \end{displaymath}
   \end{theorem}

   \renewcommand{\thesatz}{\ref{thm:existence-I}}
   \begin{theorem}
     Given $m_1,\ldots,m_r\in\N_0$, not all zero, and $z_1,\ldots,z_r\in\Sigma$,
     $r\geq 1$, in very general position.
     Let $L\in\Div(\Sigma)$ be very ample
     over $\C$, and 
     let $D\in\Div(\Sigma)$ be such that
     \begin{equationlist}
        \theequationchange{\ref{eq:existence-I:1}}
        \item
          $h^1\big(\Sigma,\kj_{X(\underline{m};\underline{z})/\Sigma}(D-L)\big)=0$, and
        \theequationchange{\ref{eq:existence-I:2}}
        \item
          $D.L-2 g(L)\geq m_i+m_j$ for all $i,j$.
        \theequationchangeback
     \end{equationlist}
     Then there exists a curve $C\in |D|_l$ with 
     ordinary singular
     points of multiplicity $m_i$ at $z_i$ for $i=1,\ldots,r$ and no
     other singular points. Furthermore,
     \begin{displaymath}
       h^1\big(\Sigma, \kj_{X(\underline{m};\underline{z})/\Sigma}(D)\big)=0,
     \end{displaymath}
     and in particular, $V_{|D|}(\underline{m})$ is T-smooth at
     $C$.

     If in addition
     \mbox{\rm(\ref{eq:existence-I:3})}
     $D^2>\sum_{i=1}^r m_i^2$,
     then $C$ can be chosen to be irreducible and reduced.
   \end{theorem}

   \renewcommand{\thesatz}{\ref{thm:existence-III}}
   \begin{theorem}[Existence]
     Let $\ks_1,\ldots,\ks_r$ be singularity types, and suppose there
     exists an irreducible curve $C\in|D|_l$ with $r+r'$ ordinary
     singular points $z_1,\ldots,z_{r+r'}$ of multiplicities
     $m_1,\ldots,m_{r+r'}$ respectively as its 
     only singularities such that $m_i=s(\ks_i)+1$, for
     $i=1,\ldots,r$,
     and 
     \begin{displaymath}
       h^1\big(\Sigma,\kj_{X(\underline{m};\underline{z})/\Sigma}(D)\big)=0.
     \end{displaymath}
     Then there exists an irreducible curve $\widetilde{C}\in|D|_l$ with $r$
     singular points of types $\ks_1,\ldots,\ks_r$ and
     $r'$ ordinary singular points of multiplicities
     $m_{r+1},\ldots, m_{r+r'}$ as its only
     singularities.\footnote{Here, of course,
       $\underline{m}=(m_1,\ldots,m_{r+r'})$ and
       $\underline{z}=(z_1,\ldots,z_{r+r'})$. See Definition
       \ref{def:Tsmooth} for the 
       definition of $s(\ks_i)$.}
   \end{theorem}
   \renewcommand{\thesatz}{\remembersatz}

   Of course, combining the vanishing theorem Theorem \ref{thm:vanishing} with 
   the existence theorems Theorem \ref{thm:existence-I} and Theorem
   \ref{thm:existence-III} we get sufficient numerical conditions for
   the existence of curves with certain singularities (see Corollaries 
   \ref{cor:existence-II} and \ref{cor:existence-IV}, and see Section
   \ref{sec:examples} for special surfaces).

   Given any scheme $X$ and any coherent sheaf $\kf$ on $X$, we will
   often write $H^\nu(\kf)$ instead of $H^\nu(X,\kf)$ when no ambiguity
   can arise. Moreover, if $\kf=\ko_X(D)$ is the invertible sheaf corresponding 
   to a divisor $D$, we will usually use the notation $H^\nu(X,D)$
   instead of $H^\nu\big(X,\ko_X(D)\big)$. Similarly when considering tensor products over the
   structure sheaf of some scheme $X$ we may sometimes just write
   $\otimes$ instead of $\otimes_{\ko_X}$.

   Section \ref{sec:vanishing} is devoted to the proof of the
   vanishing theorem Theorem \ref{thm:vanishing}, and Section
   \ref{sec:Geng-Xu} provides an important ingredient in this proof. 
   The following sections Section \ref{sec:existence-I} and Section
   \ref{sec:existence-III} are concerned with the existence theorems
   Theorem \ref{thm:existence-I} and Theorem \ref{thm:existence-III},
   while in Section \ref{sec:examples} we calculate the conditions
   which we have found in the case of ruled surfaces, products of
   elliptic curves, surfaces in $\PC^3$, and K3-surfaces.
   Finally, in the appendix we gather some well known respectively fairly easy
   facts used throughout the paper for the convenience of the reader.


   \section{The Vanishing Theorem}\label{sec:vanishing}
   \setcounter{equation}{0}

   \begin{theorem}\label{thm:vanishing}\index{Theorem!Vanishing}
     Let $m_1\geq\ldots\geq m_r\geq 0$ be non-negative integers,
     $\alpha\in\R$ with $\alpha>1$,
     $k_\alpha=\max\big\{n\in\N\;\big|\;n<\tfrac{\alpha}{\alpha-1}\big\}$ and let
     $D\in\Div(\Sigma)$ be a divisor satisfying the following three 
     conditions\tom{\footnote{The proof uses
         Kawamata--Viehweg vanishing which needs characteristic zero
         for the ground field.}}
     \begin{equationlist}
        \item[eq:vanishing:1+] $(D-K_\Sigma)^2 \geq
          \max\left\{\alpha\cdot\sum\limits_{i=1}^r(m_i+1)^2, 
          (k_\alpha\cdot m_1+k_\alpha)^2\right\}$, 
        \item[eq:vanishing:2+] $(D-K_\Sigma).B\geq k_\alpha\cdot(m_1+1)$\; for any
          irreducible curve $B$ with $B^2=0$ and $\dim|B|_a>0$, and
        \item[eq:vanishing:3+] $D-K_\Sigma$ is nef.
     \end{equationlist}
     Then for $z_1,\ldots,z_r\in\Sigma$ in very
     general position and $\nu>0$ 
     \begin{displaymath}
       H^\nu\left(\Bl_{\underline{z}}(\Sigma),\pi^*D-\sum\limits_{i=1}^r m_iE_i\right)=0.
     \end{displaymath}
     In particular,
     \begin{displaymath}
       H^\nu\big(\Sigma,\kj_{X(\underline{m};\underline{z})/\Sigma}(D)\big)=0.
     \end{displaymath}
   \end{theorem}

   \begin{proof}
     By the Kawamata--Viehweg Vanishing Theorem (cf.~\cite{Kaw82} and
     \cite{Vie82}) it suffices to show that
     $A=\big(\pi^*D-\sum_{i=1}^r m_iE_i\big)-K_{\widetilde{\Sigma}}$ is 
     big and nef, i.~e.~we have to show:
     \begin{enumerate}
        \item[(a)] $A^2>0$, and
        \item[(b)] $A.B'\geq 0$ for any irreducible curve $B'$ in
          $\widetilde{\Sigma}=\Bl_{\underline{z}}(\Sigma)$. 
     \end{enumerate}
     Note that $A=\pi^*(D-K_\Sigma)-\sum_{i=1}^r(m_i+1)E_i$, and thus
     by Hypothesis (\ref{eq:vanishing:1}) we have
     \begin{displaymath}
       A^2=(D-K_\Sigma)^2-\sum_{i=1}^r (m_i+1)^2 >0,
     \end{displaymath}
     which gives condition (a).

     For condition (b) we observe that an irreducible curve $B'$ on
     $\widetilde{\Sigma}$ is either the strict transform of an irreducible 
     curve $B$ in $\Sigma$ or is one of the exceptional curves $E_i$. In the
     latter case we have 
     \begin{displaymath}
       A.B'=A.E_i=m_i+1>0.
     \end{displaymath}
     We may, therefore, assume that $B'=\widetilde{B}$ is the strict transform of
     an irreducible curve $B$ on $\Sigma$ having multiplicity
     $\mult_{z_i}(B)=n_i$ at $z_i$, $i=1,\ldots,r$. 
     Then  
     \begin{displaymath}
       A.B'=(D-K_\Sigma).B-\sum_{i=1}^r(m_i+1)n_i,
     \end{displaymath}
     and thus condition (b) is equivalent to 
     \begin{enumerate}
     \item[(b')] $(D-K_\Sigma).B\geq\sum\limits_{i=1}^r (m_i+1)n_i$.  
     \end{enumerate}
     Since $\underline{z}$ is in very general
     position Lemma \ref{lem:Geng-Xu} applies in view of Corollary
     \ref{cor:verygeneral}. 
     Using the Hodge Index Theorem, Hypothesis (\ref{eq:vanishing:1}), Lemma
     \ref{lem:Geng-Xu}, and the Cauchy-Schwarz Inequality we get the
     following sequence of inequalities:
     \begin{displaymath}
       \renewcommand{\arraystretch}{2}
       \begin{array}{l}
         \bigl((D-K_\Sigma).B\bigr)^2 \;\;\geq \;\;
         \left(D-K_\Sigma\right)^2\cdot B^2 
         \;\;\;\geq \;\; \alpha\cdot\bigl(\sum_{i=1}^r(m_i+1)^2\bigr) \cdot
         \bigl(\sum_{i=1}^r n_i^2-n_{i_0}\big) \\
         \;\; = \;\; \sum_{i=1}^r(m_i+1)^2\cdot\sum_{i=1}^r n_i^2  
         + (\alpha-1)\cdot\!\Bigl(\sum_{i=1}^r(m_i+1)^2\cdot\bigl(\sum_{i=1}^r
         n_i^2-\tfrac{\alpha}{\alpha-1}\cdot n_{i_0}\bigr)\Bigr)\\
         \;\; \geq \;\; \bigl(\sum_{i=1}^r (m_i+1)\cdot n_i\bigr)^2 + 
         (\alpha-1)\cdot\Bigl(\sum_{i=1}^r(m_i+1)^2\cdot\bigl(\sum_{i=1}^r
         n_i^2-\tfrac{\alpha}{\alpha-1}\cdot n_{i_0}\bigr)\Bigr),
       \end{array}
       \renewcommand{\arraystretch}{1.2}
     \end{displaymath}
     where ${i_0}\in\left\{1,\ldots,r\right\}$ is such that    
     $n_{i_0}=\min\{n_i\;|\;n_i\not=0\}$. Since $D-K_\Sigma$ is nef,
     condition (b') is satisfied as soon as we have 
     \begin{displaymath}
       \sum_{i=1}^r n_i^2 \geq \tfrac{\alpha}{\alpha-1}\cdot n_{i_0}.
     \end{displaymath}
     If this is not fulfilled, then
     $n_i<\tfrac{\alpha}{\alpha-1}$ for all $i=1,\ldots,r$, and thus
     \begin{displaymath}
       \sum_{i=1}^r (m_i+1)\cdot n_i\leq k_\alpha\cdot (m_1+1).
     \end{displaymath}
     Hence, for the remaining considerations (b') may be replaced by the worst case
     \begin{displaymath}
       (D-K_\Sigma).B\geq k_\alpha\cdot (m_1+1).
     \end{displaymath}
     Note that since the $z_i$ are in very general position and 
     $z_{i_0}\in B$ we have that $B^2\geq 0$ and $\dim|B|_a>0$ 
     (cf.~Corollary \ref{cor:countable}).  If
     $B^2>0$ then we are done by the 
     Hodge Index Theorem  and Hypothesis (\ref{eq:vanishing:1}), since
     $D-K_\Sigma$ is nef: 
     \begin{displaymath}
       (D-K_\Sigma).B\geq\sqrt{(D-K_\Sigma)^2}\geq
       \sqrt{(k_\alpha\cdot m_1+k_\alpha)^2} \geq k_\alpha\cdot (m_1+1).
     \end{displaymath}
     It remains to consider the case  $B^2=0$ which is covered by 
     Hypothesis (\ref{eq:vanishing:2}).

     \problem{1}{
       For the ``in particular'' part consider the natural short exact sequence
       \begin{displaymath}
         \xymatrix@C0.45cm{
           0\ar[r]& 
           {\Ker}\ar[r] &
           {\kj_{X(m_1;z_1)/\Sigma}\otimes\cdots\otimes\kj_{X(m_r;z_r)/\Sigma}\otimes\ko_\Sigma(D)}\ar[r] &
           {\kj_{X(\underline{m};\underline{z})/\Sigma}\otimes\ko_\Sigma(D)}\ar[r] &
           0.
           }
       \end{displaymath}
       Since
       $H^\nu\big(\Sigma,\bigotimes_{i=1}^r\kj_{X(m_i;z_i)/\Sigma}\otimes\ko_\Sigma(D)\big)= 
       H^\nu\left(\widetilde{\Sigma},\pi^*D-\sum_{i=1}^r m_iE_i\right)$,
       hence vanishes, 
       and since
       $\supp(\Ker)\subseteq\{z_1,\ldots,z_r\}$ (see Lemma 
           \ref{lem:support:1})
       the corresponding long exact cohomology sequence,
       \begin{displaymath}
         \xymatrix@C0.37cm{
           0=H^\nu\bigg(\Sigma,\bigotimes\limits_{i=1}^r\kj_{X(m_i;z_i)/\Sigma}\otimes\ko_\Sigma(D)\bigg)\ar[r]&
           H^\nu\big(\Sigma,\kj_{X(\underline{m};\underline{z})/\Sigma}(D)\big)\ar[r]
           &
           H^{\nu+1}(\Sigma,\Ker)=0,
           }
       \end{displaymath}
       finishes the proof.
       }
     \problem{2}{
       For the ``in particular'' part we just note that
       \begin{displaymath}
         H^\nu\big(\Sigma,\kj_{X(\underline{m};\underline{z})/\Sigma}(D)\big)=
         H^\nu\bigg(\Sigma,\bigotimes\limits_{i=1}^r\kj_{X(m_i;z_i)/\Sigma}\otimes\ko_\Sigma(D)\bigg)=
       \end{displaymath}
       \begin{displaymath}
         H^\nu\left(\widetilde{\Sigma},\pi^*D-\sum_{i=1}^r m_iE_i\right).         
       \end{displaymath}
       }
   \end{proof}

   Choosing the constant $\alpha=2$ in Theorem \ref{thm:vanishing},
   then $\tfrac{\alpha}{\alpha-1}=2$ and thus
   $k_\alpha=1$.  We therefore get the following corollary, which has
   the advantage that the conditions look simpler, and that the
   hypotheses on the ``exceptional'' curves are not too hard.

   \begin{corollary}\label{cor:vanishing}\index{Theorem!Vanishing}
     Let $m_1,\ldots,m_r\in\N_0$, and $D\in\Div(\Sigma)$ be a divisor satisfying the following three
     conditions\tom{\footnote{The proof uses
         Kawamata--Viehweg vanishing which needs characteristic zero
         for the ground field.}}
     \begin{equationlist}
        \item[eq:vanishing:1] $(D-K_\Sigma)^2 \geq 2\cdot\sum\limits_{i=1}^r(m_i+1)^2$, 
        \item[eq:vanishing:2] $(D-K_\Sigma).B>\max\{m_i\;|\;i=1,\ldots,r\}$\; for any
          irreducible curve $B$ with $B^2=0$ and $\dim|B|_a>0$, and
        \item[eq:vanishing:3] $D-K_\Sigma$ is nef.
     \end{equationlist}
     Then for $z_1,\ldots,z_r\in\Sigma$ in very
     general position and $\nu>0$ 
     \begin{displaymath}
       H^\nu\left(\Bl_{\underline{z}}(\Sigma),\pi^*D-\sum\limits_{i=1}^r m_iE_i\right)=0.
     \end{displaymath}
     In particular,
     \begin{displaymath}
       H^\nu\big(\Sigma,\kj_{X(\underline{m};\underline{z})/\Sigma}(D)\big)=0.
     \end{displaymath}
   \end{corollary}

   \begin{remark}\label{rem:vanishing}
     Condition \eqref{eq:vanishing:3+} respectively Condition
     \eqref{eq:vanishing:3} are in several respects 
     ``expectable''. 
     First, Theorem \ref{thm:vanishing} is a corollary of the
     Kawamata--Viehweg Vanishing Theorem, and if we take all $m_i$ to
     be zero, our assumptions should basically be the same,
     i.~e.~$D-K_\Sigma$ nef and big. The latter is more or less
     just \eqref{eq:vanishing:1+} respectively
     \eqref{eq:vanishing:1}. Secondly, we want to apply the  
     theorem to an existence problem.
     A divisor being nef means it is somehow close to being
     effective, or better its linear system is close to being
     non-empty. If we want that some linear system $|D|_l$ contains a 
     curve with certain properties, then it seems not to be so
     unreasonable to restrict to systems where already $|D-K_\Sigma|_l$, or 
     even $|D-L-K_\Sigma|_l$ with $L$ some fixed divisor, is of
     positive dimension, thus nef.
     
     In many interesting examples, such as $\PC^2$, Condition
     \eqref{eq:vanishing:2+} respectively \eqref{eq:vanishing:2}
     turn out to be obsolete or easy to handle. So finally the most restrictive
     obstruction seems to be \eqref{eq:vanishing:1+} respectively
     \eqref{eq:vanishing:1}.  

     If we consider the situation where the largest multiplicity $m_1$
     occurs in a large number, more precisely, if
     $m_1=\ldots=m_{l_\alpha}$ with
     $l_\alpha=\min\big\{n\in\N\;\big|\;\alpha\cdot n\geq k_\alpha^2\big\}$, 
     then Condition \eqref{eq:vanishing:1+} comes down to
     \begin{equationlist}
       \theequationchange{\ref{eq:vanishing:1+}'}
       \item $(D-K_\Sigma)^2 \geq
         \alpha\cdot\sum\limits_{i=1}^r(m_i+1)^2$.
       \theequationchangeback
     \end{equationlist}
   \end{remark}

   \begin{remark}
     Even though we said that Condition (\ref{eq:vanishing:1}) was the
     really restrictive condition we would like to understand better
     what Condition (\ref{eq:vanishing:2}) means. We therefore show in 
     Appendix \ref{sec:condition} that an algebraic system $|B|_a$
     of dimension greater than zero with $B$ irreducible and $B^2=0$
     gives rise to a fibration $f:\Sigma\rightarrow H$ of $\Sigma$
     over a smooth projective curve $H$ whose fibres are just the
     elements of $|B|_a$.
   \end{remark}


   \section{The Lemma of Geng Xu}\label{sec:Geng-Xu}
   \setcounter{equation}{0}

   \begin{lemma}\label{lem:Geng-Xu}
     Let $\underline{z}=(z_1,\ldots,z_r)\in\Sigma^r$ be in very
     general position, $\underline{n}\in\N_0^r$, and let
     $B\subset\Sigma$ be an irreducible curve  
     with $\mult_{z_i}(B)\geq n_i$, then
     \begin{displaymath}
       B^2\geq \sum_{i=1}^r n_i^2 - \min\{n_i\;|\;n_i\not=0\}. 
     \end{displaymath}     
   \end{lemma}

   \begin{eremark}\leererpunkt
     \begin{enumerate}
        \item A proof for the above lemma in the case $\Sigma=\PC^2$ can
          be found in \cite{Xu94} and in the case $r=1$ in \cite{EL93}. Here we 
          just extend the arguments given there to the slightly more
          general situation.     
        \item For better estimates of the self intersection number of
          curves where one has some knowledge on equisingular
          deformations inside the algebraic system see \cite{GS84}. 
        \item With the notation of Lemma \ref{lem:closed} respectively 
          Corollary \ref{cor:verygeneral} the
          assumption in Lemma \ref{lem:Geng-Xu} could be 
          formulated more precisely as ``let $B\subset\Sigma\subseteq\PC^N$ be an irreducible
          curve such that 
          $V_{B,\underline{n}}=\Sigma^r$'', or ``let
          $\underline{z}\in\Sigma^r\setminus V$''.\tom{\footnote{Since $B$ is irreducible, the
              general element in $|B|_a$ will be irreducible. Since
              $V_{B,\underline{n}}=\Sigma^r$ there will 
              be some family of curves in $\Hilb_\Sigma^h$ satisfying the
              requirements of Lemma \ref{lem:deformation}.}}
        \item Note, that one cannot expect to get rid of the ``$-
          \min\{n_i\;|\;n_i\not=0\}$''. E.~g.~$\Sigma=\Bl_z\big(\PC^2\big)$, the
          projective plane blown up in a point $z$, and $B\subset\Sigma$
          the strict transform of a line through $z$. Let now $r=1$, $n_1=1$
          and $z_1\in\Sigma$ be any point. Then there is of course a
          (unique) curve $B_1\in|B|_a$ through $z_1$, but $B^2=0<1=n_1^2$.
     \end{enumerate}
   \end{eremark}

   \begin{varthm-roman}[Idea of the proof]\label{idea:Geng-Xu}
     Set $e_1:=n_1-1$ and 
     $e_i:=n_i$ for $i\not=1$, where
     w.~l.~o.~g.~$n_1=\min\{n_i|n_i\not=0\}$. 
     By assumption there is a family $\{C_t\}_{t\in\C}$ in $|B|_a$
     satisfying the requirements of Lemma \ref{lem:deformation}. 
     Setting $C:=C_0$ the proof is done in three steps:
     \begin{varthm-roman}[Step 1]
       We show that
       $H^0\big(C,\kj_{X(\underline{e};\underline{z})/\Sigma}\cdot\ko_C(C)\big)\not=0$.  
       (Lemma \ref{lem:deformation})
     \end{varthm-roman}
     \begin{varthm-roman}[Step 2]
       We deduce that
         $H^0\big(C,\pi_*\ko_{\widetilde{C}}\big(-\sum_{i=1}^re_iE_i\big)\otimes\ko_C(C)\big)\not=0$.  
         (Lemma \ref{lem:lifting})\tom{\footnote{The main idea for
             this part is the following diagram:
             \begin{displaymath}
               \xymatrix@C0.6cm{
                 0\ar[r]&\Ker(\beta)\ar[r]&
                 H^0\big(\Sigma,\pi_*\ko_{\widetilde{\Sigma}}\big(-\sum_{i=1}^re_iE_i\big)
                 \otimes_{\ko_\Sigma}\ko_C(C)\big)\ar[r]^(0.62)\beta\ar[d]^\alpha &
                 H^0\big(C,\kj_{X(\underline{e};\underline{z})/\Sigma}\cdot\ko_C(C)\big)\ar[r]&
                 0\\
                 & &
                 H^0\big(\Sigma,\pi_*\ko_{\widetilde{C}}\big(-\sum_{i=1}^re_iE_i\big)\otimes_{\ko_C}\ko_C(C)\big)
                 \ar@{-->}^{\overline{\beta}}[ur] & &
                 }
             \end{displaymath}
             i.~e.~the fact that $\Ker(\alpha)\subseteq\Ker(\beta)$,
             or in other words, that $\beta$ factorises over $\alpha$.
             }} 
     \end{varthm-roman}
     \begin{varthm-roman}[Step 3]
       It follows that
         $\deg\big(\pi_*\ko_{\widetilde{C}}\big(-\sum_{i=1}^re_iE_i\big)\otimes\ko_C(C)\big)\geq 0$,
         but this degree is just $C^2-\sum_{i=1}^re_in_i$.
     \end{varthm-roman}
     \hspace*{\fill}\mbox{$ ^{\Box}$}
   \end{varthm-roman}

   \begin{lemma}\label{lem:deformation}
     Given $e_1,\ldots,e_r\in\N_0$, $r\geq 1$.
     Let $\{ C_t\}_{t\in\C}$ be a non-trivial
     family of curves in $\Sigma$ such that 
     \begin{itemize}
        \item $\C\rightarrow\Sigma\: :\: t\mapsto z_{1,t}\in C_t$ is a smooth curve,  
        \item $\mult_{z_{1,t}}(C_t)\geq e_1+1$ for all $t\in\C$, 
        \item $z_2,\ldots,z_r\in C_t$ for any $t\in\C$, and 
        \item $\mult_{z_i}(C_t)\geq e_i$ for all $i=2,\ldots,r$
          and $t\in\C$. 
     \end{itemize}
     Then  for $z_1=z_{1,0}$ 
     \begin{displaymath}
       H^0\big(C,\kj_{X(\underline{e};\underline{z})/\Sigma}\cdot\ko_C(C)\big)\not=0,
     \end{displaymath}
     i.~e.~there is a non-trivial section of the normal bundle of $C$, 
     vanishing at $z_i$ 
     to the order of at least $e_i$ for $i=1,\ldots,r$.\tom{\footnote{Note, that 
           $H^0\big(C,\kj_{X(\underline{e};\underline{z})/\Sigma}\cdot\ko_C(C)\big) = 
           H^0\big(C,\m_{\Sigma,z_1}^{e_1}\cdots\m_{\Sigma,z_r}^{e_r}\cdot\ko_C(C)\big)$.}} 
   \end{lemma}
   \begin{proof}
     We stick to the convention $n_1=e_1+1$ and $n_i=e_i$ for $i=2,\ldots,r$, and we
     set $z_{i,t}:=z_i$ for $i=2,\ldots,r$ and $t\in\C$.
     Let $\Delta\subset\C$ be a small disc around $0$ with coordinate
     $t$, and choose coordinates $(x_i,y_i)$ on $\Sigma$ around $z_i$ such that
     \begin{itemize}
        \item $z_{i,t}=(a_i(t),b_i(t))$ for $t\in\Delta$ with
          $a_i,b_i\in \C\{t\}$, 
        \item $z_i=(a_i(0),b_i(0))=(0,0)$, and
        \item $F_i(x_i,y_i,t)=f_{i,t}(x_i,y_i)\in\C\{x_i,y_i,t\}$, where $C_t=\{f_{i,t}=0\}$
          locally at $z_{i,t}$ (for $t\in\Delta$).
     \end{itemize}
     We view $\{C_t\}_{t\in\Delta}$ as a non-trivial deformation of $C$, 
     which implies that the image of $\frac{\partial}{\partial t}_{|t=0}\in
     T_0(\Delta)$ under the 
     Kodaira-Spencer map is a non-zero section $s$ of
     $H^0\big(C,\ko_{C}(C)\big)$. $s$ is locally at $z_i$ given by
     $\frac{\partial F_i}{\partial t}_{|t=0}$.
     \begin{varthm-roman}[Idea]
     Show that
     $\frac{\partial F_i}{\partial t}_{|t=0}\in
     (x_i,y_i)^{e_i}$, which are the stalks of 
     $\kj_{X(\underline{e};\underline{z})/\Sigma}\cdot\ko_C(C)$\tom{\footnote{$\kj_{X(\underline{e};\underline{z})/\Sigma}
         \cdot\ko_C(C)=\m_{\Sigma,z_1}^{e_1}\cdots\m_{\Sigma,z_r}^{e_r}\cdot\ko_{C}(C)$.}} 
     at the $z_i$, and hence $s$ is actually
     a global section of the subsheaf 
     $\kj_{X(\underline{e};\underline{z})/\Sigma}\cdot\ko_C(C)$.
     \end{varthm-roman}
     Set $\Phi_{i,t}(x_i,y_i):=F_{i,t}(x_i+a_i(t),y_i+b_i(t),t)
     =\sum_{k=0}^{\infty}\varphi_{i,k}(x_i,y_i)\cdot t^k\in\C\{x_i,y_i,t\}$. By
     assumption for any $t\in\Delta$ the multiplicity of $\Phi_{i,t}$ at
     $(0,0)$ is at least $n_i$, i.~e.~$\Phi_{i,t}(x_i,y_i)\in (x_i,y_i)^{n_i}$
     for every fixed complex number $t\in\Delta$. Hence,
     $\varphi_{i,k}(x_i,y_i)\in (x_i,y_i)^{n_i}$ 
     for every $k$.\footnote{See Lemma \ref{lem:powerseries}.}\\
     On the other hand we have
     \begin{displaymath}
       \begin{array}{rcl}
         \varphi_{i,1}(x_i,y_i)&=&
         \frac{\partial\Phi_{i,t}(x_i,y_i)}{\partial t}_{|t=0}\\
         &=&\left\langle \big(\frac{\partial F_i}{\partial x_i}(x_i,y_i,0),
         \frac{\partial F_i}{\partial y_i}(x_i,y_i,0),
         \frac{\partial F_i}{\partial t}(x_i,y_i,0)\big),
         \big(\dot{a}_i(0),\dot{b}_i(0),1\big)\right\rangle\\
         &=&\frac{\partial f_{i,0}}{\partial x_i}(x_i,y_i)\cdot\dot{a}_i(0)+
         \frac{\partial f_{i,0}}{\partial y_i}(x_i,y_i)\cdot\dot{b}_i(0)+
         \frac{\partial F_i}{\partial t}(x_i,y_i,0).
       \end{array}
     \end{displaymath}
     Since $f_{i,0}\in (x_i,y_i)^{n_i}$, we have $\frac{\partial f_{i,0}}{\partial
       x_i}(x_i,y_i),\frac{\partial f_{i,0}}{\partial y_i}(x_i,y_i)\in
     (x_i,y_i)^{n_i-1}$, and hence $\frac{\partial F_i}{\partial
       t}(x_i,y_i,0)\in (x_i,y_i)^{e_i}$. For this note that 
     $\dot{a}_i(0)=\dot{b}_i(0)=0$, if $i\not=1$.
   \end{proof}

   \begin{lemma}\label{lem:lifting}
     Given $e_1,\ldots,e_r\in\N_0$ and $z_1,\ldots,z_r\in\Sigma$,
     $r\geq 1$.
     
     The canonical morphism\tom{\footnote{I.~e.~ 
     $\m_{\Sigma,z_1}^{e_1}\otimes\cdots\otimes\m_{\Sigma,z_r}^{e_r}\otimes\ko_C(C))  
     \longrightarrow \m_{\Sigma,z_1}^{e_1}\cdots\m_{\Sigma,z_r}^{e_r}\cdot\ko_C(C)$.}}
     $\kj_{X(e_1;z_1)/\Sigma}\otimes\cdots\otimes\kj_{X(e_r;z_r)/\Sigma}\otimes\ko_C(C)\longrightarrow 
     \kj_{X(\underline{e};\underline{z})/\Sigma}\cdot\ko_C(C)$ 
     induces a surjective morphism $\beta$ on the level of global
     sections.\tom{\footnote{$\beta$ is surjective since 
     $\supp\big(\Ker(\beta)\big)\subseteq\{z_1,\ldots,z_r\}$ by Lemma
     \ref{lem:support:1}, and hence
     $H^1(\Sigma,\Ker(\beta))=0$.}} 
     \\
     If
     $s\in H^0\big(C,\kj_{X(e_1;z_1)/\Sigma}\otimes\tom{_{\ko_\Sigma}}\cdots
     \otimes\tom{_{\ko_\Sigma}}\kj_{X(e_r;z_r)/\Sigma}\otimes\tom{_{\ko_\Sigma}}\ko_C(C)\big)$, 
     but not in $\Ker(\beta)$, 
     then $s$ induces a non-zero section $\tilde{s}$ in
     $H^0\big(C,\pi_*\ko_{\widetilde{C}}\big(-\sum_{i=1}^re_iE_i\big)\otimes_{\ko_C}\ko_C(C)\big)$.
   \end{lemma} 

   \begin{proof}
     Set $E:=-\sum_{i=1}^re_iE_i$.\\     
     We start with the structure sequence for $\widetilde{C}$:
     \begin{displaymath}
       \xymatrix@C0.6cm{
         0\ar[r] & {\ko_{\widetilde{\Sigma}}\big(-\widetilde{C}\big)}\ar[r] &
         {\ko_{\widetilde{\Sigma}}}\ar[r] & {\ko_{\widetilde{C}}}\ar[r] & 0.
         }
     \end{displaymath}
     Tensoring with the locally free sheaf $\ko_{\widetilde{\Sigma}}(E)$ and then
     applying $\pi_*$ we get a morphism:
     \begin{displaymath}
       \pi_*\ko_{\widetilde{\Sigma}}(E) \longrightarrow\pi_*\ko_{\widetilde{C}}(E). 
     \end{displaymath}
     Now tensoring by $\ko_C(C)$ over $\ko_\Sigma$ we have an exact sequence:
     \begin{displaymath}
       \xymatrix@C0.6cm{
         0\ar[r] & {\Ker(\gamma)}
         \ar[r] & {\pi_*\ko_{\widetilde{\Sigma}}(E)\otimes_{\ko_\Sigma}\ko_C(C)}
         \ar[r]^\gamma & 
         {\pi_*\ko_{\widetilde{C}}(E)\otimes_{\ko_\Sigma}\ko_C(C)}. 
         }
     \end{displaymath}
     And finally taking global sections, we end up with:
     \begin{displaymath}
       \xymatrix@C0.55cm{
         0\ar[r] &  H^0\big(\Sigma,\Ker(\gamma)\big)
         \ar[r] & H^0\big(\Sigma,\pi_*\ko_{\widetilde{\Sigma}}(E)\otimes\ko_C(C)\big)
         \ar[r]^{\alpha} & 
         H^0\big(\Sigma,\pi_*\ko_{\widetilde{C}}(E)\otimes\ko_C(C)\big).
         }
     \end{displaymath}
     Since the sheaves we look at are actually $\ko_C$-sheaves 
     and since $C$ is a closed subscheme of $\Sigma$, the global sections
     of the sheaves as sheaves on $\Sigma$ and as sheaves on $C$ coincide
     (cf.~\cite{Har77} III.2.10  - for more details, see Corollary
     \ref{cor:identifications}). Furthermore,
     $\pi_*\ko_{\widetilde{\Sigma}}(E)=\bigotimes_{i=1}^r\kj_{X(e_i;z_i)/\Sigma}$.

     Thus it suffices to show that
     $\Ker{(\alpha)}\subseteq\Ker{(\beta)}$.

     Since $\pi_|:\widetilde{\Sigma}\setminus(\bigcup_{i=1}^rE_i)\rightarrow
     \Sigma\setminus\{z_1,\ldots,z_r\}$ is an isomorphism, we have that
     $\supp\big(\Ker(\gamma)\big)\subseteq\{z_1,\ldots,z_r\}$ is
     finite.\footnote{See Lemma \ref{lem:support:2}.} 
     Hence, by Lemma
     \ref{lem:torsion}, 
     \begin{displaymath}
       \Ker(\alpha)\;\;=\;\;H^0\big(\Sigma,\Ker(\gamma)\big)\;\;\subseteq 
       \;\;H^0\Big(\Sigma,\Tor\big(\Ker(\gamma)\big)\Big)
     \end{displaymath}
     \begin{displaymath}
       \subseteq
       H^0\bigg(\Sigma,\Tor\Big(\bigotimes\nolimits_{i=1}^r\kj_{X(e_i;z_i)/\Sigma}\otimes\ko_C(C)\Big)\bigg).
     \end{displaymath}
     Let now $t\in\Ker(\alpha)$ be given. We have to show that
     $\beta(t)=0$, i.~e.~$\beta_z(t_z)=0$ for every $z\in \Sigma$. If
     $z\not\in\{z_1,\ldots,z_r\}$, then $t_z=0$. Thus
     we may assume $z=z_k$. As we have shown,
     \begin{displaymath}
       t_{z_k}\in\Tor\big(\m_{\Sigma,z_k}^{e_k}\otimes_{\ko_{\Sigma,z_k}}\ko_{C,z_k}\big)
       =\Tor\big(\m_{\Sigma,z_k}^{e_k}/f_{z_k}\m_{\Sigma,z_k}^{e_k}\big)=(f_{z_k})/f_{z_k}\m_{\Sigma,z_k}^{e_k},
     \end{displaymath}
     where $f_{z_k}$ is a local 
     equation of $C$ at $z_k$. Therefore, there exists a
     $0\not=g_{z_k}\in\ko_{\Sigma,z_k}$ such that $t_{z_k}=f_{z_k} g_{z_k}\; \big(\mod
     f_{z_k}\m_{\Sigma,z_k}^{e_k}\big)\equiv f_{z_k}\otimes g_{z_k}$
     (note that $f_{z_k}\in\m_{\Sigma,z_k}^{n_k}\subseteq\m_{\Sigma,z_k}^{e_k}$!).  
     But then $\beta_{z_k}(t_{z_k})$ is just the residue 
     class of $f_{z_k} g_{z_k}$ in
     $\m_{\Sigma,z_k}^{e_k}\ko_{C,z_k}=\m_{\Sigma,z_k}^{e_k}/(f_{z_k})$, and is thus zero.
   \end{proof}

   \begin{proof}[Proof of Lemma \ref{lem:Geng-Xu}]
     Using the notation of the idea of the proof given on page \pageref{idea:Geng-Xu},
     we have, by Lemma \ref{lem:deformation}, a non-zero section
     $s\in
     H^0\big(C,\kj_{X(\underline{e};\underline{z})/\Sigma}\cdot\ko_C(C)\big)$. 
     This lifts under the surjection $\beta$ to a section $s'\in
     H^0\big(C,\bigotimes_{i=1}^r\kj_{X(e_i;z_i)/\Sigma}\otimes\ko_C(C)\big)$ which is not in
     the kernel of $\beta$. 
     Again setting $E:=-\sum_{i=1}^re_iE_i$, by Lemma \ref{lem:lifting}, we have a non-zero section
     $\tilde{s}\in
     H^0\big(C,\pi_*\ko_{\widetilde{C}}(E)\otimes_{\ko_C}\ko_C(C)\big)$, where by
     the projection formula the latter is just
     $H^0\big(C,\pi_*\big(\ko_{\widetilde{C}}(E)\otimes_{\ko_{\widetilde{C}}}\pi^*\ko_C(C)\big)\big)
     =_{\mbox{\footnotesize def}}
     H^0\big(\widetilde{C},\ko_{\widetilde{C}}(E)\otimes_{\ko_{\widetilde{C}}}\pi^*\ko_C(C)\big)$.\\
     Since $\ko_{\widetilde{C}}(E)\otimes_{\ko_{\widetilde{C}}}\pi^*\ko_C(C)$ has 
     a global section and since $\widetilde{C}$ is
     irreducible and reduced, we get by Lemma \ref{lem:degree}:
     $$0\leq\deg\big(\ko_{\widetilde{C}}(E)\otimes_{\ko_{\widetilde{C}}}\pi^*\ko_C(C)\big)
     =E.\widetilde{C}+\deg\big(\ko_C(C)\big)
     =\sum_{i=1}^r -e_i n_i + C^2.$$
   \end{proof}


   \section{Existence Theorem for Generic Fat Point Schemes}\label{sec:existence-I}
   \setcounter{equation}{0}

   \begin{theorem}\label{thm:existence-I}
     Given $m_1,\ldots,m_r\in\N_0$, not all zero, and $z_1,\ldots,z_r\in\Sigma$,
     $r\geq 1$, in very general position.
     Let $L\in\Div(\Sigma)$ be very ample
     over $\C$, and 
     let $D\in\Div(\Sigma)$ be such that
     \begin{equationlist}
        \item[eq:existence-I:1]
          $h^1\big(\Sigma,\kj_{X(\underline{m};\underline{z})/\Sigma}(D-L)\big)=0$, and
        \item[eq:existence-I:2]
          $D.L-2 g(L)\geq m_i+m_j$ for all $i,j$.
        \theequationchangeback
     \end{equationlist}
     Then there exists a curve $C\in |D|_l$ with 
     ordinary singular
     points of multiplicity $m_i$ at $z_i$ for $i=1,\ldots,r$ and no
     other singular points. Furthermore,
     \begin{displaymath}
       h^1\big(\Sigma, \kj_{X(\underline{m};\underline{z})/\Sigma}(D)\big)=0,
     \end{displaymath}
     and in particular, $V_{|D|}(\underline{m})$ is T-smooth at
     $C$.

     If in addition
     \stepcounter{equation}\mbox{\rm(\theequation)} \mylabel{eq:existence-I:3}{\theequation}
     $D^2>\sum_{i=1}^r m_i^2$,
     then $C$ can be chosen to be irreducible and reduced.
   \end{theorem}

   \begin{varthm-roman}[Idea of the proof]
     For each $z_j$ find a curve
     $C_j\in\big|H^0\big(\kj_{X(\underline{m},\underline{z})/\Sigma}(D)\big)\big|_l$ 
     with an ordinary singular point of multiplicity $m_j$ and show that
     this linear system has no other base points than
     $z_1,\ldots,z_r$. Then the generic element is smooth outside
     $z_1,\ldots,z_r$ and has an ordinary singularity of
     multiplicity $m_j$ in $z_j$.
   \end{varthm-roman}

   \begin{proof}
     W.~l.~o.~g.~$m_i\geq 1$ for all $i=1\ldots,r$. For the
     convenience of notation we set $z_{r+1}:=z_1$ and 
     $m_{r+1}:=m_1$.
     Since $L$ is very ample, we
     may choose smooth curves 
     $L_j\in|L|_l$ through $z_j$ and $z_{j+1}$ for $j=1,\ldots,r$
     (cf.~Lemma \ref{lem:smoothcurves}).
     Writing $X$ for
     $X(\underline{m};\underline{z})$ we introduce
     zero-dimensional schemes $X_j$ for $j=1,\ldots,r$ by
     \begin{displaymath}
       \kj_{X_j/\Sigma,z}=\left\{
         \begin{array}{ll}
           \kj_{X/\Sigma,z}, & \mbox{ if }z\not=z_j,\\
           \m_{\Sigma,z_j}\cdot\kj_{X/\Sigma,z_j}, & \mbox{ if }z=z_j.
         \end{array}
         \right. 
     \end{displaymath}

     \begin{varthm-roman}[Step 1]
       $h^1\big(\kj_{X_j/\Sigma}(D)\big)=0$. 
     \end{varthm-roman}
     By Condition (\ref{eq:existence-I:2}) 
     we get
     \begin{equation}
       \label{eq:existence-I:4}
       \deg\big(X_j\cap L_j\big)\leq m_j+m_{j+1} +1 \leq D.L+1-2 g(L),
     \end{equation}
     and the exact sequence 
     \begin{displaymath}
       \xymatrix@C0.6cm{
         0\ar[r] & {\kj_{X/\Sigma}}\ar[r] &{\kj_{X_j:L_j/\Sigma}}\ar[r] &
         {\m_{\Sigma,z_{j+1}}^{m_{j+1}-1}/\m_{\Sigma,z_{j+1}}^{m_{j+1}}}\ar[r] & 0
         }
     \end{displaymath}
     implies with the aid of (\ref{eq:existence-I:1})
     \begin{equation}
       \label{eq:existence-I:5}
       h^1\big(\kj_{X_j:L_j/\Sigma}(D-L)\big)=0.
     \end{equation}
     (\ref{eq:existence-I:4}) and (\ref{eq:existence-I:5}) allow us
     to apply Lemma \ref{lem:inductionstep} in order to obtain
     \begin{displaymath}
       h^1\big(\kj_{X_j/\Sigma}(D)\big)=0.
     \end{displaymath}

     \begin{varthm-roman}[Step 2]
       For each $j=1,\ldots,r$ there exists a curve 
       $C_j\in|D|_l$ with an  
       ordinary singular point of multiplicity $m_j$ at $z_j$ and with 
       $\mult_{z_i}(C_j)\geq m_i$ for $i\not=j$. 
     \end{varthm-roman}
     Consider the exact sequence
     \begin{displaymath}
       \xymatrix@C0.6cm{
         0\ar[r] & {\kj_{X_j/\Sigma}}\ar[r] &{\kj_{X/\Sigma}}\ar[r] &
         {\m_{\Sigma,z_j}^{m_j}/\m_{\Sigma,z_j}^{m_j+1}}\ar[r] & 0
         }       
     \end{displaymath}
     twisted by $D$ and the corresponding long exact cohomology sequence
     \begin{equation}\label{eq:existence-I:6}
       \xymatrix@C0.4cm@R0.3cm{
         H^0\big(\kj_{X/\Sigma}(D)\big)\ar[r] &
         {\m_{\Sigma,z_j}^{m_j}/\m_{\Sigma,z_j}^{m_j+1}}\ar[r] &
         H^1\big(\kj_{X_j/\Sigma}(D)\big)
         \ar[r]\ar@{=}[d]^{\mbox{\tiny Step 1}} & 
         H^1\big(\kj_{X/\Sigma}(D)\big)\ar[r]&0.\\
         & & 0 & &
         }
     \end{equation}
     Thus we may choose the $C_j$ to be given by a section in 
     $H^0\big(\kj_{X/\Sigma}(D)\big)$
     where the $m_j$ tangent directions at $z_j$ are all different. 

     \begin{varthm-roman}[Step 3]
       The base locus of $\P\big(H^0\big(\kj_{X/\Sigma}(D)\big)\big)$ is
       $\{z_1,\ldots,z_r\}$. 
     \end{varthm-roman}
     Suppose $w\in\Sigma$ was an additional base point and define the 
     zero-dimensional scheme $X\cup\{w\}$ by 
     \begin{displaymath}
       \kj_{X\cup\{w\}/\Sigma,z}=\left\{
         \begin{array}{ll}
           \kj_{X/\Sigma,z}, & \mbox{ if }z\not=w,\\
           \m_{\Sigma,w}\cdot\kj_{X/\Sigma,w}, & \mbox{ if }z=w.
         \end{array}
         \right.        
     \end{displaymath}
     Choosing a generic, and thus smooth\tom{\footnote{Consider the linear subsystem of 
         $|L|_l$ of curves passing through $w$. Then $\{w\}$ is its
         base locus (since the curves are given as hyperplane sections 
         of the hyperplanes passing through $w$) and the generic
         element is thus smooth outside $w$ by Bertini's Theorem. But since $L$
         is very ample over $\C$ it defines a closed embedding, and
         thus the elements passing through $w$ separate tangent
         vectors, i.~e.~their germs in $w$ generate
         $\m_{\Sigma,w}/\m_{\Sigma,w}^2$. In particular, the generic
         element is smooth in $w$. - One could also see this in a very 
         elementary way. Since $\Sigma$ is smooth in $w$ it is locally 
         in $w$ a complete intersection embedded in $\PC^n$ via $L$, 
         and its local analytic
         ring in is w.~l.~o.~g.~of the form
         $\C\{x_1,\ldots,x_N\}/(x_1+h.o.t.,\ldots,x_{N-2}+h.o.t.)$,
         where the $x_1$ are equations of hyperplanes in
         $\PC^n$ through $w$. Thus a generic hyperplane
         $\sum_{i=1}^N a_ix_i$ (with $(a_{N-1},a_N)\not=(0,0)$) will
         give rise to a smooth curve through $w$ since its local
         analytic ring is regular.}}, 
     curve $L_w\in|L|_l$ through $w$ we may deduce
     as in Step 1 
     \begin{displaymath}
       h^1\big(\kj_{X\cup\{w\}/\Sigma}(D)\big)=0,
     \end{displaymath}
     and thus as in Step 2
     \begin{displaymath}
       h^0\big(\kj_{X/\Sigma}(D)\big)=h^0\big(\kj_{X\cup\{w\}/\Sigma}(D)\big)+1.
     \end{displaymath}
     But by assumption $w$ is a base point, and thus
     \begin{displaymath}
       h^0\big(\kj_{X/\Sigma}(D)\big)=h^0\big(\kj_{X\cup\{w\}/\Sigma}(D)\big), 
     \end{displaymath}
     which gives us the desired contradiction.

     \begin{varthm-roman}[Step 4]
       $\exists\;C\in
       \P\big(H^0\big(\kj_{X/\Sigma}(D)\big)\big)\subseteq |D|_l$ with  
       an ordinary singular
       point of multiplicity $m_i$ at $z_i$ for $i=1,\ldots,r$ and no
       other singular points. 
     \end{varthm-roman}
     Because of Step 2 the generic element in
     $\P\big(H^0\big(\kj_{X/\Sigma}(D)\big)\big)$ has an ordinary singular 
     point of multiplicity $m_i$ at $z_i$ and is by Bertini's Theorem
     (cf.~\cite{Har77} III.10.9.2) smooth outside its base locus.

     For two generic curves $C,
     C'\in\P\big(H^0\big(\kj_{X/\Sigma}(D)\big)\big)$ the
     intersection multiplicity in $z_i$ is $i(C,C';z_i)=m_i^2$.
     Thus, if Condition (\ref{eq:existence-I:3}) 
     is fulfilled then 
     $C$ and $C'$ have an additional intersection point outside
     the base locus of
     $\P\big(H^0\big(\kj_{X/\Sigma}(D)\big)\big)$,\tom{\footnote{
     Just consider the inequality
     \begin{displaymath}
       \sum_{z\in C\cap C'}i(C,C';z)=C.C'=D^2>\sum_{i=1}^r
       m_i^2=\sum_{i=1}^r i(C,C';z_i). 
     \end{displaymath}}}
     and
     Bertini's Theorem (cf.~\cite{Wae73} \S 47, Satz 4) implies
     that the generic curve in
     $\P\big(H^0\big(\kj_{X/\Sigma}(D)\big)\big)$ is irreducible.

     \begin{varthm-roman}[Step 5]
       $h^1\big(\kj_{X/\Sigma}(D)\big)=0$, by Equation (\ref{eq:existence-I:6}).
     \end{varthm-roman}

     \begin{varthm-roman}[Step 6]
       $V_{|D|}(\underline{m})$ is T-smooth at $C$.
     \end{varthm-roman}
     By \cite{GLS98a}, Lemma 2.7, we have 
     \begin{displaymath}
       \kj_{X/\Sigma}\subseteq \kj_{X^{es}(C)/\Sigma},
     \end{displaymath}
     and thus by Step 5
     \begin{displaymath}
       h^1\big(\kj_{X^{es}(C)/\Sigma}(D)\big)=0,
     \end{displaymath}
     which proves the claim.
   \end{proof}

   \begin{corollary}\label{cor:existence-II}
     Let $m_1,\ldots,m_r\in\N_0$, not all zero, $r\geq 1$, and
     let $L\in\Div(\Sigma)$ be very ample over $\C$. 
     Suppose $D\in\Div(\Sigma)$ such that
     \begin{equationlist}
        \item[eq:existence-II:1]
          $(D-L-K_{\Sigma})^2 \geq 2\sum_{i=1}^r (m_i+1)^2$,
        \item[eq:existence-II:2]
          $(D-L-K_\Sigma).B > \max\{m_1,\ldots,m_r\}\;$ for any
          irreducible curve $B\subset \Sigma$  with $B^2=0$
          and $\dim|B|_a\geq 1$, 
        \item[eq:existence-II:3]
          $D-L-K_\Sigma$ is nef, and
        \item[eq:existence-II:4]
          $D.L-2 g(L)\geq m_i+m_j$ for all $i,j$.
     \end{equationlist}
     Then for $z_1,\ldots,z_r\in\Sigma$ in very general position there exists a
     curve $C\in |D|_l$ with  
     ordinary singular
     points of multiplicity $m_i$ at $z_i$ for $i=1,\ldots,r$ and no
     other singular points. Furthermore,
     \begin{displaymath}
       h^1\big(\Sigma, \kj_{X(\underline{m};\underline{z})/\Sigma}(D)\big)=0,
     \end{displaymath}
     and in particular, $V_{|D|}(\underline{m})$ is T-smooth in
     $C$.

     If in addition
     \stepcounter{equation}\mbox{\rm(\theequation)} \mylabel{eq:existence-II:5}{\theequation}
     $D^2>\sum_{i=1}^r m_i^2$,
     then $C$ can be chosen to be irreducible and reduced.
   \end{corollary}

   \begin{proof}
     Follows from Theorem \ref{thm:existence-I} and Corollary \ref{cor:vanishing}.
   \end{proof}

   \begin{remark}\label{rem:existence-I}
     In view of Condition (\ref{eq:existence-II:1}) Condition
     (\ref{eq:existence-II:5}) is satisfied if the
     following condition is fulfilled:
     \begin{equation}\label{eq:existence-II:6}
       D^2+(2D-L-K_{\Sigma}).(L+K_{\Sigma})+4\sum_{i=1}^rm_i +2r >0
     \end{equation}
   \end{remark}
   \begin{proof}
     Suppose (\ref{eq:existence-II:5}) was not satisfied, then
     \begin{displaymath}
       2\sum_{i=1}^r m_i^2\geq 2D^2 = D^2+(D-L-K_\Sigma)^2 + (2 D 
       -L-K_{\Sigma}).(L+K_\Sigma)
     \end{displaymath}
     \begin{displaymath}
       \geq D^2+ 2\sum_{i=1}^r m_i^2 +4 \sum_{i=1}^r m_i +2r +(2 D 
       -L-K_{\Sigma}).(L+K_\Sigma).
     \end{displaymath}
     Hence,
     \begin{displaymath}
       D^2+(2D-L-K_{\Sigma}).(L+K_\Sigma)+ 4
       \sum_{i=1}^r m_i +2r \leq 0, 
     \end{displaymath}
     which implies that (\ref{eq:existence-II:6}) is sufficient. 
   \end{proof}


   \section{Existence Theorem for General Equisingularity Schemes}\label{sec:existence-III}
   \setcounter{equation}{0}

   \begin{varthm-roman-break}[Notation]
     In the following we will denote by $\C[x,y]_d$,
     respectively by $\C[x,y]_{\leq d}$ the
     $\C$-vector spaces of polynomials of degree $d$,
     respectively of degree at most $d$. If $f\in\C[x,y]_{\leq d}$ we denote
     by $f_k\in\C[x,y]_k$ for $k=0,\ldots,d$ the homogeneous part of degree $k$ of 
     $f$, and thus $f=\sum_{k=0}^d f_k$. By $\ua=(a_{k,l}|0\leq
     k+l\leq d)$ we will denote the coordinates of $\C[x,y]_{\leq d}$
     with respect to the basis $\big\{x^ky^l|0\leq k+l\leq d\big\}$.

     For any $f\in\C[x,y]_{\leq d}$ the tautological family 
     \begin{displaymath}
       \C[x,y]_{\leq d}\times \C^2 \;\supset \;\bigcup_{g\in\C[x,y]_{\leq d}}\{g\}\times
         g^{-1}(0) \longrightarrow \C[x,y]_{\leq d}
     \end{displaymath}
     induces a deformation of the plane curve singularity
     $\big(f^{-1}(0),0\big)$ whose base space is the germ $\big(\C[x,y]_{\leq 
       d},f\big)$ of $\C[x,y]_{\leq d}$ at $f$.
     Given any deformation $(X,x)\hookrightarrow (\kx,x)\rightarrow (S,s)$ of a plane
     curve singularity $(X,x)$, we will denote by $S^{es}=(S^{es},s)$
     the germ of the equisingular stratum of $(S,s)$. 
     \tom{\footnote{That is, $S^{es}$ is the analytical space germ parametrising the
     subfamily of $(\kx,x)\rightarrow (S,s)$ of singularities which
     are topologically equivalent to $(X,x)$.}} 
     Thus, fixed an $f\in\C[x,y]_{\leq d}$, $\C[x,y]_{\leq
       d}^{es}=\big(\C[x,y]_{\leq d}^{es},f\big)$ is the (local) equisingular
     stratum of $\C[x,y]_{\leq d}$ at $f$. 
   \end{varthm-roman-break}

   \begin{edefinition}\label{def:Tsmooth}\leererpunkt
     \begin{enumerate}
        \item We say the family $\C[x,y]_{\leq d}$ is
          \emph{T-smooth} at $f\in\C[x,y]_{\leq d}$ if for any $e\geq
          d$ there exists a $\Lambda\subset\big\{(k,l)\in\N_0^2\;\big|\; 0\leq
          k+l\leq d\big\}$ with $\#\Lambda=\tau^{es}$ such that
          $\C[x,y]_{\leq e}^{es}$ is given by equations
          \begin{displaymath}
            a_{k,l}=\phi_{k,l}\big(\uaii,\uaiii\big),\;\;\; (k,l)\in\Lambda,
          \end{displaymath}
          with $\phi_{k,l}\in\C\big\{\uaii,\uaiii\big\}$ where
          $\uai=(a_{k,l}|(k,l)\in\Lambda)$, $\uaii=(a_{k,l}\;|\;0\leq
          k+l\leq d, (k,l)\not\in\Lambda)$, and $\uaiii=(a_{k,l}\;|\;
          d+1\leq k+l\leq e)$, and where
          $\tau^{es}=\dim_{\C}\big(\C\{x,y\}/I^{es}\big(f^{-1}(0),0\big)\big)$ is the
          codimension of the equisingular stratum in the base space of 
          the semiuniversal deformation of $\big(f^{-1}(0),0\big)$.
        \item A polynomial $f\in\C[x,y]_{\leq d}$ is said to be
          a \emph{good representative} of the singularity type $\ks$
          in $\C[x,y]_{\leq d}$
          if it meets the following conditions:
          \begin{enumerate}
             \item[(a)] $\Sing\big(f^{-1}(0)\big)=\left\{p\in\C^2\;|\;f(p)=0,\frac{\partial
                   f}{\partial x}(p)=0,\frac{\partial
                   f}{\partial y}(p)=0\right\}=\{0\}$,
             \item[(b)] $\big(f^{-1}(0),0\big)\sim_t\ks$,
             \item[(c)] $f_d$ is reduced, and
             \item[(d)] $\C[x,y]_{\leq d}$ is T-smooth at $f$.
          \end{enumerate}
        \item Given a singularity type $\ks$ we define $s(\ks)$ to be
          the minimal number $d$ such that $\ks$ has a good representative
          of degree $d$.
     \end{enumerate}
   \end{edefinition}

   \begin{eremark}\leererpunkt 
     \begin{enumerate}
        \item The condition for T-smoothness just means that for any
          $e\geq d$ the equisingular stratum $\C[x,y]_{\leq e}^{es}$
          is smooth at the point $f$ of the expected codimension in $\big(\C[x,y]_{\leq
            e},f\big)$.
        \item Note that for a polynomial of degree $d$ the
          highest homogeneous part $f_d$ defines
          the normal cone, i.~e.~the intersection of the curve $\{\hat{f}=0\}$
          with the line at infinity in $\PC^2$, where $\hat{f}$ is
          the homogenisation of $f$. Thus the condition
          ``$f_d$ reduced'' in the definition of a good representative 
          just means that the line at infinity intersects the curve
          transversally in $d$ different points.
        \item If $f\in\C[x,y,z]_d$ is an irreducible polynomial such
          that $(0:0:1)$ is the only singular point of the plane curve 
          $\{f=0\}\subset \PC^2$, then a linear change of
          coordinates of the type $(x,y,z)\mapsto (x,y,z+ax+by)$ will
          ensure that the dehomogenisation $\check{f}$ of $f$
          satisfies ``$\check{f}_d$ reduced''. Note for this that the
          coordinate change corresponds to choosing a line in
          $\PC^2$, not passing through $(0:0:1)$ and meeting the
          curve in $d$ distinct points. 
          Therefore, the bounds for $s(\ks)$ given in \cite{Los98}
          Theorem 4.2 and Remark 4.3 do apply here.
        \item For refined results using the techniques of the
          following proof we refer to \cite{Shu99}.
     \end{enumerate}
   \end{eremark}

   \begin{theorem}[Existence]\label{thm:existence-III}
     Let $\ks_1,\ldots,\ks_r$ be singularity types, and suppose there
     exists an irreducible curve $C\in|D|_l$ with $r+r'$ ordinary
     singular points $z_1,\ldots,z_{r+r'}$ of multiplicities
     $m_1,\ldots,m_{r+r'}$ respectively as its 
     only singularities such that $m_i=s(\ks_i)+1$, for $i=1,\ldots,r$,
     and 
     \begin{displaymath}
       h^1\big(\Sigma,\kj_{X(\underline{m};\underline{z})/\Sigma}(D)\big)=0.
     \end{displaymath}
     Then there exists an irreducible curve $\widetilde{C}\in|D|_l$ with $r$
     singular points of types $\ks_1,\ldots,\ks_r$  and
     $r'$ ordinary singular points of multiplicities
     $m_{r+1},\ldots, m_{r+r'}$ as its only
     singularities.\footnote{Here, of course,
       $\underline{m}=(m_1,\ldots,m_{r+r'})$ and $\underline{z}=(z_1,\ldots,z_{r+r'})$.}
   \end{theorem}

   \begin{varthm-roman}[Idea of the proof]
     The basic idea is to glue locally at the $z_i$ equations of good
     representatives for the $\ks_i$ into the curve $C$. Let us now
     explain more detailed what we mean by this.

     If $g_i=\sum_{k+l=0}^{m_i-1}a_{k,l}^\ifix x_i^ky_i^l$, $i=1,\ldots,r$, are
     good representatives of the $\ks_i$, then we are
     looking for a family $F_t$, $t\in(\C,0)$, in
     $H^0\big(\Sigma,\ko_\Sigma(D)\big)$ which in local coordinates
     $x_i, y_i$ at $z_i$ looks like
     \begin{displaymath}
       F_t^\i= \sum_{k+l=0}^{m_i-1}t^{m_i-1-k-l}\ta(t) x_i^ky_i^l + h.o.t.,
     \end{displaymath}
     where the $\ta(t)$ should be convergent power series in $t$ with
     $\ta(0)= a_{k,l}^\ifix$. Replacing $g_i$ by some arbitrarily
     small multiple $\lambda_i g_i$ the curve defined by $F_0$ is
     an arbitrarily small deformation of $C$ inside some suitable
     linear system, thus it is smooth outside $z_1,\ldots,z_{r+r'}$
     and has ordinary singular points in $z_1,\ldots,z_{r+r'}$. For 
     $t\not=0$, on the other hand, $F_t^\i$ can be transformed, by
     $(x_i,y_i)\mapsto (tx_i,ty_i)$, into a
     member of some family
     \begin{displaymath}
       \tilde{F}_t^\i= \sum_{k+l=0}^{m_i-1}\ta(t) x_i^ky_i^l +
       h.o.t.,\;\; t\in\C,
     \end{displaymath}
     with
     \begin{displaymath}
       \tilde{F}_0^\i=g_i.
     \end{displaymath}
     Using now the T-smoothness property of $g_i$, $i=1,\ldots,r$,
     we can choose the $\ta(t)$ such that this family is
     equisingular. Hence, for small $t\not=0$, the curve given by
     $F_t$ will have the 
     right singularities at the $z_i$. Finally, the
     knowledge on the singularities of the curve defined by $F_0$ and
     the conservation of Milnor numbers will ensure that the
     curve given by $F_t$ has no further singularities, for $t\not=0$
     sufficiently small.

     The proof will be done in several steps. First of all we are
     going to fix some notation by
     choosing a basis of $H^0\big(\Sigma,\ko_\Sigma(D)\big)$
     which reflects the ``independence'' of the coordinates at the
     different $z_i$ ensured by 
     $h^1\big(\Sigma,\kj_{X(\underline{m};\underline{z})/\Sigma}(D)\big)=0$ (Step 1.1),
     and by choosing good representatives for the $\ks_i$ (Step 1.2). In 
     a second step we are making an ``Ansatz'' for the family $F_t$,
     and, for the local investigation of the singularity type, we are
     switching to some other families $\tilde{F}_t^\i$,
     $i=1,\ldots,r$ (Step 2.1). We, then, reduce the problem of $F_t$, 
     for $t\not=0$ small, having the right singularities to a question 
     about the equisingular strata of some families of polynomials
     (Step 2.2),
     which in Step 2.3 will be solved. The final step serves to
     show that the curves $F_t$ have only the singularities which we
     controlled in the previous steps.
   \end{varthm-roman}

   \begin{proof}
     \begin{varthm-roman}[Step 1.1]
       Parametrise $|D|_l=\P\big(H^0\big(\ko_\Sigma(D)\big)\big)$.
     \end{varthm-roman}
     Consider the following exact sequence:
     \begin{displaymath}
       0\longrightarrow
       \kj_{X(\underline{m};\underline{z})/\Sigma}(D) \longrightarrow 
       \ko_\Sigma(D) \longrightarrow
       \bigoplus_{i=1}^{r+r'}\ko_{\Sigma,z_i}/\m_{\Sigma,z_i}^{m_i}
       \longrightarrow 0.
     \end{displaymath}
     Since
     $h^1\big(\kj_{X(\underline{m};\underline{z})/\Sigma}(D)\big)=0$,
     the long exact cohomology sequence gives
     \begin{displaymath}
       H^0\big(\ko_\Sigma(D)\big) =
       \bigoplus_{i=1}^{r+r'} \C\{x_i,y_i\}/(x_i,y_i)^{m_i}\oplus
       H^0\big(\kj_{X(\underline{m};\underline{z})/\Sigma}(D)\big), 
     \end{displaymath}
     where $x_i,y_i$ are local coordinates of $(\Sigma,z_i)$.

     We, therefore, can find a basis 
     $\big\{s_{k,l}^\i, s_j \;\big|\; 
     j=1,\ldots,e,\; i= 1,\ldots, r+r',\; 0\leq k+l\leq m_i-1\big\}$ 
     of $H^0\big(\ko_\Sigma(D)\big)$, with
     $e=h^0\big(\kj_{X(\underline{m};\underline{z})/\Sigma}(D)\big)$,
     such that\footnote{Throughout this proof we will use the 
       multi index notation 
       $\alpha=(\alpha_1,\alpha_2)\in\N^2$ and
       $|\alpha|=\alpha_1+\alpha_2$.} 
     \begin{itemize}
        \item $C$ is the curve defined by $s_1$,
        \item $(s_j)_{z_i}=\sum_{|\alpha|\geq m_i} B_\alpha^{j,i}
          x_i^{\alpha_1}y_i^{\alpha_2}$ for $j=1,\ldots,e$,
          $i=1,\ldots,r+r'$,
        \item $\big(s_{k,l}^\j\big)_{z_i}=\left\{
            \begin{array}{ll}
              x_i^ky_i^l+\sum\limits_{|\alpha|\geq m_i}
              A_{\alpha,k,l}^{\i,\i} x_i^{\alpha_1}y_i^{\alpha_2}, & \mbox{
                if } i=j, \\
              \sum\limits_{|\alpha|\geq m_i}
              A_{\alpha,k,l}^{\j,\i} x_i^{\alpha_1}y_i^{\alpha_2}, & \mbox{
                if } i\not=j.
            \end{array}
            \right.$
     \end{itemize}
     Let us now denote the coordinates of $H^0\big(\ko_\Sigma(D)\big)$ 
     w.~r.~t.~this basis by
     $(\ua,\underline{b})=\big(\ua^\ii{1},\ldots,\ua^\ii{r+r'},\underline{b}\big)$ with 
     $\ua^\i=\big(a_{k,l}^\i\;|\; 0\leq k+l\leq m_i-1\big)$ and
     $\underline{b}=(b_j \;|\;j=1,\ldots,e)$. 

     Thus the family 
     \begin{displaymath}
       F_{(\ua,\underline{b})}=\sum_{i=1}^{r+r'}\sum_{k+l=0}^{m_i-1}
       a_{k,l}^\i s_{k,l}^\i + \sum_{j=1}^e b_j s_j, \;\;\;
       (\ua,\underline{b})\in\C^N \mbox{ with } N=e+\sum_{i=1}^{r+r'}
       \binom{m_i+1}{2},
     \end{displaymath}
     parametrises $H^0\big(\ko_\Sigma(D)\big)$.

     \begin{varthm-roman}[Step 1.2]
       By the definition of $s(\ks_i)$ and since $s(\ks_i)=m_i-1$, we
       may choose good representatives  
       \begin{displaymath}
         g_i=\sum_{k+l=0}^{m_i-1} a_{k,l}^\ifix x_i^ky_i^l\in
         \C[x_i,y_i]_{\leq m_i-1}  
       \end{displaymath}
       for the $\ks_i$, $i=1,\ldots,r$. Let
       $\ua^\ifix=\big(a_{k,l}^\ifix\; \big| \; 0\leq k+l\leq m_i-1\big)$ and
       $\ua^{fix}=\big(\ua^{\ii{1},fix},\ldots,\ua^{\ii{r},fix}\big)$. 
       We should remark here that for any $\lambda_i\not=0$ the
       polynomial $\lambda_i g_i$ is also a good representative, and
       thus, replacing $g_i$ by $\lambda_i g_i$, we may assume that the 
       $a_{k,l}^\ifix$ are arbitrarily close to $0$.
     \end{varthm-roman}

     \begin{varthm-roman}[Step 2]
       We are going to glue the good representatives for the $\ks_i$
       into the curve $C$. More precisely, we are constructing a subfamily
       $F_t$, $t\in (\C,0)$, in
       $H^0\big(\ko_\Sigma(D)\big)$ such that, if $C_t\in
       |D|_l$ denotes the curve defined by $F_t$, 
       \begin{enumerate}
          \item[(1)] $z_1,\ldots,z_{r+r'}$ are the only singular points of
            the irreducible reduced curve $C_0$,
            and they are ordinary singularities of multiplicities $m_i-1$,
            for $i=1,\ldots,r$, and $m_i$, for $i=r+1,\ldots,r+r'$
            respectively,  
          \item[(2)] locally in $z_i$, $i=1,\ldots,r$, the
            $F_t$, for small $t\not=0$, can be transformed into members
            of a fixed $\ks_i$-equisingular family, 
          \item[(3)] while for $i=r+1,\ldots,r+r'$ and $t\not=0$ small
            $C_t$ has an ordinary singularity of multiplicity $m_i$ in 
            $z_i$.
       \end{enumerate}
     \end{varthm-roman}

     \begin{varthm-roman}[Step 2.1]
       ``Ansatz'' and first reduction for a local investigation.
     \end{varthm-roman}
     Let us make the following ``Ansatz'':
     \begin{displaymath}
       b_1=1, \;b_2=\ldots= b_e=0,\; \ua^\i =0,\; \mbox{ for }
       i=r+1,\ldots,r+r',
     \end{displaymath}
     \begin{displaymath}
       a_{k,l}^\i = t^{m_i-1-k-l}\cdot \ta,\;\mbox{ for } i=1,\ldots,r,\;
       0\leq k+l\leq m_i-1.
     \end{displaymath}
     This gives rise to a family 
     \begin{displaymath}
       F_{(t,\tua)} = s_1 + \sum_{i=1}^r \sum_{k+l=0}^{m_i-1}
       t^{m_i-1-k-l} \ta s_{k,l}^\i \in H^0\big(\ko_\Sigma(D)\big) 
     \end{displaymath}
     with $t\in\C$ and $\tua=\big(\tua^\ii{1},\ldots,\tua^\ii{r}\big)$ where
     $\tua^\i = \big(\ta\;|\; 0\leq k+l\leq m_i-1\big)\in\C^{N_i}$
     with $N_i=\binom{m_i+1}{2}$.

     Fixing $i\in \{1,\ldots,r\}$, in local coordinates at $z_i$ 
     the family looks like
     \begin{displaymath}
       F_{(t,\tua)}^\i:= \big(F_{(t,\tua)}\big)_{z_i} = 
       \sum_{k+l=0}^{m_i-1} t^{m_i-1-k-l} \ta x_i^ky_i^l +
       \sum_{|\alpha|\geq m_i} \varphi_\alpha^\i(t,\tua)\: x_i^{\alpha_1}y_i^{\alpha_2},
     \end{displaymath}
     with
     \begin{displaymath}
       \varphi_\alpha^\i(t,\tua) = B_\alpha^{1,i} +
         \sum_{j=1}^r\sum_{k+l=0}^{m_j-1}t^{m_j-1-k-l}\taj
         A_{\alpha,k,l}^{\j,\i}.
     \end{displaymath}

     For $t\not=0$ the transformation
     $\psi_t^i:(x_i,y_i)\mapsto (tx_i,ty_i)$ is indeed a coordinate
     transformation, and thus $F_{(t,\tua)}^\i$ is contact
     equivalent\footnote{Let $f,g\in\ko_n=\C\{x_1,\ldots,x_n\}$ be two
       convergent power series in $n$ indeterminates. We call $f$ and
       $g$ \emph{contact equivalent}, if $\ko_n/(f)\cong\ko_n/(g)$,
       and we write in this case $f\sim_c g$. Equivalently, we could
       ask the germs $\big(V(f),0\big)$ and $\big(V(g),0\big)$ to
       be isomorphic, that is, ask the
       singularities to be
       \emph{analytically equivalent}. C.~f.~\cite{DP00} Definition
       9.1.1 and Definition 3.4.19.} 
     to 
     \begin{displaymath}
       \tilde{F}_{(t,\tua)}^\i := t^{-m_i+1} \cdot
       F_{(t,\tua)}^\i(tx_i,ty_i) =
       \sum_{k+l=0}^{m_i-1} \ta x_i^ky_i^l + \sum_{|\alpha|\geq
         m_i}t^{1+|\alpha|-m_i}\varphi_\alpha^\i(t,\tua)\: x_i^{\alpha_1}y_i^{\alpha_2}.
     \end{displaymath}
     Note that for this new family in $\C\{x_i,y_i\}$ we have
     \begin{displaymath}
       \tilde{F}_{(0,\ua^{fix})}^\i = \sum_{k+l=0}^{m_i-1}
       a_{k,l}^\ifix x_i^ky_i^l = g_i,
     \end{displaymath}
     and hence it gives rise to a deformation of
     $\big(g_i^{-1}(0),0\big)$. 

     \begin{varthm-roman}[Step 2.2]
       Reduction to the investigation of the equisingular strata of certain
       families of polynomials.
     \end{varthm-roman}

     It is basically our aim to verify the $\tua$
     as convergent power series in $t$ such that the corresponding
     family is equisingular. However, since the
     $\tilde{F}_{(t,\tua)}^\i$ are power series in $x_i$ and $y_i$, we 
     cannot right away apply the T-smoothness property of $g_i$, but
     we rather have to reduce to polynomials. For this let $e_i$ be
     the determinacy bound\footnote{A power series
       $f\in\ko_n=\C\{x_1,\ldots,x_n\}$ (respectively the singularity
       $\big(V(f),0\big)$ defined by $f$) is said to be \emph{finitely
         determined} with respect to some equivalence relation $\sim$ if there exists
       some positive integer $e$ such that $f\sim g$ whenever $f$ and
       $g$ have the same $e$-jet. If $f$ is finitely determined,
       the smallest possible $e$ is called the \emph{determinacy
         bound}. Isolated singularities are finitely determined with
       respect to analytical equivalence and hence, for $n=2$, as well with respect to
       topological equivalence. C.~f.~\cite{DP00} Theorem 9.1.3 and
       Footnote \ref{foot:top-equiv}.} 
     of $\ks_i$ and define
     \begin{displaymath}
       \hat{F}_{(t,\tua)}^\i := 
       \sum_{k+l=0}^{m_i-1} \ta x_i^ky_i^l + \sum_{|\alpha| =
         m_i}^{e_i}t^{1+|\alpha|-m_i}\varphi_\alpha^\i(t,\tua)\:
       x_i^{\alpha_1}y_i^{\alpha_2}
       \equiv \tilde{F}_{(t,\tua)}^\i \big(\mod (x_i,y_i)^{e_i+1}\big).
     \end{displaymath}
     Thus $\hat{F}_{(t,\tua)}^\i$ 
     is a family in $\C[x_i,y_i]_{\leq
       e_i}$, and still
     \begin{displaymath}
        \hat{F}_{(0,\ua^{fix})}^\i = \tilde{F}_{(0,\ua^{fix})}^\i = g_i.
     \end{displaymath}

     We claim that it suffices to find $\tua(t)\in\C\{t\}$ with
     $\tua(0)=\big(a_{k,l}^\ifix\;\big|\;i=1,\ldots, r,\; 0\leq k+l\leq m_i-1\big)$,
     such that the families $\hat{F}_t^\i := \hat{F}_{(t,\tua(t))}^\i$,
     $t\in (\C,0)$, are in the equisingular strata $\C[x_i,y_i]_{\leq
       e_i}^{es}$, for $i=1,\ldots, r$. 

     Since then we have, for small\footnote{\label{foot:top-equiv} Here $f\sim_t g$, for
       two convergent power series $f,g\in\ko_2=\C\{x,y\}$, means that
       the singularities $\big(V(f),0\big)$ and $\big(V(g),0\big)$ are
       \emph{topologically equivalent}, that is, there exists a
       homeomorphism
       $\Phi:\big(\C^2,0\big)\rightarrow\big(\C^2,0\big)$ with
       $\Phi\big(V(f),0\big)=\big(V(g),0\big)$, which of course means,
       that this is correct for suitably chosen representatives. Note
       that if $f$ and $g$ are contact equivalent, then there exists
       even an analytic coordinate change $\Phi$, that is, $f\sim_c g$
       implies $f\sim_t g$.} 
     $t\not=0$,
     \begin{displaymath}
       g_i=\hat{F}_0^\i \sim \hat{F}_t^\i \sim \tilde{F}_{(t,\tua(t))}^\i 
       \sim F_{(t,\tua(t))}^\i = \big(F_{(t,\tua(t))}\big)_{z_i(t)},
     \end{displaymath}
     by the $e_i$-determinacy and since $\psi_t^i$ is a coordinate
     change for $t\not=0$, which proves condition (2). 
     Note that the singular points $z_i$ will move with $t$. 

     It remains to verify conditions (1) and (3). Setting $F_t:= F_{(t,\tua(t))}\in
     H^0\big(\ko_\Sigma(D)\big)$, for $t\in (\C,0)$, we find that
     \begin{displaymath}
       F_0=s_1+\sum_{j=1}^r \sum_{k+l=m_j-1} a_{k,l}^\jfix s_{k,l}^\j
     \end{displaymath}
     is an element inside the linear system $\kd=\{\lambda_0 s_1 +
     \sum_{j=1}^r \lambda_j s^\j \; |\;
     (\lambda_0:\ldots:\lambda_r)\in\PC^r\}$, where $s^\j=
     \sum_{k+l=m_j-1} a_{k,l}^\jfix s_{k,l}^\j$. 
     Locally at $z_i$, $i=1,\ldots,{r+r'}$, $\kd$ induces a
     deformation of $(C,z_i)$ with equations
     \begin{displaymath}
       \lambda_i\cdot \big(g_i\big)_{m_i-1} + h.o.t.,\;\;\mbox{ if } i=1,\ldots,r,
     \end{displaymath}
     and
     \begin{displaymath}
       \lambda_0\cdot\left(\sum_{|\alpha|= m_i}B_\alpha^{1,i}
       x_i^{\alpha_1}y_i^{\alpha_2}\right) +
       \sum_{j=1}^r \lambda_j \cdot \left(\sum_{k+l=m_j-1}a_{k,l}^\jfix 
       \sum_{|\alpha|=m_i}A_{\alpha,k,l}^{\j,\i}
       x_i^{\alpha_1}y_i^{\alpha_2}\right) + h.o.t.,
     \end{displaymath}
     \begin{displaymath}
       \mbox{ if } i=r+1,\ldots,r+r',
     \end{displaymath}
     respectively. Thus any element of $\kd$ has ordinary
     singularities of multiplicity $m_i-1$ at $z_i$ for $i=1,\ldots,r$, 
     and since $s_1$ has an ordinary
     singularity of multiplicity $m_i$ at $z_i$ for
     $i=r+1,\ldots,r+r'$, so has a generic element of $\kd$. 
     Moreover, a generic
     element of $\kd$ has not more singular points than the special
     element $s_1$ and has thus 
     singularities precisely in $\{z_1,\ldots,z_{r+r'}\}$.  
     Replacing the $g_i$ by some suitable multiples, we
     may assume that the curve defined by $F_0$ is a generic element
     of $\kd$, which proves (1). Similarly, we note that $F_t$ in
     local coordinates at $z_i$, for $i=r+1,\ldots,r+r'$, looks like
     \begin{displaymath}
       \sum_{|\alpha|= m_i}B_\alpha^{1,i} x_i^{\alpha_1}y_i^{\alpha_2}
       + 
       \sum_{j=1}^r \sum_{k+l=0}^{m_j-1}
       t^{m_j-1-k-l} \taj(t) 
       \sum_{|\alpha|=m_i}A_{\alpha,k,l}^{\j,\i}
       x_i^{\alpha_1}y_i^{\alpha_2}
       + h.o.t.
     \end{displaymath}
     \begin{displaymath}
       = \sum_{|\alpha|=m_i}\left(B_\alpha^{1,i} + \sum_{j=1}^r
         \sum_{k+l=0}^{m_j-1} t^{m_j-1-k-l} \taj(t)
         A_{\alpha,k,l}^{\j,\i}\right)\: x_i^{\alpha_1}y_i^{\alpha_2}
       + h.o.t.,
     \end{displaymath}
     and thus, for $t\not=0$ sufficiently small, the singularity of $F_t$ at
     $z_i$ will be an ordinary singularity of multiplicity
     $m_i$, which gives (3).

     \begin{varthm-roman}[Step 2.3]
       Find $\tua(t)\in\C\{t\}^n$ with 
       $\tua(0)=\big(a_{k,l}^\ifix$, $i=1,\ldots, r,\; 0\leq k+l\leq
       m_i-1\big)$, $n={\sum_{i=1}^r\binom{m_i+1}{2}}$,
       such that the families $\hat{F}_t^\i = \hat{F}_{(t,\tua(t))}^\i$,
       $t\in (\C,0)$, are in the equisingular strata $\C[x_i,y_i]_{\leq
         e_i}^{es}$, for $i=1,\ldots r$.
     \end{varthm-roman}
     In the sequel we adopt the notation of definition \ref{def:Tsmooth} adding
     indices $i$ in the obvious way.

     Since $\C[x_i,y_i]_{\leq m_i-1}$ is T-smooth at $g_i$, for
     $i=1,\ldots,r$, there exist $\Lambda_i\subseteq \{(k,l)\;|\; 0\leq 
       k+l\leq m_i-1\}$ and power series $\phi_{k,l}^\i\in
       \C\big\{\tuaii^\i,\tuaiii^\i\big\}$, for $(k,l)\in\Lambda_i$,
       such that the equisingular 
       stratum $\C[x_i,y_i]_{\leq e_i}^{es}$ is given by the
       $\tau^{es,\i}=\#\Lambda_i$ equations 
       \begin{displaymath}
         \ta = \phi_{k,l}^\i \big(\tuaii^\i,\tuaiii^\i\big), \;\;\mbox{
           for } (k,l)\in\Lambda_i.
       \end{displaymath}
       Setting $\Lambda=\bigcup_{j=1}^r\{j\}\times\Lambda_j$ we use
       the notation        
       $\tuai=\big(\tuai^\ii{1},\ldots,\tuai^\ii{r}\big) =\big(\ta \;
       \big|\; (i,k,l)\in\Lambda\big)$ and, similarly $\tuaii$,
       $\tuaiii$, $\uaii^\ifix$, $\uai^{fix}$, and $\uaii^{fix}$. 
       Moreover, setting $\tilde{\varphi}^\i\big(t,\tuai\big)=
       \big(t^{|\alpha|-m_i}\varphi_\alpha^\i\big(t,\tuai,\uaii^{fix}\big) \;\big|\;
       m_i\leq |\alpha| \leq e_i\big)$, 
       we define an analytic map germ 
       \begin{displaymath}
         \Phi :
         \Big(\C\times\C^{\tau^{es,\ii{1}}}\times\cdots\times\C^{\tau^{es,\ii{r}}},\big(0,\uai^{fix}\big)\Big)
         \rightarrow
         \Big(\C^{\tau^{es,\ii{1}}}\times\cdots\times\C^{\tau^{es,\ii{r}}},0\Big)
       \end{displaymath}
       by 
       \begin{displaymath}
         \Phi_{k,l}^\i \big(t,\tuai\big) = \ta - \phi_{k,l}^\i 
         \Big(\uaii^\ifix, t\cdot\tilde{\varphi}^\i\big(t,\tuai\big)\Big),\;\;\mbox{ for }
         (i,k,l)\in\Lambda,
       \end{displaymath}
       and we consider the system of equations
       \begin{displaymath}
         \Phi_{k,l}^\i \big(t,\tuai\big) = 0, \;\;\mbox{ for }
         (i,k,l)\in\Lambda.
       \end{displaymath}
       One easily verifies that
       \begin{displaymath}
         \left(\frac{\partial\Phi_{k,l}^\i}{\partial\tal}
         \big(0,\uai^{fix}\big)\right)_{(i,k,l),(j,\kappa,\lambda)\in\Lambda} 
         = \id_{\C^n}.
       \tom{\footnote{
       \begin{varthm-roman}[Step 2.4]
         Show that 
         $\left(\frac{\partial\Phi_{k,l}^\i}{\partial\tal} 
         \big(0,\uai^{fix}\big)\right)_{(i,k,l),(j,\kappa,\lambda)\in\Lambda} 
         = \id_n$. 
       \end{varthm-roman} 
       For this it suffices to show that 
       \begin{displaymath}
         \frac{\partial\phi_{k,l}^\i\big(\uaii^\ifix,t\cdot
         \tilde{\varphi}^\i\big(t,\tuai\big)\Big)}{\partial\tal}_{\Big| t=0,\tuai=\uai^{fix}}
          = 0,
       \end{displaymath} 
       for any  $(i,k,l),(j,\kappa,\lambda)\in\Lambda$, which is
       fulfilled since the map 
       \begin{displaymath}
         \theta: \big(t,\tuai\big)\mapsto \big(\uaii^{fix},t\cdot\tilde{\varphi}\big(t,\tuai\big)\Big)
       \end{displaymath}
       satisfies
       \begin{displaymath}
         \frac{\partial\theta}{\partial\tal}_{\big| t=0,\tuai=\uai^{fix}}
         = \left(\underline{0},\;
           t\cdot\frac{\partial\tilde{\varphi}^\i\big(t,\tuai\big)}{\partial\tal}_{\big| t=0,\tuai=\uai^{fix}}\right) 
         = 0.
       \end{displaymath}
       }}
       \end{displaymath}
       Thus by the Inverse Function Theorem there exist
       $\ta(t)\in\C\{t\}$ with $\ta(0)=a_{k,l}^\ifix$ such that 
       \begin{displaymath}
         \ta(t) = \phi_{k,l}^\i
         \Big(\uaii^\ifix,t\cdot\tilde{\varphi}^\i\big(t,\tuai(t)\big)\Big),
         \;\;\;(i,k,l)\in\Lambda.
       \end{displaymath}
       Now, setting $\tuaii(t)\equiv \uaii^{fix}$, the families
       $\hat{F}_t^\i=\hat{F}_{(t,\tua(t))}^\i$ are in the equisingular strata
       $\C[x_i,y_i]_{\leq e_i}^{es}$, for $i=1,\ldots,r$.

       \begin{varthm-roman}[Step 3]
         It finally remains to show that $F_t$, for
         small $t\not=0$, has  no other singular points than
         $z_1(t),\ldots,z_{r}(t),z_{r+1},\ldots,z_{r+r'}$. 
       \end{varthm-roman}
       Since for any $i=1,\ldots,r+r'$ the family $F_t$, $t\in(\C,0)$,
       induces a deformation of the 
       singularity $(C_0,z_i)$  there are, by the Conservation of
       Milnor Numbers\footnote{Recall the definition of the Milnor
         number of a holomorphic map $f\in\ko(U)$ respectively of 
         $f^{-1}\big(f(z)\big)$ at a point $z\in U\subset\C^2$:
         $\mu(f,z)=\mu\big(f^{-1}\big(f(z)\big),z\big)=
         \dim_\C\Big(\ko_{U,z} \big/ \big(\frac{\partial
           f}{\partial x}(z),\frac{\partial f}{\partial y}(z)\big)\Big).$
         } 
       (cf.~\cite{DP00}, Chapter 6), (Euclidean) open
       neighbourhoods $U(z_i)\subset \Sigma$ and $V(0)\subset\C$ such
       that for any $t\in V(0)$
       \begin{equationlist}
          \item[eq:existence-III:1] $\Sing(C_t) \subset \bigcup_{i=1}^{r+r'} 
            U(z_i)$, i.~e.~singular points of $C_t$ come from singular 
            points of $C_0$, 
          \item[eq:existence-III:2] $\mu(C_0,z_i)=\sum_{z\in\Sing(F_t^\i)\cap
              U(z_i)} \mu\big(F_t^\i,z\big)$,\;\; $i=1,\ldots,r+r'$.
       \end{equationlist}
       For $i=r+1,\ldots,r+r'$ condition (\ref{eq:existence-III:2}) implies
       \begin{displaymath}
         (m_i-1)^2=\mu(C_0,z_i)\geq\mu\big(F_t^i,z_i\big)=(m_i-1)^2,
       \end{displaymath}
       and thus $z_i$ must be the only critical point of $F_t^\i$ in 
       $U(z_i)$, in particular,
       \begin{displaymath}
         \Sing(C_t)\cap U(z_i) = \{z_i\}.
       \end{displaymath}

       Let now $i\in\{1,\ldots,r\}$. For $t\not=0$ fixed, we consider the transformation
       defined by the coordinate change $\psi_t^i$,
       \begin{displaymath}
         \begin{array}{rcl}
           \C^2\;\;\supset\;\; U(z_i)\;\; &\longrightarrow & \;\;U_t(z_i)\;\;\subset\;\;\C^2\\
           {\begin{sideways}$\in$\end{sideways}}\;\;\;\;\; & &\;\;\;\;\;{\begin{sideways}$\in$\end{sideways}}\\
           (x_i,y_i) & \mapsto & \big(\frac{1}{t}x_i,\frac{1}{t}y_i\big),
         \end{array}
       \end{displaymath}
       and the transformed equations
       \begin{displaymath}
         \tilde{F}_t^\i(x_i,y_i)= t^{-m_i+1}F_t^i(tx_i,ty_i)=0.
       \end{displaymath}
       Condition (\ref{eq:existence-III:2}) then implies,
       \begin{displaymath}
         (m_i-2)^2 = \mu(C_0,z_i)
         =\sum_{z\in\Sing(F_t^\i)\cap U(z_i)} \mu\big(F_t^\i,z\big)
         =\sum_{z\in\Sing(\tilde{F}_t^\i)\cap U_t(z_i)} \mu\big(\tilde{F}_t^\i,z\big).
       \end{displaymath}
       \tom{\footnote{
       Let now $i\in\{1,\ldots,r\}$. For $t\not=0$ fixed we consider the transformation
       ${\psi_t'}^i$ defined by the coordinate change $\psi_t^i$,
       \begin{displaymath}
         \begin{array}{rcl}
           \C^2\;\;\supset\;\; U(z_i)\;\; &\longrightarrow & \;\;U_t(z_i)\;\;\subset\;\;\C^2\\
           {\begin{sideways}$\in$\end{sideways}}\;\;\;\;\; & &\;\;\;\;\;{\begin{sideways}$\in$\end{sideways}}\\
           (x_i,y_i) & \mapsto & \big(\frac{1}{t}x_i,\frac{1}{t}y_i\big),
         \end{array}
       \end{displaymath}
       and the transformed equations
       \begin{displaymath}
         \tilde{F}_t^\i(x_i,y_i)= t^{-m_i+1}F_t^i(tx_i,ty_i)=0.
       \end{displaymath}
       Then obviously,
       \begin{displaymath}
         F_t^\i(z)=0 \;\;\;\Longleftrightarrow\;\;\;\big(\tilde{F}_t^\i\circ{\psi_t'}^i\big)(z)=0,
       \end{displaymath}
       and 
       \begin{displaymath}
         \nabla{\psi_t'}^i\equiv\frac{1}{t}\id_{\C^2}.
       \end{displaymath}
       Thus we have
       \begin{displaymath}
         (m_i-2)^2 = \mu(C_0,z_i)
         =\sum_{z\in\Sing(F_t^\i)\cap U(z_i)} \mu\big(F_t^\i,z\big)
         =\sum_{z\in\Sing(\tilde{F}_t^\i)\cap U_t(z_i)} \mu(\tilde{F}_t^\i,z).
       \end{displaymath}
       }}
       For $t\not=0$ very small $U_t(z_i)$ becomes very large, so
       that, by shrinking $V(0)$ we may suppose that for any
       $0\not=t\in V(0)$
       \begin{displaymath}
         \Sing(g_i)\subset U_t(z_i),  
       \end{displaymath}
       and that for any $z\in\Sing(g_i)$
       there is an open neighbourhood $U(z)\subset U_t(z_i)$ such that
       \begin{displaymath}
         \mu(g_i,z) = \sum_{{z'}\in\Sing(\tilde{F}_t^\i)\cap U(z)}\mu\big(\tilde{F}_t^\i,{z'}\big).
       \end{displaymath}
       If we now take into account that $g_i$ has precisely one
       critical point, $z_i$, on its zero level, and that the critical points
       on the zero level of $\tilde{F}_t^\i$ all contribute to the
       Milnor number $\mu(g_i,z_i)$, then we get
       the following sequence of inequalities:\tom{\footnote{Note, an
           ordinary plane curve singularity $(X,x)$ of multiplicity
           $m$ has Milnor number
           $\mu(X,x)=\dim_\C\big(\C\{x,y\}/(x^{m-1},y^{m-1})\big)=(m-1)^2$. 
           And, for an affine plane curve given by an equation $g$
           such that the equation has no critical point at infinity we 
           have by B\'ezout's Theorem: 
           \begin{displaymath}
             (m-1)^2=\sum_{p\in\A^2}i\left(\frac{\partial g}{\partial
               x},\frac{\partial g}{\partial y};p\right),
           \end{displaymath}
           where by definition $i\big(\frac{\partial g}{\partial
             x},\frac{\partial g}{\partial y};p\big)$ and $\mu(g,p)$ 
           are defined as the dimension of the same vector spaces.}}
       \begin{displaymath}
         (m_i-2)^2-\mu(\ks_i) =
         \sum_{z\in\Sing(g_i)}\mu(g_i,z) - \sum_{z\in\Sing(g_i^{-1}(0))}\mu(g_i,z) 
       \end{displaymath}
       \begin{displaymath}
         \leq \sum_{z\in\Sing(\tilde{F}_t^\i)\cap U_t(z_i)}\mu\big(\tilde{F}_t^i,z\big)
         - \sum_{z\in\Sing((\tilde{F}_t^\i)^{-1}(0))\cap U_t(z_i)}\mu\big(\tilde{F}_t^i,z\big)
       \end{displaymath}
       \begin{displaymath}
         = \sum_{z\in\Sing(F_t^\i)\cap U(z_i)}\mu(F_t^i,z)
         - \sum_{z\in\Sing((F_t^\i)^{-1}(0))\cap U(z_i)}\mu(F_t^i,z)
       \end{displaymath}
       \begin{displaymath}
         \leq \mu(C_0,z_i)-\mu\big(F_t^\i,z_i\big)
         =(m_i-2)^2-\mu(\ks_i).
       \end{displaymath}
       Hence all inequalities must have been equalities, and, in
       particular, 
       \begin{displaymath}
         \Sing(C_t)\cap U(z_i)=\Sing\big((F_t^\i)^{-1}(0)\big)\cap U(z_i) =\{ z_i\},
       \end{displaymath}
       which in view of Condition (\ref{eq:existence-III:1}) finishes
       the proof.

       Note that $C_t$, being a small deformation of the irreducible
       reduced curve $C_0$, will again be irreducible and reduced.
   \end{proof}

   \begin{corollary}\label{cor:existence-IV}
     Let $L\in\Div(\Sigma)$ be very ample over
     $\C$. Suppose that $D\in\Div(\Sigma)$ and $\ks_1,\ldots,\ks_r$ are
     topological singularity types with 
     $\mu(\ks_1)\geq\ldots\geq\mu(\ks_r)$  such that 
     \begin{equationlist}
       \item[eq:existence-IV:1] $(D-L-K_{\Sigma})^2
         \geq \frac{414}{5}\sum\limits_{\mu(\ks_i)\leq 38} \mu(\ks_i) +
         58\sum\limits_{\mu(\ks_i)\geq 39} \Big(\sqrt{\mu(\ks_i)}+\frac{13}{2\sqrt{29}}\Big)^2$,
       \item[eq:existence-IV:2] $(D-L-K_\Sigma).B >
         \left\{
           \begin{array}{ll}
             \sqrt{\frac{207}{5}}\sqrt{\mu(\ks_1)}-1, & \mbox{ if 
               }\mu(\ks_1)\leq 38,\\[0.3cm]
             \sqrt{29} \sqrt{\mu(\ks_1)}+\frac{11}{2}, & \mbox{
               if }\mu(\ks_1)\geq 39,
           \end{array}
         \right.$
         \\[0.3cm]
         for any irreducible curve $B$ with $B^2=0$ and $\dim|B|_a>0$,
       \item[eq:existence-IV:3] $D-L-K_\Sigma$ is nef,
       \item[eq:existence-IV:4] $D.L-2 g(L)\geq 
         \left\{
           \begin{array}{ll}
             \sqrt{\frac{207}{5}}\Big(\sqrt{\mu(\ks_1)}+\sqrt{\mu(\ks_2)}\Big) -2, & \mbox{ if 
               }\mu(\ks_1)\leq 38,\\[0.3cm]
             \sqrt{\frac{207}{5}}\sqrt{\mu(\ks_1)}+\sqrt{29}
             \sqrt{\mu(\ks_1)}+\frac{9}{2}, & \mbox{ if }
             \mu(\ks_1)\geq 39\\[0.2cm]
             & \mbox{ \& }\mu(\ks_2)\leq 38,\\[0.3cm]
             \sqrt{29} \Big(\sqrt{\mu(\ks_1)}+\sqrt{\mu(\ks_2)}\Big)+ 11, & \mbox{
               if }\mu(\ks_2)\geq 39,
           \end{array}
         \right.$
       \item[eq:existence-IV:5] $D^2
         \geq \frac{207}{5}\sum\limits_{\mu(\ks_i)\leq 38} \Big(\sqrt{\mu(\ks_i)}-\sqrt{\frac{5}{207}}\Big)^2+
         29\sum\limits_{\mu(\ks_i)\geq 39} \Big(\sqrt{\mu(\ks_i)}+\frac{11}{2\sqrt{29}}\Big)^2$,
     \end{equationlist}
     then there is an irreducible reduced curve $C$ in $|D|_l$ with
     $r$ singular points of topological types $\ks_1,\ldots,\ks_r$ as
     its only singularities.
   \end{corollary}

   \begin{proof}
     This follows right away from Corollary \ref{cor:existence-II},
     Theorem \ref{thm:existence-III}, and \cite{Los98} Theorem
     4.2.\tom{\footnote{Note that in \cite{Los98} Theorem 4.2 the
         first condition is better for $\mu(\ks_i)\geq 39$ while the
         second one is better for $\mu(\ks_i)\leq 38$, and note that
         the strict inequalities in Corollary \ref{cor:existence-II} are
         absorbed by \cite{Los98} Theorem 4.2.}}
   \end{proof}

   \begin{remark}
     One could easily simplify the above formulae by not
     distinguishing the cases $\mu(\ks_i)\geq 39$ and $\mu(\ks_i)\leq
     38$. However, one would loose information. 

     On the other hand, knowing something more about the singularity
     type one could achieve much better results, applying the
     corresponding bounds for the $s(\ks_i)$. We leave it to the
     reader to apply the bounds. (Cf.~\cite{Los98} Remarks 4.3, 4.8, and
     4.15) 

     As we have already mentioned earlier\tom{ (see Remark
       \ref{rem:vanishing})} the most restrictive of the above \emph{sufficient} conditions is 
     (\ref{eq:existence-IV:1}), which could be characterised 
     as a condition of the type
     \begin{displaymath}
       \sum_{i=1}^r\mu(\ks_i)\leq \alpha D^2 + \beta D.K + \gamma,
     \end{displaymath}
     where $K$ is some fixed divisor class, $\alpha,\beta$ and
     $\gamma$ are some constants.  
     
     There are also \emph{necessary} conditions of this type, e.~g.~
     \begin{displaymath}
       \sum_{i=1}^r\mu(\ks_i)\leq D^2 +D.K_\Sigma+2,
     \end{displaymath}
     which follows from the genus formula.\footnote{If $D$ is an
       irreducible curve with precisely $r$ singular points of types
       $\ks_1,\ldots,\ks_r$
       and $\nu:\widetilde{D}\rightarrow D$ its normalisation, then
       $p_a(D)=g\big(\widetilde{D}\big)+\delta(D)\geq\delta(D)$, where
       $\delta(D)=\dim_\C\big(\nu_*\ko_{\widetilde{D}}/\ko_D\big)$ is
       the delta-invariant of $D$ (cf.~\cite{BPV84}
         II.11). Moreover, by definition
       $\delta(D)=\sum_{z\in\Sing(D)}\delta(D,z)$ where
       $\delta(D,z)=\dim_\C\big(\nu_*\ko_{\widetilde{D}}/\ko_D\big)_z$ 
       is the local delta-invariant at $z$, and it is well known that
       $2\delta(D,z)=\mu(D,z)+r(D,z)-1\geq\mu(D,z)$, where $r(D,z)$ is
       the number of branches of the curve singularity $(D,z)$ and
       $\mu(D,z)$ is its Milnor number (cf.~\cite{Mil68}
         Chapter 10).
       Using now the genus formula we get
       \begin{displaymath}
         D^2+D.K_\Sigma+2=2p_a(D)\geq 2\sum_{i=1}^r\delta(\ks_i)\geq\sum_{i=1}^r\mu(\ks_i).
       \end{displaymath}
       }
     
     See \cite{Los98} Section 4.1 for considerations on the
     asymptotical optimality of the constant $\alpha$.
   \end{remark}


   \section{Examples}\label{sec:examples}
   \setcounter{equation}{0}

   In this section we are going to examine the conditions in the
   vanishing theorem (Corollary \ref{cor:vanishing}) and in the
   corresponding existence results  for various types of
   surfaces. Unless otherwise stated, $r\geq 1$ is a positive
   integer, and $m_1,\ldots,m_r\in\N_0$ are non-negative, while at
   least one $m_i$ is positive whenever we consider conditions for 
   existence theorems. 

   \subsection{The Classical Case - $\Sigma=\PC^2$}

   Since in $\PC^2$ there are no irreducible curves of
   self-intersection number zero, Condition (\ref{eq:vanishing:2}) is
   redundant. Moreover, Condition (\ref{eq:vanishing:3}) takes in view 
   of (\ref{eq:vanishing:1}) the form $d+3\geq \sqrt{2}$. Corollary
   \ref{cor:vanishing} thus takes the following 
   form, where $L\in|\ko_{\PC^2}(1)|_l$ is a generic line.  

   \renewcommand{\thesatz}{\ref{cor:vanishing}\alph{subsection}}
   \begin{corollary}
     Let $d$ be any integer such that 
     \begin{equationlist}
        \theequationchange{\ref{eq:vanishing:1}\alph{subsection}}
        \item $(d+3)^2 \geq 2\sum_{i=1}^r(m_i+1)^2$, 
        \theequationchange{\ref{eq:vanishing:3}\alph{subsection}}
        \item $d\geq -1$.
        \theequationchangeback
     \end{equationlist}
     Then for $z_1,\ldots,z_r\in\PC^2$ in very general position and $\nu>0$
     \begin{displaymath}
       H^\nu\left(\Bl_{\underline{z}}\big(\PC^2\big),d\pi^*L-\sum_{i=1}^r m_iE_i\right)=0.
     \end{displaymath}
   \end{corollary}
   \renewcommand{\thesatz}{\remembersatz}

   Now turning to the existence theorem Corollary \ref{cor:existence-II} for
   generic fat point schemes, we, of course, find that
   Condition (\ref{eq:existence-II:2}) is 
   obsolete, and so is (\ref{eq:existence-II:3}), taking into account
   that (\ref{eq:existence-II:4}) implies $d>0$. But then Conditions
   (\ref{eq:existence-II:4}) and  
   (\ref{eq:existence-II:5}) become also redundant in view of Condition
   (\ref{eq:existence-II:1}) and equation
   (\ref{eq:existence-II:6}).\tom{\footnote{Using the second
       binomial formula, (\ref{eq:existence-II:1}) implies, 
       \begin{displaymath}
         (d+2)^2\geq 2\big((m_i+1)^2+(m_j+1)^2)\geq(m_i+1+m_j+1)^2,
       \end{displaymath}
       and thus (\ref{eq:existence-II:4}) follows.\\
       For $d=1,2,3$ it follows easily from
       (\ref{eq:existence-II:1}\alph{subsection}) that there are
       irreducible curves of degree $d$ with the desired
       singularities. We therefore may assume $d\geq 4$, 
       and thus (\ref{eq:existence-II:6}), taking the form 
       \begin{displaymath}
         d(d-4)+4\left(\sum_{i=1}^rm_i-1\right)+2r >0,
       \end{displaymath}
       is fulfilled and implies (\ref{eq:existence-II:5}).
       }} 

   Thus the conditions in Corollary \ref{cor:existence-II} reduce
   to $d>0$ and
   \begin{equationlist}
      \theequationchange{\ref{eq:existence-II:1}\alph{subsection}}
      \item $(d+2)^2\geq 2\sum_{i=1}^r 
        (m_i+1)^2$, 
      \theequationchangeback
   \end{equationlist}
   and, similarly, the conditions in Corollary \ref{cor:existence-IV} 
   reduce to $d>7$ and
   \theequationchange{\ref{eq:existence-IV:1}\alph{subsection}}
   \begin{equationlist}
      \item $(d+2)^2
         \geq \frac{414}{5}\sum\limits_{\mu(\ks_i)\leq 38} \mu(\ks_i) +
         58\sum\limits_{\mu(\ks_i)\geq 39}
         \Big(\sqrt{\mu(\ks_i)}+\frac{13}{2\sqrt{29}}\Big)^2$. 
   \end{equationlist} 
   \theequationchangeback
   These results are much weaker than the previously known ones
   (e.~g.~\cite{Los98} Proposition 4.11, where the factor $2$
   is replaced by $\frac{10}{9}$) which use the Vanishing
   Theorem of Geng Xu (cf.~\cite{Xu95} Theorem 3), particularly
   designed for $\PC^2$. -- Using $L\in|\ko_\Sigma(l)|_l$ with $l>1$ 
   instead of $\ko_\Sigma(1)$ in Corollary \ref{cor:existence-II} 
   does not improve the conditions.\tom{\footnote{As soon as we consider $d\geq 
       \frac{l}{2}+2$  again only
       Condition (\ref{eq:existence-II:1}) remains but is worse than
       in the case $l=1$.
       For this, note that as above (\ref{eq:existence-II:1})
       implies $d-l+1\geq m_i+m_j$. Taking this into account, 
       Condition (\ref{eq:existence-II:5}),
       $dl-\frac{(l-1)(l-2)}{2}\geq 
       m_i+m_j$, is fulfilled as soon as $d\geq\frac{l}{2}+2$.}}


   \subsection{Geometrically Ruled Surfaces}\label{subsec:ruled-surfaces}

   Let $\xymatrix@C0.6cm{\Sigma=\P(\ke)\ar[r]^(0.65)\pi & C}$ be a
   geometrically ruled surface with normalised bundle $\ke$ (in the
   sense of \cite{Har77} V.2.8.1). The N\'eron-Severi group of $\Sigma$ is 
   \begin{displaymath}
     \NS(\Sigma) = C_0\Z\oplus F\Z,
   \end{displaymath}
   with intersection matrix
   \begin{displaymath}
     \left(\!\begin{array}{rc}-e & 1 \\ 1 & 0\end{array}\right),
   \end{displaymath}
   where $F\cong\PC^1$ is a fibre of $\pi$, $C_0$ a section of $\pi$
   with $\ko_\Sigma(C_0)\cong\ko_{\P(\ke)}(1)$, and
   $e=-\deg(\Lambda^2\ke)\geq -g$.\footnote{By \cite{Nag70} Theorem 1
     there is some section $D\sim_a C_0+bF$ with $g\geq
     D^2=2b-e$. Since $D$ is irreducible, by \cite{Har77} V.2.20/21
     $b\geq 0$, and thus $-g\leq e$.} 
   For the canonical divisor we have 
   \begin{displaymath}
     K_\Sigma \sim_a -2C_0+ (2g-2-e)F,
   \end{displaymath}
   where $g=g(C)$ is the genus of the base curve $C$.

   In order to understand Condition (\ref{eq:vanishing:2}) we have to
   examine special irreducible curves on $\Sigma$.

   \problem{1}{ \textbf{!!! DER BEWEIS DES LEMMAS IST FALSCH !!!}

     \begin{lemma}\label{lem:ruled-surface}
       Let $B\in|aC_0+bF|_a$ be an irreducible curve
       with $B^2=0$ and $\dim|B|_a\geq 1$. Then we are in one of the
       following two cases 
       \begin{equationlist}
       \item[eq:ruled-surface:1] $a=0$, $b=1$, and $B\sim_a F$, or
       \item[eq:ruled-surface:2] $a=1$, $b=0$, $B\sim_a C_0$, and
         $\Sigma\cong C_0\times\PC^1$.
       \end{equationlist}
     \end{lemma}
     \begin{proof}
       Obviously $F$ and $C_0$ are irreducible curves satisfying the cases
       (\ref{eq:ruled-surface:1}) and (\ref{eq:ruled-surface:2})
       respectively. 
       
       Let us now assume that $B\subset \Sigma$ is an irreducible curve with
       $B\not\sim_a F$. We have to show that we 
       are in case (\ref{eq:ruled-surface:2}).
       
       By Proposition \ref{prop:condition} $B$ gives rise to a second
       fibration $\pi':\Sigma\rightarrow H$, where the elements of
       $|B|_a$ are just the fibres of $\pi'$. Then\tom{, by
         \cite{Har77} V.2.21,} $a=B.F>0$ and thus the restriction 
       $\pi'_|:F\rightarrow H$
       of $\pi'$ to $F$ is a finite morphism of smooth projective
       curves. Hurwitz's formula therefore gives
       \begin{displaymath}
         0=g(F)\geq a\big(g(H)-1\big)+1\geq 0.
       \end{displaymath}
       Hence, $g(H)=0$ and $a=1$. In particular, \tom{since $B.F=1$, }the elements of $|B|_a$ 
       are disjoint sections of the fibration $\pi$. 
       
       Thus, by Lemma \ref{lem:ruledproduct}, $\Sigma\cong
       C\times\PC^1$, and hence $e=0$. 
       But then $0=B^2=a\big(b-\frac{a}{2}e\big)=b$, and we are in case
       (\ref{eq:ruled-surface:2}). 
     \end{proof}

     \begin{corollary}
       A geometrically ruled surface admits a second fibration (with
       irreducible fibres) if and only if it is trivial, i.~e.~is
       isomorphic to the product of the base curve with $\PC^1$.  
     \end{corollary}     
     }
   \problem{2}{

     \begin{lemma}\label{lem:ruled-surface}
       Let $B\in|aC_0+bF|_a$ be an irreducible curve
       with $B^2=0$ and $\dim|B|_a\geq 1$. Then we are in one of the
       following cases 
       \begin{equationlist}
       \item[eq:ruled-surface:1] $a=0$, $b=1$, and $B\sim_a F$, 
       \item[eq:ruled-surface:2] $e=0$, $a\geq 1$, $b=0$, and $B\sim_a aC_0$, or
       \item[eq:ruled-surface:3] $e<0$, $a\geq 2$, $b=\frac{a}{2}e<0$,
         and $B\sim_a aC_0+\frac{a}{2}eF$.  
       \end{equationlist}
       Moreover, if $a=1$, then $\Sigma\cong C_0\times\PC^1$.
     \end{lemma}
     \begin{proof}
       Since $B$ is irreducible, we have 
       \begin{equation}\label{eq:nachtrag:1}
         0\leq B.F=a \;\;\;\;\mbox{ and }\;\;\;\; 0\leq B.C_0=b-ae.
       \end{equation}
       If $a=0$, then $|B|_a=|bF|_a$, but since the general element of
       $|B|_a$ is irreducible, $b$ has to be one, and we are in case
       (\ref{eq:ruled-surface:1}). 

       We, therefore, may assume that $a\geq 1$. Since $B^2=0$ we have
       \begin{equation}\label{eq:nachtrag:2}
         0=B^2=2a\left(b-\tfrac{a}{2}e\right),\;\;\;\mbox{ hence }\;\;\;b=\tfrac{a}{2}e.
       \end{equation}
       Combining this with (\ref{eq:nachtrag:1}) we get
       $e\leq 0$. 

       Moreover, if $e=0$, then of course $b=0$, while, if $e<0$, then
       $a\geq 2$ by \cite{Har77} V.2.21, since otherwise $b$ would
       have to be non-negative. This brings us down to the cases
       (\ref{eq:ruled-surface:2}) and (\ref{eq:ruled-surface:3}). 

       It remains to show, that $B.F=a=1$ implies $\Sigma\cong
       C_0\times\PC^1$. But by assumption the elements of $|B|_a$ 
       are disjoint sections of the fibration $\pi$. 
       Thus, by Lemma \ref{lem:ruledproduct}, $\Sigma\cong
       C\times\PC^1$. 
     \end{proof}     
     }

   \begin{lemma}\label{lem:ruledproduct}
     If $\pi:\Sigma\rightarrow C$ has three disjoint sections, then
     $\Sigma$ is isomorphic to $C\times\PC^1$ as a ruled surface,
     i.~e.~there is an isomorphism $\alpha:\Sigma\rightarrow
     C\times\PC^1$ such that the following diagram is commutative: 
     \begin{displaymath}
       \xymatrix@C0.6cm{
         \Sigma\ar[rr]^(0.45)\alpha\ar[dr]_\pi&& C\times\PC^1\ar[dl]^\pr\\
         &C.&
         }
     \end{displaymath}
   \end{lemma}
   \begin{proof}
     See \cite{IS96} p.~229.

     $\pi$ is a locally trivial $\PC^1$-bundle, thus $C$ is covered by a finite 
     number of open affine subsets $U_i\subset C$ with 
     trivialisations 
     \begin{displaymath}
       \xymatrix@C0.6cm{
         \pi^{-1}(U_i)\ar[rr]^{\cong}_{\varphi_i}\ar[dr]_{\pi} &&
         U_i\times\PC^1\ar[dl]^{\pr}\\
         & U_i,&
         }
     \end{displaymath}
     which are linear on the fibres.

     The three disjoint sections on $\Sigma$, say $S_0$, $S_1$, and
     $S_\infty$, give rise to three sections $S^i_0$, $S^i_1$, and
     $S^i_\infty$ on $U_i\times \PC^1$. For each point $z\in U_i$
     there is a unique linear projectivity on the fibre
     $\{z\}\times\PC^1$ mapping the three points
     $p_{0,z}=S^i_0\cap\big(\{z\}\times\PC^1\big)$,
     $p_{1,z}=S^i_1\cap\big(\{z\}\times\PC^1\big)$, and
     $p_{\infty,z}=S^i_\infty\cap\big(\{z\}\times\PC^1\big)$ to the
     standard basis $0\equiv\big(z,(1:0)\big)$, $1\equiv\big(z,(1:1)\big)$, and
     $\infty\equiv\big(z,(0:1)\big)$ of $\PC^1\cong\{z\}\times\PC^1$.
     If $p_{0,z}=\big(z,(x_0:y_0)\big)$, $p_{1,z}=\big(z,(x_1:y_1)\big)$, and
     $p_{\infty,z}=\big(z,(x_\infty:y_\infty)\big)$, the projectivity is
     given by the matrix
     \begin{displaymath}
       A=\left(
         \begin{array}{cc}
           \frac{(x_0y_1-y_0x_1)y_\infty}{y_0y_1x_\infty^2-y_0x_1x_\infty y_\infty-x_0y_1x_\infty y_\infty+x_0x_1y_\infty} &
           \frac{(x_0y_1-y_0x_1)x_\infty}{y_0y_1x_\infty^2-y_0x_1x_\infty y_\infty-x_0y_1x_\infty y_\infty+x_0x_1y_\infty}\\
           \frac{y_0}{x_0y_\infty-y_0x_\infty}& \frac{x_0}{x_0y_\infty-y_0x_\infty}
         \end{array}
       \right),
     \end{displaymath}
     whose entries are rational functions in the coordinates of $p_{0,z}$, $p_{1,z}$, and
     $p_{\infty,z}$.
     Inserting for the coordinates local equations of the sections, 
     $A$ finally gives rise to an isomorphism of $\PC^1$-bundles\tom{,
     i.~e.~a morphism which is a linear isomorphism on the fibres,}
     \begin{displaymath}
       \alpha_i:U_i\times\PC^1\rightarrow U_i\times\PC^1
     \end{displaymath}
     mapping the sections $S^i_0$, $S^i_1$, and $S^i_\infty$ to the
     trivial sections. 

     The transition maps 
     \begin{displaymath}
       \xymatrix@C1cm{
         U_{ij}\times\PC^1\ar[r]^{\alpha_{i|}^{-1}} &
         U_{ij}\times\PC^1\ar[r]^{\varphi_{i|}^{-1}} &
         \pi^{-1}(U_{ij})\ar[r]^{\varphi_{j|}} &
         U_{ij}\times\PC^1\ar[r]^{\alpha_{j|}} &
         U_{ij}\times\PC^1
         },
     \end{displaymath}
     with $U_{ij}=U_i\cap U_j$,
     are linear on the fibres and fix the three trivial sections.
     Thus they must be the identity maps, which implies
     that the $\alpha_i\circ\varphi_i$, $i=1,\ldots,r$, glue together to an isomorphism
     of ruled surfaces:     
     \begin{displaymath}
       \xymatrix@C0.6cm{
         \Sigma\ar[rr]^(0.45)\alpha\ar[dr]_\pi&& C\times\PC^1\ar[dl]^\pr\\
         &C.&
         }       
     \end{displaymath}
   \end{proof}

   Knowing the algebraic equivalence classes of irreducible curves
   in $\Sigma$ which satisfy the assumptions in Condition
   (\ref{eq:vanishing:2}) we can give a better formulation of the
   vanishing theorem in 
   the case of geometrically ruled surfaces. 

   In order to do the
   same for the existence theorems, we have
   to study very ample divisors on $\Sigma$. These, however,
   depend very much on the structure of the base curve
   $C$,\tom{\footnote{Cf.~\cite{Har77} V.2.22.2, Ex.~V.2.11 and
       Ex.~V.2.12.}} 
   and the general results which we give may be not the best
   possible. Only in the case 
   $C=\PC^1$ we can give a
   complete investigation.

   The geometrically ruled surfaces with base curve $\PC^1$ are, up to 
   isomorphism, just the Hirzebruch surfaces
   $\F_e=\P\big(\ko_{\PC^1}\oplus\ko_{\PC^1}(-e)\big)$, $e\geq
   0$. Note that $\Pic(\F_e)=\NS(\F_e)$, that is, algebraic 
   equivalence and linear equivalence coincide.  
   Moreover, by \cite{Har77} V.2.18 a divisor class $L=\alpha
   C_0+\beta F$ is very ample over $\C$ if and  
   only if $\alpha>0$ and $\beta>\alpha e$.
   The conditions throughout the existence theorems turn out to be optimal
   if we work with $L=C_0+(e+1)F$, while for other choices of $L$ they 
   become more
   restrictive.\footnote{\setcounter{subsubsection}{1}
       \setcounter{sss}{\value{subsubsection}}\stepcounter{sss}
       Let $L'=\alpha C_0+\beta F$, then 
       $D-L'-K_{\F_e}=(a+1-\alpha)C_0+(b+1+e-\beta)F$, and thus the 
       optimality of the conditions follows from
       \begin{equationlist}
       \theequationchange{\ref{eq:existence-II:1}\alph{subsection}.\roman{subsubsection}/\roman{sss}}
       \item $(D-L'-K_{\F_e})^2=
         (a+1-\alpha)\big(2(b+1+e-\beta)-(a+1-\alpha)e\big) \leq
         a \big((2b-ae)+(\alpha e+e+2-2\beta)\big)\leq a(2b-ae)=(D-L-K_{\F_e})^2$, 
       \theequationchange{\ref{eq:existence-II:2}\alph{subsection}.\roman{subsubsection}/\roman{sss}}
       \item $(D-L'-K_{\F_e}).F=a+1-\alpha\leq a=(D-L-K_{\F_e}).F$, and for $e=0$,
         $(D-L'-K_{\F_e}).C_0=b+1-\beta\leq b=(D-L-K_{\F_e}).C_0$, and
       \theequationchange{\ref{eq:existence-II:3}\alph{subsection}.\roman{sss}}
       \item $b+1+e-\beta \geq e(a+1-\alpha)$ implies $b\geq
         b+e\alpha+1-\beta \geq ae$.
       \theequationchangeback
       \end{equationlist}
       \addtocounter{subsubsection}{-1}
       }

   In the case $C\not\cong\PC^1$, we may choose an integer
   $l\geq\max\{e+1,2\}$ such that the algebraic  
   equivalence class $|C_0+lF|_a$ contains a very ample divisor
   $L$, e.~g.~$l=e+3$ will do, if $C$ is an elliptic
   curve.\footnote{$l$ will be the degree of a suitable
     very ample divisor $\mathfrak{d}$ on $C$. Now $\mathfrak{d}$
     defines an embedding of $C$ into some $\PC^N$ such that 
     the degree of the image $C'$ is just $\deg(\mathfrak{d})$. Therefore
     $\deg(\mathfrak{d})\geq 2$, unless $C'$ is linear\tom{
       (cf.~\cite{Har77} Ex.~I.7.6)}, which implies $C\cong\PC^1$.}
   In particular, $l\geq 2$ as soon as $\Sigma\not\cong\PC^1\times\PC^1$.

   With the above choice of $L$ we have
   $g(L)=1+\frac{L^2+L.K_{\Sigma}}{2}=1+\frac{(-e+2l)+(e-2l+2g-2)}{2}=g$, and
   hence the generic curve in $|L|_l$ is a smooth curve whose genus
   equals the genus of the base curve.
   
\problem{1}{

   \subsubsection{The trivial geometrically ruled surfaces
     $\Sigma\simeq C\times\PC^1$}

   Since the case $\Sigma\cong C\times\PC^1$ is exceptional 
   by Lemma \ref{lem:ruled-surface}, we will consider it
   separately. In view of
   (\ref{eq:vanishing:2}\alph{subsection}.\roman{subsubsection}) the
   Condition 
   (\ref{eq:vanishing:3}) becomes obsolete, and Corollary
   \ref{cor:vanishing} takes the following form.

   \renewcommand{\thesatz}{\ref{cor:vanishing}\alph{subsection}.\roman{subsubsection}}
   \begin{corollary}
     Let $a,b\in\Z$ be two 
     integers satisfying 
     \begin{equationlist}
       \theequationchange{\ref{eq:vanishing:1}\alph{subsection}.\roman{subsubsection}}
       \item $ab\geq \sum_{i=1}^r(m_i+1)^2$, and
       \theequationchange{\ref{eq:vanishing:2}\alph{subsection}.\roman{subsubsection}}
       \item $a,b>\max\{m_i\;|\;i=1,\ldots,r\}$,
       \theequationchangeback
     \end{equationlist}
     then for $z_1,\ldots,z_r\in C\times\PC^1$ in very general
     position and $\nu>0$
     \begin{displaymath}
       H^\nu\left(\Bl_{\underline{z}}(C\times\PC^1),\;
         (a-2)\pi^*C_0+(b-2+2g)\pi^*F-\sum_{i=1}^r m_iE_i\right)=0. 
     \end{displaymath}
   \end{corollary}
   \renewcommand{\thesatz}{\remembersatz}

   \begin{proof}
     Setting $D=(a-2)C_0+(b-2+2g)F$
     we have
     \begin{displaymath}
       \Big(D-K_{C\times\PC^1}\Big)^2=(aC_0+bF)^2=2ab\geq 2\sum_{i=1}^r(m_i+1)^2,
     \end{displaymath}
     which is just (\ref{eq:vanishing:1}).
     Similarly (\ref{eq:vanishing:2}\alph{subsection}.\roman{subsubsection}) and Lemma
     \ref{lem:ruled-surface} ensure that Condition (\ref{eq:vanishing:2}) is
     satisfied. 

     By \cite{Har77} V.2.20, if a curve $B\sim_a (a'C_0+b'F)$ is
     irreducible, then either $a'=0$ and $b'=1$, or $a'>0$ and $b'\geq 
     0$. Thus $\big(D-K_{C\times\PC^1}\big).B$ is either $a$ or $ab'+ba'$,
     but in any case positive by
     (\ref{eq:vanishing:2}\alph{subsection}.\roman{subsubsection}). Hence
     $D-K_{C\times\PC^1}$ is nef (even ample), which is Condition
     (\ref{eq:vanishing:3}). 
   \end{proof}

   In order to obtain nice formulae we considered $D=(a-2)C_0+(b-2+2g)F$
   in the formulation of the vanishing theorem. For the existence
   theorems it turns out that the formulae look best if we work with
   $D=(a-1)C_0+(b+l+2g-2)F$ instead. In the case of Hirzebruch
   surfaces this is just $D=(a-1)C_0+(b-1)F$.
   Some calculations show that (\ref{eq:existence-II:3}),
   (\ref{eq:existence-II:4}), and (\ref{eq:existence-II:5}) are
   obsolete,\footnote{\setcounter{sss}{2} See the proof of Corollary
     \ref{cor:existence-II}\alph{subsection}.\roman{sss} below.}
   and we thus receive the following theorem.

   \renewcommand{\thesatz}{\ref{cor:existence-II}\alph{subsection}.\roman{subsubsection}}
   \begin{corollary}
     Given integers $a,b\in\Z$ satisfying 
     \begin{equationlist}
       \theequationchange{\ref{eq:existence-II:1}\alph{subsection}.\roman{subsubsection}}
       \item $ab\geq\sum_{i=1}^r (m_i+1)^2$, and 
       \theequationchange{\ref{eq:existence-II:2}\alph{subsection}.\roman{subsubsection}}
       \item $a,b>\max\{m_i\;|\;i=1,\ldots,r\}$,  
       \theequationchangeback
     \end{equationlist}
     then for $z_1,\ldots,z_r\in C\times\PC^1$ in very 
     general position there is an irreducible reduced curve
     $C\in|(a-1)C_0+(b+l+2g-2)F|_a$ with ordinary singularities of
     multiplicities $m_i$ at the $z_i$ as only singularities. 
     Moreover, $V_{|C|}(\underline{m})$ is T-smooth at $C$. 
   \end{corollary}
   \renewcommand{\thesatz}{\remembersatz}

   With the same $D$ and $L$ as above the conditions in the existence 
   theorem Corollary \ref{cor:existence-IV} reduce to 
   \begin{equationlist}
     \theequationchange{\ref{eq:existence-IV:1}\alph{subsection}.\roman{subsubsection}}
     \item $ab 
         \geq \frac{207}{5}\sum\limits_{\mu(\ks_i)\leq 38} \mu(\ks_i) +
         29\sum\limits_{\mu(\ks_i)\geq 39}
         \Big(\sqrt{\mu(\ks_i)}+\frac{13}{2\sqrt{29}}\Big)^2$, and
     \theequationchange{\ref{eq:existence-IV:2}\alph{subsection}.\roman{subsubsection}}
     \item $a,b >
         \left\{
           \begin{array}{ll}
             \sqrt{\frac{207}{5}}\sqrt{\mu(\ks_1)}-1, & \mbox{ if 
               }\mu(\ks_1)\leq 38,\\[0.3cm]
             \sqrt{29} \sqrt{\mu(\ks_1)}+\frac{11}{2}, & \mbox{
               if }\mu(\ks_1)\geq 39.
           \end{array}
         \right.$ 
     \theequationchangeback
   \end{equationlist}

   \subsubsection{The non-trivial geometrically ruled surfaces $\Sigma\not\simeq C\times\PC^1$} 

   Let us now consider the case of a non-trivial geometrically ruled surface.

   \renewcommand{\thesatz}{\ref{cor:vanishing}\alph{subsection}.\roman{subsubsection}}
   \begin{corollary}
     Given two integers $a,b\in\Z$ satisfying 
     \begin{equationlist}
       \theequationchange{\ref{eq:vanishing:1}\alph{subsection}.\roman{subsubsection}}
       \item $a\big(b-\big(\frac{a}{2}-1\big)e\big)\geq \sum_{i=1}^r(m_i+1)^2$, 
       \theequationchange{\ref{eq:vanishing:2}\alph{subsection}.\roman{subsubsection}}
       \item $a>\max\{m_i\;|\;i=1,\ldots,r\}$, and
       \theequationchange{\ref{eq:vanishing:3}\alph{subsection}.\roman{subsubsection}}
       \item $b\geq (a-1)e$, if $e>0$.
       \theequationchangeback
     \end{equationlist}
     For $z_1,\ldots,z_r\in\Sigma$ in very general position and
     $\nu>0$ 
     \begin{displaymath}
       H^\nu\left(\Bl_{\underline{z}}(\Sigma),(a-2)\pi^*C_0+(b-2+2g)\pi^*F-\sum_{i=1}^r m_iE_i\right)=0.
     \end{displaymath}
   \end{corollary}
   \renewcommand{\thesatz}{\remembersatz}

   \begin{proof}
     Note that if the invariant $e$ is non-positive, then
     $\big(b-\big(\frac{a}{2}-1\big)e\big)>0$ implies 
     \begin{equation} 
       \label{eq:vanishing:3*}
       b\geq(a-1)e,
     \end{equation}
     so that this inequality is fulfilled for any choice of $e$.

     Setting $D=(a-2)C_0+(b-2+2g)F$
     we have
     \begin{displaymath}
       (D-K_{\Sigma})^2=\big(aC_0+(b+e)F\big)^2=2a\bigg(b-\Big(\frac{a}{2}-1\Big)e\bigg)\geq 2\sum_{i=1}^r(m_i+1)^2,
     \end{displaymath}
     which is just (\ref{eq:vanishing:1}).
     Similarly, by (\ref{eq:vanishing:2}\alph{subsection}.\roman{subsubsection}) and Lemma
     \ref{lem:ruled-surface} Condition (\ref{eq:vanishing:2}) is
     satisfied. Finally Condition (\ref{eq:vanishing:3}\alph{subsection}.\roman{subsubsection}) implies
     that $D-K_{\Sigma}$ is nef. 

     In order to see the last statement, we
     have to consider two cases.
     \begin{varthm-roman}[Case 1]
       $\ke$ is decomposable, in particular $e\geq 0$.
     \end{varthm-roman}
     If $B\in|a'C_0+b'F|_a$ is irreducible, then we are in one of the
     following situations, by \cite{Har77} V.2.20:
     \begin{enumerate}
        \item $a'=0$ and $b'=1$, which, considering
          (\ref{eq:vanishing:2}\alph{subsection}.\roman{subsubsection}), implies
          \begin{displaymath}
            (D-K_\Sigma).B=a> 0.
          \end{displaymath}
        \item $a'=1$ and $b'=0$, which by (\ref{eq:vanishing:3*})
          leads to
          \begin{displaymath}
            (D-K_\Sigma).B=b-(a-1)e\geq 0.
          \end{displaymath}
        \item $a'>0$ and $b'\geq a'e$, which in view of
          (\ref{eq:vanishing:2}\alph{subsection}.\roman{subsubsection}), 
          (\ref{eq:vanishing:3*}), and $e\geq 0$  gives 
          \begin{displaymath}
            (D-K_\Sigma).B=-aa'e+ab'+(b+e)a'\geq (b+e)a'\geq 0.
          \end{displaymath}
     \end{enumerate}
     Hence, $D-K_\Sigma$ is nef.
     \begin{varthm-roman}[Case 2]
       $\ke$ is indecomposable\tom{, in particular $-2g\leq e\leq 2g-2$}.
     \end{varthm-roman}
     In this case we may apply \cite{Har77} V.2.21 and find that if
     $B\in|a'C_0+b'F|_a$ is irreducible, then we are in one of the 
     following situations:
     \begin{enumerate}
        \item $a'=0$ and $b'=1$, which is
          treated as in Case 1.
        \item $a'=1$ and $b'\geq 0$, which, considering
          (\ref{eq:vanishing:2}\alph{subsection}.\roman{subsubsection}) and
          (\ref{eq:vanishing:3*}), implies
          \begin{displaymath}
            (D-K_\Sigma).B=b-(a-1)e+ab'\geq 0.
          \end{displaymath}
        \item $a'\geq 2$ and $b'\geq \frac{1}{2}a'e$, which leads to
          \begin{displaymath}
            (D-K_\Sigma).B=-aa'e+ab'+(b+e)a'\geq
            -\frac{1}{2}aa'e+(b+e)a'
          \end{displaymath}
          \begin{displaymath}
            =\bigg(b-\Big(\frac{a}{2}-1\Big)e\bigg)a'>
            0, 
          \end{displaymath}
          since by (\ref{eq:vanishing:1}\alph{subsection}.\roman{subsubsection}) and
          (\ref{eq:vanishing:2}\alph{subsection}.\roman{subsubsection}) 
          $\big(b-\big(\frac{a}{2}-1\big)e\big)>0$.  
     \end{enumerate}
     Hence, $D-K_\Sigma$ is nef.     
   \end{proof}
   
   As in the case $\Sigma\cong C\times\PC^1$ we consider for the existence theorems 
   a different algebraic divisor class $D$, namely $D=(a-1)C_0+(b+l+2g-2-e)F$, 
   so that again $D-L-K_{\Sigma}=aC_0+bF$. For Hirzebruch surfaces $D$ 
   again looks somewhat nicer, $D=(a-1)C_0+(b-1)F$.

   \renewcommand{\thesatz}{\ref{cor:existence-II}\alph{subsection}.\roman{subsubsection}}
   \begin{corollary}
     Given integers $a,b\in\Z$ satisfying
     \begin{equationlist}
       \theequationchange{\ref{eq:existence-II:1}\alph{subsection}.\roman{subsubsection}}
       \item $a\big(b-\frac{a}{2}e\big)\geq\sum_{i=1}^r (m_i+1)^2$, 
       \theequationchange{\ref{eq:existence-II:2}\alph{subsection}.\roman{subsubsection}}
       \item $a>\max\{m_i\;|\;i=1,\ldots,r\}$, and 
       \theequationchange{\ref{eq:existence-II:3}\alph{subsection}.\roman{subsubsection}}
       \item $b\geq ae$, if $e>0$,
       \theequationchangeback
     \end{equationlist}
     then for $z_1,\ldots,z_r\in\Sigma$ in very
     general position there is an irreducible reduced curve
     $C\in|(a-1)C_0+(b+l+2g-2-e)F|_a$ with ordinary singularities of
     multiplicities $m_i$ at the $z_i$ as only singularities.
     Moreover, $V_{|C|}(\underline{m})$ is T-smooth at $C$.
   \end{corollary}
   \renewcommand{\thesatz}{\remembersatz}

   \begin{proof}
     \setcounter{sss}{\value{subsubsection}}\addtocounter{subsubsection}{-1}
     Note that by
     (\ref{eq:existence-II:1}\alph{subsection}.\roman{subsubsection}/\roman{sss}) and
     (\ref{eq:existence-II:2}\alph{subsection}.\roman{subsubsection}/\roman{sss}) 
     $b> \frac{a}{2}e\geq ae$, if $e\leq 0$, and thus the
     inequality 
     \begin{equation}
       \label{eq:existence-II:3*}
       b\geq ae,
     \end{equation}
     is fulfilled no matter what $e$ is.

     During the proof we consider the case of an arbitrary
     geometrically ruled surface $\Sigma$. 

     It is clear, that the
     Conditions (\ref{eq:existence-II:1}) and
     (\ref{eq:existence-II:2}) take the form
     (\ref{eq:existence-II:1}\alph{subsection}.\roman{subsubsection}/\roman{sss}) 
     respectively 
     (\ref{eq:existence-II:2}\alph{subsection}.\roman{subsubsection}/\roman{sss}).
     It, therefore, remains to show that (\ref{eq:existence-II:4}) and
     (\ref{eq:existence-II:5}) are obsolete, and that
     (\ref{eq:existence-II:3}) takes the form
     (\ref{eq:existence-II:3}\alph{subsection}.\roman{sss}), which in
     particular means that it is obsolete in the case $\Sigma\cong
     C\times\PC^1$.

     \begin{varthm-roman}[Step 1]
       (\ref{eq:existence-II:4}) is obsolete.
     \end{varthm-roman}
     If $\Sigma\not\cong\PC^1\times\PC^1$, then $l\geq 2$. Since, moreover, $g(L)=g$
     and $D.L=a(l-e)+b+2g-2$, Condition
     (\ref{eq:existence-II:2}\alph{subsection}.\roman{subsubsection}/\roman{sss}) 
     and (\ref{eq:existence-II:3*})
     imply  (\ref{eq:existence-II:4}), i.~e.~for all $i,j$
     \begin{displaymath}
       D.L-2g(L)=a(l-e)+b-2\geq 
       \left\{
         \begin{array}{l}
           a+b-2\geq m_i+m_j, \mbox{ if } \Sigma\cong
           \PC^1\times\PC^1,\\
           2a+(b-ae)-2 \geq m_i+m_j, \mbox{ else.}
         \end{array}
       \right.
     \end{displaymath}

     \begin{varthm-roman}[Step 2]
       (\ref{eq:existence-II:3}) takes the form
       (\ref{eq:existence-II:3}\alph{subsection}.\roman{sss}).
     \end{varthm-roman}
     We have to consider two cases.
     \begin{varthm-roman}[Case 1]
       $\ke$ is decomposable, in particular $e\geq 0$.
     \end{varthm-roman}
     If $B\in|a'C_0+b'F|_a$ is irreducible, then we are in one of the
     following situations, by \cite{Har77} V.2.20:
     \begin{enumerate}
        \item $a'=0$ and $b'=1$, which, considering
          (\ref{eq:existence-II:2}\alph{subsection}.\roman{subsubsection}/\roman{sss}), implies
            \begin{displaymath}
              (D-L-K_\Sigma).B=a> 0.
            \end{displaymath}
          \item $a'=1$ and $b'=0$, which by
            (\ref{eq:existence-II:3*})
            leads to
            \begin{displaymath}
              (D-L-K_\Sigma).B=b-ae\geq 0.
            \end{displaymath}
        \item $a'>0$ and $b'\geq a'e$, which in view of
          (\ref{eq:existence-II:2}\alph{subsection}.\roman{subsubsection}/\roman{sss}), 
          (\ref{eq:existence-II:3*}), and $e\geq 0$ gives 
          \begin{displaymath}
            (D-L-K_\Sigma).B=-aa'e+ab'+ba'\geq ba'\geq 0.
          \end{displaymath}
     \end{enumerate}
     Hence, $D-L-K_\Sigma$ is nef. 
     \begin{varthm-roman}[Case 2]
       $\ke$ is indecomposable\tom{, in particular $-2g\leq e\leq 2g-2$}.
     \end{varthm-roman}
     In this case we may apply \cite{Har77} V.2.21 and find that if
     $B\in|a'C_0+b'F|_a$ is irreducible, then we are in one of the 
     following situations:
     \begin{enumerate}
        \item $a'=0$ and $b'=1$, which is
          treated as in Case 1.
        \item $a'=1$ and $b'\geq 0$, which, considering
          (\ref{eq:existence-II:2}\alph{subsection}.\roman{subsubsection}/\roman{sss}) and
          (\ref{eq:existence-II:3*}), implies
          \begin{displaymath}
            (D-L-K_\Sigma).B=b-ae+ab'\geq 0.
          \end{displaymath}
        \item $a'\geq 2$ and $b'\geq \frac{1}{2}a'e$, which leads to
          \begin{displaymath}
            (D-L-K_\Sigma).B=-aa'e+ab'+ba'\geq
            -\frac{1}{2}aa'e+ba'
            =\bigg(b-\frac{a}{2}e\bigg)a'>
            0, 
          \end{displaymath}
          since by (\ref{eq:existence-II:1}\alph{subsection}.\roman{subsubsection}/\roman{sss}) and
          (\ref{eq:existence-II:2}\alph{subsection}.\roman{subsubsection}/\roman{sss}) 
          $\big(b-\frac{a}{2}e\big)>0$.  
     \end{enumerate}
     Hence, $D-L-K_\Sigma$ is nef.                  

     \begin{varthm-roman}[Step 3]
       (\ref{eq:existence-II:6}) is satisfied, and thus
       (\ref{eq:existence-II:5}) is obsolete. 
     \end{varthm-roman}
     We have
     \begin{displaymath}
       D^2=-e(a-1)^2+2(a-1)(b+l+2g-2-e),
     \end{displaymath}
     and
     \begin{displaymath}
       (2D-L-K_{\Sigma}).(L+K_{\Sigma})= e+2al+4ag+4-2b-4a-4g-2l.
     \end{displaymath}
     Hence Condition (\ref{eq:existence-II:6}) is equivalent to
     \begin{equation}\label{eq:existence-II:6*}
       4b+8a+4l+a^2e+ 8g < 2ab+4al+2e+8ag+ 8+4\sum_{i=1}^rm_i +2r.
     \end{equation}
     If $\Sigma\cong \PC^1\times\PC^1$, then the situation is symmetric and we may
     w.~l.~o.~g.~assume that $b\geq a$. Since by
     (\ref{eq:existence-II:2}\alph{subsection}.\roman{subsubsection}/\roman{sss}) 
     $a\geq 2$ we have to consider the   
     following cases:

     \begin{list}{}{\leftmargin1.3cm \labelwidth1.3cm}
       \item[$a\geq 4$:] 
         \begin{list}{}{}
            \item[$g=0,\; e=0$:] 
              Then $\Sigma\cong\PC^1\times\PC^1$, and by
              assumption $b\geq a\geq 4$ and $l=e+1=1$. We thus have
              $2ab+4al=ab+ab+4a\geq 4b+8a$ and $8> 4l$, which implies 
              (\ref{eq:existence-II:6*}).
            \item[$g=0,\; e>0$ or $g\geq 1,\; e\geq 0$:]
              By (\ref{eq:existence-II:3*}) we get 
              $2ab\geq 4b+ab \geq 4b+a^2e$. 
              \begin{list}{}{}
                 \item[$g=0,\; e\in\{1,2\}$:]
                   Then $l=e+1$, and hence $4al\geq 8a$ and $8+2e\geq 4l$.
                 \item[$g=0,\; e\geq 3$:]
                   Thus $l=e+1\geq 4$, which implies $2al\geq 8a$ and
                   $2al\geq 4l$.
                 \item[$g\geq 1$:] 
                   Then $l\geq 2$, and thus $2al+4ag \geq 
                   8a$, $2al\geq 4l$, and $4ag\geq  8g$.
              \end{list}
              In any of the above cases (\ref{eq:existence-II:6*}) is
              satisfied. 
            \item[$g\geq 1,\; e<0$:] 
              Then $l\geq 2$ and  $2ag\geq  8g$. We therefore
              consider the following cases:
              \begin{list}{}{}
                 \item[$b\geq 0$:]
                   Thus $2ab\geq 4b$, $2al+4ag\geq 8a$, $2al\geq
                   4l$, and $2e\geq a^2e$.
                 \item[$b<0$:]
                   By
                   (\ref{eq:existence-II:1}\alph{subsection}.\roman{subsubsection}/\roman{sss}) and 
                   (\ref{eq:existence-II:2}\alph{subsection}.\roman{subsubsection}/\roman{sss}) 
                   $2ab\geq a^2e$, and of course $0>4b$. Moreover, since $e\geq -g$, we have
                   $ag+2e\geq 0$. And finally, $3al+5ag\geq 8a$ and
                   $al\geq 4l$.
              \end{list}
              These considerations together ensure that
              (\ref{eq:existence-II:6*}) is fulfilled.                 
         \end{list}
       \item[$a=3$:] 
         In this case (\ref{eq:existence-II:6*}) comes down to 
         \tom{
         \begin{displaymath}
           4b+24+4l+9e+ 8g<6b+12l+2e+24g+ 8+4\sum_{i=1}^rm_i +2r.           
         \end{displaymath}
         or, equivalently,}
         \begin{equation}\label{eq:existence-II:6**}
           16+7e<2b+8l+16g+4\sum_{i=1}^rm_i +2r.
         \end{equation}
         \begin{list}{}{}
            \item[$e>0$, or $e=0$ and $g=0$:]
              Then $b\geq a=3$. Thus $2b+8l+4\sum_{i=1}^rm_i\geq
              6+8(e+1)+4>16+7e$, so that the inequality 
              (\ref{eq:existence-II:6**}) is 
              certainly satisfied.
            \item[$e<0$, or $e=0$ and $g\geq 1$:]
              Then $g\geq 1$ and $16g\geq 16+7e$, so that again the inequality
              (\ref{eq:existence-II:6**}) is fulfilled.
         \end{list}
       \item[$a=2$:] 
         (\ref{eq:existence-II:6*}) reads just
         \tom{
         \begin{displaymath}
           4b+16+4l+4e+8g<4b+8l+2e+16g+8+4\sum_{i=1}^rm_i +2r,
         \end{displaymath}
         or, equivalently,}
         \begin{equation}
           \label{eq:existence-II:6***}
           8+2e<4l+8g+4\sum_{i=1}^rm_i +2r.
         \end{equation}
         \begin{list}{}{}
            \item[$e<0$:]
              Then $g\geq 1$, and thus $8g\geq 8+2e$, which implies
              (\ref{eq:existence-II:6***}). 
            \item[$e\geq 0$:] 
              Then $4l+4\sum_{i=1}^rm_i\geq 4(e+1)+4\geq
              8+2e$, and hence (\ref{eq:existence-II:6***}) is
              fulfilled. 
         \end{list}         
    \end{list}
    \stepcounter{subsubsection}
  \end{proof}

   With the same $D$ and $L$ as above the conditions in the existence
   theorem Corollary \ref{cor:existence-IV} reduce to
   \begin{equationlist}
     \theequationchange{\ref{eq:existence-IV:1}\alph{subsection}.\roman{subsubsection}}
     \item $a(b-\frac{a}{2}e)
         \geq \frac{207}{5}\sum\limits_{\mu(\ks_i)\leq 38} \mu(\ks_i) +
         29\sum\limits_{\mu(\ks_i)\geq 39} \Big(\sqrt{\mu(\ks_i)}+\frac{13}{2\sqrt{29}}\Big)^2$,
     \theequationchange{\ref{eq:existence-IV:2}\alph{subsection}.\roman{subsubsection}}
     \item $a >
         \left\{
           \begin{array}{ll}
             \sqrt{\frac{207}{5}}\sqrt{\mu(\ks_1)}-1, & \mbox{ if 
               }\mu(\ks_1)\leq 38,\\[0.3cm]
             \sqrt{29} \sqrt{\mu(\ks_1)}+\frac{11}{2}, & \mbox{
               if }\mu(\ks_1)\geq 39,\mbox{ and,} 
           \end{array}
         \right.$
     \theequationchange{\ref{eq:existence-IV:3}\alph{subsection}.\roman{subsubsection}}
     \item $b\geq ae$, if $e>0$.
     \theequationchangeback
   \end{equationlist}

}

\problem{2}{

   \renewcommand{\thesatz}{\ref{cor:vanishing}\alph{subsection}}
   \begin{corollary}
     Given two integers $a,b\in\Z$ satisfying 
     \begin{equationlist}
       \theequationchange{\ref{eq:vanishing:1}\alph{subsection}}
       \item $a\big(b-\big(\frac{a}{2}-1\big)e\big)\geq \sum_{i=1}^r(m_i+1)^2$, 
       \theequationchange{\ref{eq:vanishing:2}\alph{subsection}.i}
       \item $a>\max\{m_i\;|\;i=1,\ldots,r\}$, 
       \theequationchange{\ref{eq:vanishing:2}\alph{subsection}.ii}
       \item 
         $b>\max\{m_i\;|\;i=1,\ldots,r\}$, if
         $e=0$, 
       \theequationchange{\ref{eq:vanishing:2}\alph{subsection}.iii}
       \item 
         $2\big(b-\big(\frac{a}{2}-1\big)e\big)>\max\{m_i\;|\;i=1,\ldots,r\}$, if
         $e< 0$, and
       \theequationchange{\ref{eq:vanishing:3}\alph{subsection}}
       \item $b\geq (a-1)e$, if $e>0$.
       \theequationchangeback
     \end{equationlist}
     For $z_1,\ldots,z_r\in\Sigma$ in very general position and
     $\nu>0$ 
     \begin{displaymath}
       H^\nu\left(\Bl_{\underline{z}}(\Sigma),(a-2)\pi^*C_0+(b-2+2g)\pi^*F-\sum_{i=1}^r m_iE_i\right)=0.
     \end{displaymath}
   \end{corollary}
   \renewcommand{\thesatz}{\remembersatz}

   \begin{proof}
     Note that if the invariant $e$ is non-positive, then
     $\big(b-\big(\frac{a}{2}-1\big)e\big)>0$ implies 
     \begin{equation} 
       \label{eq:vanishing:3*}
       b\geq(a-1)e,
     \end{equation}
     so that this inequality is fulfilled for any choice of $e$.

     Setting $D=(a-2)C_0+(b-2+2g)F$
     we have
     \begin{displaymath}
       (D-K_{\Sigma})^2=\big(aC_0+(b+e)F\big)^2=2a\bigg(b-\Big(\frac{a}{2}-1\Big)e\bigg)\geq 2\sum_{i=1}^r(m_i+1)^2,
     \end{displaymath}
     which is just (\ref{eq:vanishing:1}).
     Similarly, by (\ref{eq:vanishing:2}\alph{subsection}) and Lemma
     \ref{lem:ruled-surface} Condition (\ref{eq:vanishing:2}.i/ii/iii) is
     satisfied.\footnote{To see this, let $B\sim_a a'C_0+b'F$ be an
       irreducible curve with $B^2=0$. Then by Lemma
       \ref{lem:ruled-surface} either $a'=0$ and $b'=1$, or $e=0$,
       $a'\geq 1$ and $b'=0$, or $e<0$,
       $a'\geq 2$, and $b'=\tfrac{a'}{2}e<0$. In the first case,
       $(D-K_{\Sigma}).B=a>\max\{m_i\;|\;i=1,\ldots,r\}$ by
       (\ref{eq:vanishing:2}\alph{subsection}.i). In the second case,
       $(D-K_{\Sigma}).B=ba'\geq b>\max\{m_i\;|\;i=1,\ldots,r\}$ by
       (\ref{eq:vanishing:2}\alph{subsection}.ii). 
       And finally, in the third case, we have
       $(D-K_{\Sigma}).B=a'\cdot\big(b-\big(\frac{a}{2}-1\big)e)\big) 
       >\max\{m_i\;|\;i=1,\ldots,r\}$ by
       (\ref{eq:vanishing:2}\alph{subsection}.iii).}
     Finally Condition (\ref{eq:vanishing:3}\alph{subsection}) implies
     that $D-K_{\Sigma}$ is nef. 

     In order to see the last statement, we
     have to consider two cases.
     \begin{varthm-roman}[Case 1]
       $e\geq 0$.
     \end{varthm-roman}
     If $B\in|a'C_0+b'F|_a$ is irreducible, then we are in one of the
     following situations, by \cite{Har77} V.2.20:
     \begin{enumerate}
        \item $a'=0$ and $b'=1$, which, considering
          (\ref{eq:vanishing:2}\alph{subsection}.i), implies
          \begin{displaymath}
            (D-K_\Sigma).B=a> 0.
          \end{displaymath}
        \item $a'=1$ and $b'=0$, which by (\ref{eq:vanishing:3*})
          leads to
          \begin{displaymath}
            (D-K_\Sigma).B=b-(a-1)e\geq 0.
          \end{displaymath}
        \item $a'>0$ and $b'\geq a'e$, which in view of
          (\ref{eq:vanishing:2}\alph{subsection}.i), 
          (\ref{eq:vanishing:3*}), and $e\geq 0$  gives 
          \begin{displaymath}
            (D-K_\Sigma).B=-aa'e+ab'+(b+e)a'\geq (b+e)a'\geq 0.
          \end{displaymath}
     \end{enumerate}
     Hence, $D-K_\Sigma$ is nef.
     \begin{varthm-roman}[Case 2]
       $e<0$.
     \end{varthm-roman}
     In this case we may apply \cite{Har77} V.2.21 and find that if
     $B\in|a'C_0+b'F|_a$ is irreducible, then we are in one of the 
     following situations:
     \begin{enumerate}
        \item $a'=0$ and $b'=1$, which is
          treated as in Case 1.
        \item $a'=1$ and $b'\geq 0$, which, considering
          (\ref{eq:vanishing:2}\alph{subsection}.i) and
          (\ref{eq:vanishing:3*}), implies
          \begin{displaymath}
            (D-K_\Sigma).B=b-(a-1)e+ab'\geq 0.
          \end{displaymath}
        \item $a'\geq 2$ and $b'\geq \frac{1}{2}a'e$, which in view of
          (\ref{eq:vanishing:2}\alph{subsection}.iii)
          leads to
          \begin{displaymath}
            (D-K_\Sigma).B=-aa'e+ab'+(b+e)a'\geq
            -\frac{1}{2}aa'e+(b+e)a'
          \end{displaymath}
          \begin{displaymath}
            =\bigg(b-\Big(\frac{a}{2}-1\Big)e\bigg)a'>
            0.
          \end{displaymath}
     \end{enumerate}
     Hence, $D-K_\Sigma$ is nef.     
   \end{proof}
   
   In order to obtain nice formulae we considered $D=(a-2)C_0+(b-2+2g)F$
   in the formulation of the vanishing theorem. For the existence
   theorems it turns out that the formulae look best if we work with
   $D=(a-1)C_0+(b+l+2g-2-e)F$ instead. In the case of Hirzebruch
   surfaces this is just $D=(a-1)C_0+(b-1)F$.

   \renewcommand{\thesatz}{\ref{cor:existence-II}\alph{subsection}}
   \begin{corollary}
     Given integers $a,b\in\Z$ satisfying
     \begin{equationlist}
       \theequationchange{\ref{eq:existence-II:1}\alph{subsection}}
       \item $a\big(b-\frac{a}{2}e\big)\geq\sum_{i=1}^r (m_i+1)^2$, 
       \theequationchange{\ref{eq:existence-II:2}\alph{subsection}.i}
       \item $a>\max\{m_i\;|\;i=1,\ldots,r\}$, 
       \theequationchange{\ref{eq:existence-II:2}\alph{subsection}.ii}
       \item 
         $b>\max\{m_i\;|\;i=1,\ldots,r\}$, if
         $e=0$, 
       \theequationchange{\ref{eq:existence-II:2}\alph{subsection}.iii}
       \item 
         $2\big(b-\frac{a}{2}e\big)>\max\{m_i\;|\;i=1,\ldots,r\}$, if
         $e< 0$, and
       \theequationchange{\ref{eq:existence-II:3}\alph{subsection}}
       \item $b\geq ae$, if $e>0$,
       \theequationchangeback
     \end{equationlist}
     then for $z_1,\ldots,z_r\in\Sigma$ in very
     general position there is an irreducible reduced curve
     $C\in|(a-1)C_0+(b+l+2g-2-e)F|_a$ with ordinary singularities of
     multiplicities $m_i$ at the $z_i$ as only singularities.
     Moreover, $V_{|C|}(\underline{m})$ is T-smooth at $C$.
   \end{corollary}
   \renewcommand{\thesatz}{\remembersatz}

   \begin{proof}
     Note that by
     (\ref{eq:existence-II:1}\alph{subsection}) and
     (\ref{eq:existence-II:2}\alph{subsection}.i) 
     $b> \frac{a}{2}e\geq ae$, if $e\leq 0$, and thus the
     inequality 
     \begin{equation}
       \label{eq:existence-II:3*}
       b\geq ae,
     \end{equation}
     is fulfilled no matter what $e$ is.


     Noting that $D-L-K_\Sigma\sim_a aC_0+bF$, it is in view of Lemma
     \ref{lem:ruled-surface} clear, that the
     Conditions (\ref{eq:existence-II:1}) and
     (\ref{eq:existence-II:2}) take the form
     (\ref{eq:existence-II:1}\alph{subsection}) 
     respectively 
     (\ref{eq:existence-II:2}\alph{subsection}).
     It, therefore, remains to show that (\ref{eq:existence-II:4}) and
     (\ref{eq:existence-II:5}) are obsolete, and that
     (\ref{eq:existence-II:3}) takes the form
     (\ref{eq:existence-II:3}\alph{subsection}), which in
     particular means that it is obsolete in the case $\Sigma\cong
     C\times\PC^1$.

     \begin{varthm-roman}[Step 1]
       (\ref{eq:existence-II:4}) is obsolete.
     \end{varthm-roman}
     If $\Sigma\not\cong\PC^1\times\PC^1$, then $l\geq 2$. Since, moreover, $g(L)=g$
     and $D.L=a(l-e)+b+2g-2$, Condition
     (\ref{eq:existence-II:2}\alph{subsection}.i) 
     and (\ref{eq:existence-II:3*})
     imply  (\ref{eq:existence-II:4}), i.~e.~for all $i,j$
     \begin{displaymath}
       D.L-2g(L)=a(l-e)+b-2\geq 
       \left\{
         \begin{array}{l}
           a+b-2\geq m_i+m_j, \mbox{ if } \Sigma\cong
           \PC^1\times\PC^1,\\
           2a+(b-ae)-2 \geq m_i+m_j, \mbox{ else.}
         \end{array}
       \right.
     \end{displaymath}

     \begin{varthm-roman}[Step 2]
       (\ref{eq:existence-II:3}) takes the form
       (\ref{eq:existence-II:3}\alph{subsection}).
     \end{varthm-roman}
     We have to consider two cases.
     \begin{varthm-roman}[Case 1]
       $e\geq 0$.
     \end{varthm-roman}
     If $B\in|a'C_0+b'F|_a$ is irreducible, then we are in one of the
     following situations, by \cite{Har77} V.2.20:
     \begin{enumerate}
        \item $a'=0$ and $b'=1$, which, considering
          (\ref{eq:existence-II:2}\alph{subsection}.i), implies
            \begin{displaymath}
              (D-L-K_\Sigma).B=a> 0.
            \end{displaymath}
          \item $a'=1$ and $b'=0$, which by
            (\ref{eq:existence-II:3*})
            leads to
            \begin{displaymath}
              (D-L-K_\Sigma).B=b-ae\geq 0.
            \end{displaymath}
        \item $a'>0$ and $b'\geq a'e$, which in view of
          (\ref{eq:existence-II:2}\alph{subsection}.i), 
          (\ref{eq:existence-II:3*}), and $e\geq 0$ gives 
          \begin{displaymath}
            (D-L-K_\Sigma).B=-aa'e+ab'+ba'\geq ba'\geq 0.
          \end{displaymath}
     \end{enumerate}
     Hence, $D-L-K_\Sigma$ is nef. 
     \begin{varthm-roman}[Case 2]
       $e<0$.
     \end{varthm-roman}
     In this case we may apply \cite{Har77} V.2.21 and find that if
     $B\in|a'C_0+b'F|_a$ is irreducible, then we are in one of the 
     following situations:
     \begin{enumerate}
        \item $a'=0$ and $b'=1$, which is
          treated as in Case 1.
        \item $a'=1$ and $b'\geq 0$, which, considering
          (\ref{eq:existence-II:2}\alph{subsection}.i) and
          (\ref{eq:existence-II:3*}), implies
          \begin{displaymath}
            (D-L-K_\Sigma).B=b-ae+ab'\geq 0.
          \end{displaymath}
        \item $a'\geq 2$ and $b'\geq \frac{1}{2}a'e$, which in view of
          (\ref{eq:existence-II:2}\alph{subsection}.iii) leads to
          \begin{displaymath}
            (D-L-K_\Sigma).B=-aa'e+ab'+ba'\geq
            -\frac{1}{2}aa'e+ba'
            =\bigg(b-\frac{a}{2}e\bigg)a'>
            0.
          \end{displaymath}
     \end{enumerate}
     Hence, $D-L-K_\Sigma$ is nef.                  

     \begin{varthm-roman}[Step 3]
       (\ref{eq:existence-II:6}) is satisfied, and thus
       (\ref{eq:existence-II:5}) is obsolete. 
     \end{varthm-roman}
     We have
     \begin{displaymath}
       D^2=-e(a-1)^2+2(a-1)(b+l+2g-2-e),
     \end{displaymath}
     and
     \begin{displaymath}
       (2D-L-K_{\Sigma}).(L+K_{\Sigma})= e+2al+4ag+4-2b-4a-4g-2l.
     \end{displaymath}
     Hence Condition (\ref{eq:existence-II:6}) is equivalent to
     \begin{equation}\label{eq:existence-II:6*}
       4b+8a+4l+a^2e+ 8g < 2ab+4al+2e+8ag+ 8+4\sum_{i=1}^rm_i +2r.
     \end{equation}
     If $\Sigma\cong \PC^1\times\PC^1$, then the situation is symmetric and we may
     w.~l.~o.~g.~assume that $b\geq a$. Since by
     (\ref{eq:existence-II:2}\alph{subsection}.i) 
     $a\geq 2$ we have to consider the   
     following cases:

     \begin{list}{}{\leftmargin1.3cm \labelwidth1.3cm}
       \item[$a\geq 4$:] 
         \begin{list}{}{}
            \item[$g=0,\; e=0$:] 
              Then $\Sigma\cong\PC^1\times\PC^1$, and by
              assumption $b\geq a\geq 4$ and $l=e+1=1$. We thus have
              $2ab+4al=ab+ab+4a\geq 4b+8a$ and $8> 4l$, which implies 
              (\ref{eq:existence-II:6*}).
            \item[$g=0,\; e>0$ or $g\geq 1,\; e\geq 0$:]
              By (\ref{eq:existence-II:3*}) we get 
              $2ab\geq 4b+ab \geq 4b+a^2e$. 
              \begin{list}{}{}
                 \item[$g=0,\; e\in\{1,2\}$:]
                   Then $l=e+1$, and hence $4al\geq 8a$ and $8+2e\geq 4l$.
                 \item[$g=0,\; e\geq 3$:]
                   Thus $l=e+1\geq 4$, which implies $2al\geq 8a$ and
                   $2al\geq 4l$.
                 \item[$g\geq 1$:] 
                   Then $l\geq 2$, and thus $2al+4ag \geq 
                   8a$, $2al\geq 4l$, and $4ag\geq  8g$.
              \end{list}
              In any of the above cases (\ref{eq:existence-II:6*}) is
              satisfied. 
            \item[$g\geq 1,\; e<0$:] 
              Then $l\geq 2$ and  $2ag\geq  8g$. We therefore
              consider the following cases:
              \begin{list}{}{}
                 \item[$b\geq 0$:]
                   Thus $2ab\geq 4b$, $2al+4ag\geq 8a$, $2al\geq
                   4l$, and $2e\geq a^2e$.
                 \item[$b<0$:]
                   By
                   (\ref{eq:existence-II:2}\alph{subsection}.i) and 
                   (\ref{eq:existence-II:2}\alph{subsection}.iii) 
                   $2ab\geq a^2e$, and of course $0>4b$. Moreover, since $e\geq -g$, we have
                   $ag+2e\geq 0$. And finally, $3al+5ag\geq 8a$ and
                   $al\geq 4l$.
              \end{list}
              These considerations together ensure that
              (\ref{eq:existence-II:6*}) is fulfilled.                 
         \end{list}
       \item[$a=3$:] 
         In this case (\ref{eq:existence-II:6*}) comes down to 
         \tom{
         \begin{displaymath}
           4b+24+4l+9e+ 8g<6b+12l+2e+24g+ 8+4\sum_{i=1}^rm_i +2r.           
         \end{displaymath}
         or, equivalently,}
         \begin{equation}\label{eq:existence-II:6**}
           16+7e<2b+8l+16g+4\sum_{i=1}^rm_i +2r.
         \end{equation}
         \begin{list}{}{}
            \item[$e>0$, or $e=0$ and $g=0$:]
              Then $b\geq a=3$. Thus $2b+8l+4\sum_{i=1}^rm_i\geq
              6+8(e+1)+4>16+7e$, so that the inequality 
              (\ref{eq:existence-II:6**}) is 
              certainly satisfied.
            \item[$e<0$, or $e=0$ and $g\geq 1$:]
              Then $g\geq 1$ and $16g\geq 16+7e$, so that again the inequality
              (\ref{eq:existence-II:6**}) is fulfilled.
         \end{list}
       \item[$a=2$:] 
         (\ref{eq:existence-II:6*}) reads just
         \tom{
         \begin{displaymath}
           4b+16+4l+4e+8g<4b+8l+2e+16g+8+4\sum_{i=1}^rm_i +2r,
         \end{displaymath}
         or, equivalently,}
         \begin{equation}
           \label{eq:existence-II:6***}
           8+2e<4l+8g+4\sum_{i=1}^rm_i +2r.
         \end{equation}
         \begin{list}{}{}
            \item[$e<0$:]
              Then $g\geq 1$, and thus $8g\geq 8+2e$, which implies
              (\ref{eq:existence-II:6***}). 
            \item[$e\geq 0$:] 
              Then $4l+4\sum_{i=1}^rm_i\geq 4(e+1)+4\geq
              8+2e$, and hence (\ref{eq:existence-II:6***}) is
              fulfilled. 
         \end{list}         
    \end{list}
  \end{proof}

   With the same $D$ and $L$ as above the conditions in the existence
   theorem Corollary \ref{cor:existence-IV} reduce to
   \begin{equationlist}
     \theequationchange{\ref{eq:existence-IV:1}\alph{subsection}}
     \item $a(b-\frac{a}{2}e)
         \geq \frac{207}{5}\sum\limits_{\mu(\ks_i)\leq 38} \mu(\ks_i) +
         29\sum\limits_{\mu(\ks_i)\geq 39} \Big(\sqrt{\mu(\ks_i)}+\frac{13}{2\sqrt{29}}\Big)^2$,
     \theequationchange{\ref{eq:existence-IV:2}\alph{subsection}.i}
     \item $a >
         \left\{
           \begin{array}{ll}
             \sqrt{\frac{207}{5}}\sqrt{\mu(\ks_1)}-1, & \mbox{ if 
               }\mu(\ks_1)\leq 38,\\[0.3cm]
             \sqrt{29} \sqrt{\mu(\ks_1)}+\frac{11}{2}, & \mbox{
               if }\mu(\ks_1)\geq 39,
           \end{array}
         \right.$
       \theequationchange{\ref{eq:existence-IV:2}\alph{subsection}.ii}
       \item 
         $b>
         \left\{
           \begin{array}{ll}
             \sqrt{\frac{207}{5}}\sqrt{\mu(\ks_1)}-1, & \mbox{ if 
               }\mu(\ks_1)\leq 38,\\[0.3cm]
             \sqrt{29} \sqrt{\mu(\ks_1)}+\frac{11}{2}, & \mbox{
               if }\mu(\ks_1)\geq 39, 
           \end{array}
         \right.$
         if $e=0$, 
       \theequationchange{\ref{eq:existence-IV:2}\alph{subsection}.iii}
       \item 
         $2\big(b-\frac{a}{2}e\big)>
         \left\{
           \begin{array}{ll}
             \sqrt{\frac{207}{5}}\sqrt{\mu(\ks_1)}-1, & \mbox{ if 
               }\mu(\ks_1)\leq 38,\\[0.3cm]
             \sqrt{29} \sqrt{\mu(\ks_1)}+\frac{11}{2}, & \mbox{
               if }\mu(\ks_1)\geq 39, 
           \end{array}
         \right.$
         if $e< 0$, and 
     \theequationchange{\ref{eq:existence-IV:3}\alph{subsection}}
     \item $b\geq ae$, if $e>0$.
     \theequationchangeback
   \end{equationlist}

}


   \subsection{Products of Curves}\label{subsec:product-curves}

   Let $C_1$ and $C_2$ be two smooth projective curves of genuses $g_1\geq 1$ and
   $g_2\geq 1$ respectively. The surface $\Sigma=C_1\times
   C_2$ is naturally equipped with two fibrations
   $\pr_i:\Sigma\rightarrow C_i$, $i=1,2$, and by 
   abuse of notation we denote two generic fibres $\pr_2^{-1}(p_2)=C_1\times\{p_2\}$
   resp.~$\pr_1^{-1}(p_1)=\{p_1\}\times C_2$ again by $C_1$ resp.~$C_2$.

   One can show that for a generic choice of the curves
   $C_1$ and $C_2$ the  Neron-Severi group $\NS(\Sigma)=C_1\Z\oplus
   C_2\Z$ of $\Sigma$ is two-dimensional\footnote{In the case that
     $C_1$ and $C_2$ are elliptic curves, generic means precisely,
     that they are not isogenous - see Section
     \ref{subsec:elliptic-curves}. For a further investigation of the
     Neron-Severi group of a product of two curves we refer to
     Appendix \ref{app:product-surfaces}.}
   with intersection matrix 
   \begin{displaymath}
     (C_i.C_j)_{i,j}=\left(
     \begin{array}[m]{cc}
       0&1\\1&0
     \end{array}
     \right).
   \end{displaymath}
   Thus, the only irreducible curves $B\subset \Sigma$ with
   selfintersection $B^2=0$ are the fibres $C_1$ and $C_2$, and for
   any irreducible curve $B\sim_a aC_1+bC_2$ the coefficients $a$ and
   $b$ must be non-negative. Taking
   into account that $K_\Sigma\sim_a (2g_2-2)C_1+(2g_1-2)C_2$ Corollary
   \ref{cor:vanishing} comes down to the following.

   \renewcommand{\thesatz}{\ref{cor:vanishing}\alph{subsection}}
   \begin{corollary}
     Let $C_1$ and $C_2$ be two generic curves with $g(C_i)=g_i\geq
     1$, $i=1,2$, and let $a,b\in\Z$ be
     integers satisfying 
     \begin{equationlist}
       \theequationchange{\ref{eq:vanishing:1}\alph{subsection}}
       \item $(a-2g_2+2)(b-2g_1+2)\geq \sum_{i=1}^r(m_i+1)^2$, and
       \theequationchange{\ref{eq:vanishing:2}\alph{subsection}}
       \item $(a-2g_2+2),(b-2g_1+2)>\max\{m_i\;|\;i=1,\ldots,r\}$,
       \theequationchangeback
     \end{equationlist}
     then for $z_1,\ldots,z_r\in\Sigma= C_1\times C_2$ in very general
     position and $\nu>0$
     \begin{displaymath}
       H^\nu\left(\Bl_{\underline{z}}(\Sigma),\;
         a\pi^*C_1+b\pi^*C_2-\sum_{i=1}^r m_iE_i\right)=0. 
     \end{displaymath}
   \end{corollary}
   \renewcommand{\thesatz}{\remembersatz}

   We know that $C_1+C_2$ has positive self-intersection and
   intersects any irreducible curve positive, is thus ample by
   Nakai-Moishezon. But then we may find some integer $l\geq 3$ such
   that $L=lC_1+lC_2$ is very ample. We choose $l$ minimal with this
   property for the existence theorem Corollary
   \ref{cor:existence-II}, 
   and we claim
   that the Conditions (\ref{eq:existence-II:3}), 
   (\ref{eq:existence-II:4}) and (\ref{eq:existence-II:5}) become
   obsolete, while
   (\ref{eq:existence-II:1}) and 
   (\ref{eq:existence-II:2}) take the form
   \begin{equationlist}
     \theequationchange{\ref{eq:existence-II:1}\alph{subsection}}
     \item
       $(a-l-2g_2+2)(b-l-2g_1+2) \geq \sum_{i=1}^r (m_i+1)^2$, and
     \theequationchange{\ref{eq:existence-II:2}\alph{subsection}}
     \item
       $(a-l-2g_2+2),(b-l-2g_1+2)>\max\{m_i\;|\;i=1,\ldots,r\}$.
   \end{equationlist}
   That is, under these hypotheses there is an irreducible curve
   in $|D|_l$, for any $D\sim_a aC_1+bC_2$, with precisely $r$
   ordinary singular  points of multiplicities $m_1,\ldots,m_r$.
   
   (\ref{eq:existence-II:3}) becomes redundant in view of
   (\ref{eq:existence-II:2}\alph{subsection}) and since an irreducible curve
   $B\sim_a a'C_1+b'C_2$ has non-negative coefficients $a'$ and $b'$.
   For (\ref{eq:existence-II:5}) we look at (\ref{eq:existence-II:6}),
   which in this case takes the form
   \begin{displaymath}
     2ab+\big((2a-l-2g_2+2)(l+2g_1-2)+(2b-l-2g_1+2)(l+2g_2-2)\big)
   \end{displaymath}
   \begin{displaymath}
     +4\sum_{i=1}^rm_i+2r>0.
   \end{displaymath}
   However, in view of (\ref{eq:existence-II:2}\alph{subsection}) the
   factors and summands on the left-hand side are all positive, so
   that the inequality is fulfilled.

   It remains to show that $D.L-g(L)\geq m_i+m_j$ for
   all $i,j$. However, by the adjunction formula
   $g(L)=1+\frac{1}{2}(L^2+L.K_\Sigma)=1+l\cdot(l+g_1+g_2-2)$, and by
   (\ref{eq:existence-II:2}\alph{subsection}) $D.L-g(L)>
   l\cdot\big((a-l-2g_2+2)+(b-l-2g_1+2)\big)>3(m_i+m_j)\geq
   m_i+m_j$. Thus the claim is proved. 

   From these considerations we at once deduce the conditions for the
   existence of an irreducible curve in $|D|_l$, 
   $D\sim_aaC_1+bC_2$, with prescribed
   singularities of arbitrary type, i.~e.~the conditions in Corollary
   \ref{cor:existence-IV}. They come down to 

   \begin{equationlist}
     \theequationchange{\ref{eq:existence-IV:1}\alph{subsection}}
   \item \hspace*{-0.5cm}$(a-l-2g_2+2)(b-l-2g_1+2) \geq
     \frac{207}{5}\!\!\!\!\!\!\sum\limits_{\mu(\ks_i)\leq 38}\!\!\!\! \mu(\ks_i) +
     29\!\!\!\!\!\!\sum\limits_{\mu(\ks_i)\geq 39}
     \Big(\sqrt{\mu(\ks_i)}+\frac{13}{2\sqrt{29}}\Big)^2$,
   \end{equationlist}
   and
   \begin{equationlist}
     \theequationchange{\ref{eq:existence-IV:2}\alph{subsection}}
   \item \hspace*{-0.5cm}$(a-l-2g_2+2),(b-l-2g_1+2) > \left\{
           \begin{array}{ll}
             \sqrt{\frac{207}{5}}\sqrt{\mu(\ks_1)}-1, & \mbox{ if 
               }\mu(\ks_1)\leq 38,\\[0.3cm]
             \sqrt{29} \sqrt{\mu(\ks_1)}+\frac{11}{2}, & \mbox{
               if }\mu(\ks_1)\geq 39.
           \end{array}
         \right.$ \theequationchangeback
   \end{equationlist}


   \subsection{Products of Elliptic Curves}\label{subsec:elliptic-curves}

   Let $C_1=\C/\Lambda_1$ and $C_2=\C/\Lambda_2$ be two elliptic 
   curves, where $\Lambda_i=\Z\oplus\tau_i\Z\subset\C$ is a lattice 
   and $\tau_i$ is in the upper half plane of $\C$. We denote the
   natural group structure on each of the $C_i$ by $+$ and the neutral 
   element by $0$. 

   Our interest lies in the study of the surface $\Sigma=C_1\times
   C_2$. This surface is naturally equipped with two fibrations
   $\pr_i:\Sigma\rightarrow C_i$, $i=1,2$, and by 
   abuse of notation we denote the fibres $\pr_2^{-1}(0)=C_1\times\{0\}$
   resp.~$\pr_1^{-1}(0)=\{0\}\times C_2$ again by $C_1$ resp.~$C_2$. The group
   structures on $C_1$ and $C_2$ extend to $\Sigma$ so that $\Sigma$
   itself is an abelian variety. Moreover, for $p=(p_1,p_2)\in \Sigma$ the mapping
   $\tau_p:\Sigma\rightarrow\Sigma:(a,b)\mapsto(a+p_1,b+p_2)$ is an
   automorphism of abelian varieties. Due
   to these translation morphisms we know that for any 
   curve $B\subset \Sigma$ the algebraic family of curves $|B|_a$
   covers the whole of $\Sigma$, and in particular $\dim|B|_a\geq
   1$. This also implies $B^2\geq 0$.

   Since $\Sigma$ is an abelian surface, $\NS(\Sigma)=\Num(\Sigma)$, $K_\Sigma=0$,
   and the Picard number $\rho=\rho(\Sigma)\leq 4$ (cf.~\cite{LB92}
   4.11.2 and Ex.~2.5). But the N\'eron-Severi group of $\Sigma$ contains the
   two independent elements $C_1$ and $C_2$, so that $\rho\geq 
   2$. The general case\footnote{The abelian surfaces with $\rho\geq 2$ 
     possessing a principle polarisation are parametrized by
     a countable number of surfaces in a three-dimensional space, and
     the Picard number of such an abelian surface is two unless it is
     contained the intersection of two or three of these surfaces
     (cf.~\cite{IS96} 11.2). See also \cite{GH94} p.~286 and
     Proposition \ref{prop:product-surfaces}.}
   is indeed $\rho=2$, however $\rho$ might also be larger (see Example
   \ref{ex:product}), in which case the additional generators may be
   chosen to be graphs of surjective morphisms from $C_1$ to $C_2$
   (cf.~\cite{IS96} 3.2 Example 3). That is, $\rho(\Sigma)=2$  if
   and only if $C_1$ and $C_2$ are not isogenous.

   \begin{lemma}\label{lem:product}
     Let $B\subset\Sigma$ be an irreducible curve, $B\not\sim_a C_k$,
     $k=1,2$, and $\{i,j\}=\{1,2\}$.
     \begin{enumerate}
       \item If $B^2=0$, then $B$ is smooth, $g(B)=1$, and
         $\pr_{i|}:B\rightarrow C_i$ is an unramified covering of
         degree $B.C_j$. 
       \item If $B^2=0$, then $\#\big(B\cap
         \tau_p(C_i)\big) =B.C_j$ for any $p\in\Sigma$, and if
         $q,q'\in B$, then $\tau_{q-q'}(B)=B$.
       \item If $B^2=0$, then the base curve $H$ in the fibration
         $\pi:\Sigma\rightarrow H$ with fibre $B$, which exists
         according to Proposition \ref{prop:condition}, is an elliptic
         curve.  
       \item If $B.C_i=1$, then $B^2=0$ and $C_j\cong B$.
       \item If $B.C_i=1=B.C_j$, then $C_1\cong C_2$.
       \item If $B$ is the graph of a morphism $\alpha:C_i\rightarrow
         C_j$, then $B.C_j=1$ and $B^2=0$.
     \end{enumerate}
   \end{lemma}
   \begin{proof}\leererpunkt
     \begin{enumerate}
       \item The adjunction formula gives
         \begin{displaymath}
           p_a(B)=1+\frac{B^2+K_\Sigma.B}{2}=1.
         \end{displaymath}
         \\
         Since $|C_2|_a$ covers the whole of $\Sigma$ 
         and $B\not\sim_a C_2$, the two irreducible curves $B$ and
         $C_2$ must intersect properly, that is, $B$ is not a fibre of 
         $\pr_1$. But then
         the mapping $\pr_{1|}:B\rightarrow C_1$ is a 
         finite surjective morphism of degree $B.C_2$. If $B$ was a
         singular curve its normalisation would have to have
         arithmetical genus $0$ and the composition of the
         normalisation with $\pr_{1|}$ would give rise to a surjective
         morphism from $\PC^1$ to an elliptic curve, contradicting
         Hurwitz's formula. Hence, $B$ is smooth and
         $g(B)=p_a(B)=1$. We thus may
         apply the formula of Hurwitz to $\pr_{1|}$ and the degree of the 
         ramification divisor $R$ turns out to be
         \begin{displaymath}
           \deg(R)= 2\big(g(B)-1+(g(C_1)-1)\deg(\pr_{1|})\big)=0.
         \end{displaymath}
         \\
         The remaining case is treated analogously.
       \item W.~l.~o.~g.~$i=2$. 
         For $p=(p_1,p_2)\in\Sigma$ we have
         $\tau_p(C_2)=\pr_1^{-1}(p_1)$ is a fibre of $\pr_1$, and
         since $\pr_{1|}$ is unramified,
         $\#\big(B\cap\tau_p(C_2)\big)=\deg(\pr_{1|})=B.C_2$.  
         \\
         Suppose $q,q'\in B$ with
         $\tau_{q-q'}(B)\not=B$. Then $q=\tau_{q-q'}(q')\in
         B\cap\tau_{q-q'}(B)$, and hence $B^2=B.\tau_{q-q'}(B)>0$,
         which contradicts the assumption $B^2=0$. 
       \item Since $\chi(\Sigma)=0$, \cite{FM94} Lemma I.3.18 and
         Proposition I.3.22 imply that
         $g(H)=p_g(\Sigma)=h^0(\Sigma,K_\Sigma)=1$. 
       \item W.~l.~o.~g.~$B.C_2=1$. Let $0\not=p\in C_2$. We claim
         that $B\cap\tau_p(B)=\emptyset$, and hence
         $B^2=B.\tau_p(B)=0$. 
         \\
         Suppose $(a,b)\in B\cap\tau_p(B)$, then there is an
         $(a',b')\in B$ such that
         $(a,b)=\tau_p(a',b')=(a',b'+p)$, i.~e.~$a=a'$ and
         $b=b'+p$. Hence,  $(0,b),(0,b')\in 
         \tau_{-a}(B)\cap C_2$. But, $\tau_{-a}(B).C_2=B.C_2=1$, and
         thus $b'=b=b'+p$ in contradiction to the choice of $p$.
         \\
         $C_1\cong B$ via $\pr_{1|}$ follows from (i).
       \item By (iv) we have $C_1\cong B\cong C_2$.
       \item $\pr_{i|}:B\rightarrow C_i$ is an isomorphism\tom{ (cf.~\cite{EGA} I.5.1.4)}, 
         and has thus degree one. But $\deg(\pr_{i|})=B.C_j$. Thus we
         are done with (iv).
         \oldversion{
           Let $0\not=p\in C_2$ and suppose $(a,b)\in
           B\cap\tau_p(B)$. Then, as above there is a $b'\in C_2$ such
           that $(a,b')\in B$ and $b=b'+p$. Being the graph of a
           morphism, $(a,b),(a,b')\in B$ implies $b'=b=b'+p$ in
           contradiction to the choice of $p$. Hence,
           $B^2=B.\tau_p(B)=0$. 
           }
     \end{enumerate}  
   \end{proof}

   \begin{eexample}\label{ex:product}\leererpunkt
     \begin{enumerate}
       \item Let $C_1=C_2=C=\C/\Lambda$ with $\Lambda=\Z\oplus\tau\Z$, 
         and $\Sigma=C_1\times C_2=C\times C$. The Picard number $\rho(\Sigma)$
         is then either three or four, depending on whether the
         group $\End_0(C)$ of endomorphisms of $C$ fixing $0$ is just
         $\Z$ or larger. Using 
         \cite{Har77} Theorem IV.4.19 and Exercise IV.4.11 we find the 
         following classification.
         \begin{varthm-roman}[Case 1]
           $\exists\; d\in\N$ such that
             $\tau\in\Q[\sqrt{-d}]$, i.~e.~$\Z\subsetneqq\End_0(C)$.
         \end{varthm-roman}
         \noindent
         Then $\rho(\Sigma)=4$ and $\NS(\Sigma)=C_1\Z\oplus
         C_2\Z\oplus C_3\Z\oplus C_4\Z$ where $C_3$ is the
         diagonal in $\Sigma$\tom{, i.~e.~the graph of the
           identity map from $C$ to $C$,} and $C_4$ is the graph of 
         the morphism $\alpha:C\rightarrow C:p\mapsto (b\tau)\cdot
         p$ of degree $|b\tau|^2$, where $0\not=b\in\N$ minimal
         with $b(\tau+\overline{\tau})\in\Z$ and 
         $b\tau\overline{\tau}\in\Z$.  
         Setting $a:=C_3.C_4\geq 1$, the intersection matrix is  
         \begin{displaymath}
           \left(C_j.C_k\right)_{j,k=1,\ldots,4}=
           \left(\begin{array}{cccc}0&1&1&|b\tau|^2\\1&0&1&1\\1&1&0&a\\|b\tau|^2&1&a&0\end{array}\right).
         \end{displaymath}
         \\
         If e.~g.~$\tau=i$, then $C_4=\big\{(c,ic)\;|\;c\in C\big\}$ and 
         \begin{displaymath}
           \left(C_j.C_k\right)_{j,k=1,\ldots,4}=
           \left(\begin{array}{cccc}0&1&1&1\\1&0&1&1\\1&1&0&1\\1&1&1&0\end{array}\right).               
         \end{displaymath}
         \begin{varthm-roman}[Case 2]
           $\nexists
           d\in\N$ such that $\tau\in\Q[\sqrt{-d}]$,  i.~e.~$\Z=\End_0(C)$.
         \end{varthm-roman}
         \noindent
         Then $\rho(\Sigma)=3$ and $\NS(\Sigma)=C_1\Z\oplus
         C_2\Z\oplus C_3\Z$ where again $C_3$ is the diagonal in
         $\Sigma$. The intersection matrix in this case is
         \begin{displaymath}
           \left(C_j.C_k\right)_{j,k=1,2,3}=
           \left(\begin{array}{ccc}0&1&1\\1&0&1\\1&1&0\end{array}\right).                              
         \end{displaymath}
       \item Let $C_1=\C/\Lambda_1$ and $C_2=\C/\Lambda_2$ with
         $\Lambda_1=\Z\oplus\tau_1\Z$, $\tau_1=i$, and
         $\Lambda_2=\Z\oplus\tau_2\Z$, $\tau_2=\frac{1}{2}i$. Then
         $C_1\not\cong C_2$.\tom{\footnote{Suppose $C_1\cong C_2$,
             then there are integers $a,b,c,d\in\Z$ with $ad-bc=\pm 
             1$ such that $\frac{1}{2}i=\frac{ai+b}{ci+d}$
             (cf.~\cite{Har77} Theorem IV.15B). But this leads to
             $ai+b=\frac{1}{2}di-\frac{1}{2}c$, i.~e.~$2a=d$ and
             $-2b=c$. Inserting these relations in the determinant
             equation we get $\pm 1=2a^2+2b^2=2(a^2+b^2)$  
             which would say that two is a divisor of one.}} 
         \\
         We consider the surjective morphisms $\alpha_j:C_1\rightarrow
         C_2$, $j=3,4$, induced by 
         multiplication with the complex numbers
         $\alpha_3=1$ and $\alpha_4=i$ respectively. Denoting by $C_j$
         the graph of $\alpha_j$,  
         we claim, $C_1.C_3=\deg(\alpha_3)=2$ and
         $C_1.C_4=\deg(\alpha_4)=2$. $\alpha_j$ being an unramified 
         covering, we can calculate its degree by counting the
         preimages of $0$. If $p=[a+ib]\in\C/\Lambda_1=C_1$ with
         $0\leq a,b<1$, then
         \begin{displaymath}
           \begin{array}{r@{\;\Leftrightarrow\;}l}
             \alpha_3(p)=0 &
             a+ib=\alpha_3\cdot(a+ib)\in\Lambda_2\\
             & \exists\;r,s\in\Z\;:\;a=r \mbox{ and }
             b=\frac{1}{2}s\\
             & a=0 \mbox{ and }b\in\big\{0,\frac{1}{2}\big\}.
           \end{array}
         \end{displaymath}
         and
         \begin{displaymath}
           \begin{array}{r@{\;\Leftrightarrow\;}l}
             \alpha_4(p)=0 &
             ia-b=\alpha_4\cdot(a+ib)\in\Lambda_2\\
             & \exists\;r,s\in\Z\;:\;-b=r \mbox{ and }
             a=\frac{1}{2}s\\
             & b=0 \mbox{ and }a\in\big\{0,\frac{1}{2}\big\}.
           \end{array}
         \end{displaymath}         
         Moreover, the graphs $C_3$ and $C_4$ intersect only in the
         point $(0,0)$ and the intersection is obviously transversal,
         so $C_3.C_4=1$.
         \\
         Thus $\Sigma=C_1\times C_2$ is an example for a product of
         non-isomorphic elliptic curves with $\rho(\Sigma)=4$,
         $\NS(\Sigma)=C_1\Z\oplus C_2\Z\oplus C_3\Z\oplus C_4\Z$, and
         intersection matrix 
         \begin{displaymath}
           \left(C_j.C_k\right)_{j,k=1,\ldots,4}=
           \left(\begin{array}{cccc}0&1&2&2\\1&0&1&1\\2&1&0&1\\2&1&1&0\end{array}\right).           
         \end{displaymath}
       \item See \cite{HR98} p.~4 for examples $\Sigma=C_1\times C_2$
         with $\rho(\Sigma)=3$ and intersection matrix 
         $$\left(\begin{array}{cccc}0&1&a\\1&0&1\\a&1&0\end{array}\right),\;\;a\not=1.$$           
     \end{enumerate}
   \end{eexample}

   \begin{eremark}\label{rem:isogenous}\leererpunkt
     \begin{enumerate}
       \item If $C_1$ and $C_2$ are isogenous, then there are irreducible
         curves $B\subset\Sigma$ with $B.C_i$ arbitrarily large. 
         \\
         For this just note, that we have a curve $\Gamma\subset\Sigma$
         which is the graph of an isogeny $\alpha:C_1\rightarrow
         C_2$. Denoting by $n_{C_2}:C_2\rightarrow C_2$ the morphism induced
         by the multiplication with $n\in\N$, we have a morphism
         $n_{C_2}\circ \alpha$ whose degree is just
         $n^2\deg(\alpha)$. But the degree is the intersection number of
         the graph with $C_1$. The dual morphism of $n_{C_2}\circ \alpha$ 
         has the the same degree, which then is the intersection
         multiplicity of its graph with $C_2$. (cf.~\cite{Har77}
         Ex.~IV.4.7)
       \item If $C_1$ and $C_2$ are isogenous, then $\Sigma$ might
         very well contain smooth irreducible elliptic curves $B$ which are neither
         isomorphic to $C_1$ nor to $C_2$, and hence cannot be the
         graph of an isogeny between $C_1$ and $C_2$. But being an elliptic curve 
         we have $B^2=0$ by the adjunction formula. If now
         $\NS(\Sigma)=\bigoplus_{i=1}^{\rho(\Sigma)}C_i\Z$, where the
         additional generators are graphs, then
         $B\sim_a\sum_{i=1}^{\rho(\Sigma)} n_iC_i$ with some $n_i<0$. 
         (cf.~\cite{LB92} Ex.~10.6)
     \end{enumerate}
   \end{eremark}

   Throughout the remaining part of the subsection we will restrict
   our attention to the general case, that is that $C_1$ and
   $C_2$ are \emph{not} isogenous. This makes the formulae look much
   nicer, since then $\NS(\Sigma)=C_1\Z\oplus C_2\Z$.

   \begin{lemma}\label{lem:non-isogenous}
     Let $C_1$ and $C_2$ be non-isogenous elliptic curves,
     $D\in\Div(\Sigma)$ with
     $D\sim_aaC_1+bC_2$. 
     \begin{enumerate}
       \item $D^2=0$ if and only if $a=0$ or $b=0$.
       \item If $D$ is an irreducible curve, then we are in one of the
         following cases:
         \begin{enumerate}
           \item $a=0$ and $b=1$,
           \item $a=1$ and $b=0$,
           \item $a,b>0$,
         \end{enumerate}
         and if we are in one of these cases, then there is an
         irreducible curve algebraically equivalent to $D$.
       \item If $D$ is an irreducible curve and $D^2=0$, then either
         $D\sim_aC_1$ or $D\sim_aC_2$.
       \item $D$ is nef if and only if $a,b\geq 0$.
       \item $D$ is ample if and only if $a,b>0$.
       \item $D$ is very ample if and only if $a,b\geq 3$.
     \end{enumerate}
   \end{lemma}
   \begin{proof}\leererpunkt
     \begin{enumerate}
       \item $0=D^2=2ab$ if and only if $a=0$ or $b=0$.
       \item Let us first consider the case that $D$ is irreducible.
         \\
         If $a=0$ or $b=0$, then $D$ is algebraically equivalent
         to a multiple of a fibre of one of the projections $\pr_i$,
         $i=1,2$. In this situation $D^2=0$ and thus the irreducible
         curve $D$ does not intersect any of the fibres
         properly. Hence it must be a union of several fibres, and
         being irreducible it must be a fibre. That is we are in one
         of the first two cases.
         \\
         Suppose now that $a,b\not=0$. Thus $D$ intersects $C_i$
         properly, and $0<D.C_1=b$ and $0<D.C_2=a$.
         \\
         It now remains to show that the mentioned algebraic systems
         contain irreducible curves, which is clear for the first
         two of them. Let therefore $a$ and $b$ be positive. Then
         obviously the linear system $|aC_1+bC_2|_l$ contains no fixed
         component, and being ample by (v) its general element is
         irreducible according to \cite{LB92} Theorem 4.3.5.
       \item Follows from (i) and (ii).
       \item By definition $D$ is nef if and only if $D.D'\geq0$ for every
         irreducible curve $D'\subset\Sigma$. Thus the claim is an
         immediate consequence of (ii).
       \item Since by the Nakai-Moishezon-Criterion ampleness depends only on the
         numerical class of a divisor, we may assume that $D=aC_1+bC_2$.
         Moreover, by \cite{LB92} Proposition 4.5.2 $D$ is ample if and only if $D^2>0$
         and $|D|_l\not=\emptyset$.  
         \\
         If $a,b>0$, then $D^2=2ab>0$ and
         the effective divisor $aC_1+bC_2\in|D|_l$, thus $D$ is ample.
         Conversely, if $D$ is ample, then $0<D^2=2ab$ and
         $0<D.C_1=b$, thus $a,b>0$.
         \oldversion{
           Due to the Nakai-Moishezon-Criterion, $D$ is ample if and
           only if $D^2>0$ and $D.D'>0$ for any irreducible curve. Thus
           the result follows from (ii).
           }
       \item By \cite{LB92} Corollary 4.5.3 and (v) $L=3C_1+3C_2$ is
         very ample. If $a,b\geq3$, then the system
         $|(a-3)C_1+(b-3)C_2|_l$ is basepoint free, which is an
         immediate consequence of the existence of the translation
         morphisms $\tau_p$, $p\in\Sigma$. But then
         $L'=(a-3)C_1+(b-3)C_2$ is globally generated and $D=L+L'$ is
         very ample\tom{ (cf.~\cite{Har77} Lemma II.7.9 and
           Ex.~II.7.5)}.
         \\
         Conversely, if $a<3$, then $D\cap C_2$ is a divisor of degree 
         $D.C_2=a<3$ on the elliptic curve $C_2$ and hence not very
         ample (cf.~\cite{Har77} Example IV.3.3.3). But then $D$ is
         not very ample. Analogously if $b<3$.
     \end{enumerate}
   \end{proof}

   In view of
   (\ref{eq:vanishing:2}\alph{subsection}) and
   Lemma \ref{lem:non-isogenous} (iv) the Condition 
   (\ref{eq:vanishing:3}) becomes obsolete, and Corollary
   \ref{cor:vanishing} has the following form, taking Lemma
   \ref{lem:non-isogenous} (iii) and $K_\Sigma=0$ into account.

   \renewcommand{\thesatz}{\ref{cor:vanishing}\alph{subsection}}
   \begin{corollary}
     Let $C_1$ and $C_2$ be two non-isogenous elliptic curves, $a,b\in\Z$ be
     integers satisfying 
     \begin{equationlist}
       \theequationchange{\ref{eq:vanishing:1}\alph{subsection}}
       \item $ab\geq \sum_{i=1}^r(m_i+1)^2$, and
       \theequationchange{\ref{eq:vanishing:2}\alph{subsection}}
       \item $a,b>\max\{m_i\;|\;i=1,\ldots,r\}$,
       \theequationchangeback
     \end{equationlist}
     then for $z_1,\ldots,z_r\in\Sigma= C_1\times C_2$ in very general
     position and $\nu>0$
     \begin{displaymath}
       H^\nu\left(\Bl_{\underline{z}}(\Sigma),\;
         a\pi^*C_1+b\pi^*C_2-\sum_{i=1}^r m_iE_i\right)=0. 
     \end{displaymath}
   \end{corollary}
   \renewcommand{\thesatz}{\remembersatz}

   As for the existence theorem Corollary \ref{cor:existence-II} we
   work with the very ample divisor class $L=3C_1+3C_2$, and we claim
   that the Conditions (\ref{eq:existence-II:3}), 
   (\ref{eq:existence-II:4}) and (\ref{eq:existence-II:5}) become
   obsolete, while, in view of Lemma \ref{lem:non-isogenous} (iii),
   (\ref{eq:existence-II:1}) and 
   (\ref{eq:existence-II:2}) take the form
   \begin{equationlist}
     \theequationchange{\ref{eq:existence-II:1}\alph{subsection}}
     \item
       $(a-3)(b-3) \geq \sum_{i=1}^r (m_i+1)^2$, and
     \theequationchange{\ref{eq:existence-II:2}\alph{subsection}}
     \item
       $(a-3),(b-3)>\max\{m_i\;|\;i=1,\ldots,r\}$.
   \end{equationlist}
   That is, under these hypotheses there is an irreducible curve
   in $|D|_l$, for any $D\sim_a aC_1+bC_2$, with precisely $r$
   ordinary singular  points of multiplicities $m_1,\ldots,m_r$.
   
   (\ref{eq:existence-II:3}) becomes redundant in view of
   (\ref{eq:existence-II:2}\alph{subsection}) and Lemma
   \ref{lem:non-isogenous} (iv), while (\ref{eq:existence-II:5})
   is fulfilled in view of (\ref{eq:existence-II:6}) and
   $K_\Sigma=0$. It remains to show that $D.L-g(L)\geq m_i+m_j$ for
   all $i,j$. However, by the adjunction formula
   $g(L)=1+\frac{1}{2}L^2=10$, and by
   (\ref{eq:existence-II:2}\alph{subsection}) $D.L-g(L)>
   3(a-3+b-3)>3(m_i+m_j)\geq m_i+m_j$. Thus the claim is proved.

   From these considerations we at once deduce the conditions for the
   existence of an irreducible curve in $|D|_l$, 
   $D\sim_aaC_1+bC_2$, with prescribed
   singularities of arbitrary type, i.~e.~the conditions in Corollary
   \ref{cor:existence-IV}. They come down to 
   \begin{equationlist}
     \theequationchange{\ref{eq:existence-IV:1}\alph{subsection}}
     \item $(a-3)(b-3)
         \geq \frac{207}{5}\sum\limits_{\mu(\ks_i)\leq 38} \mu(\ks_i) +
         29\sum\limits_{\mu(\ks_i)\geq 39}
         \Big(\sqrt{\mu(\ks_i)}+\frac{13}{2\sqrt{29}}\Big)^2$, and
     \theequationchange{\ref{eq:existence-IV:2}\alph{subsection}}
     \item $(a-3),(b-3) >
         \left\{
           \begin{array}{ll}
             \sqrt{\frac{207}{5}}\sqrt{\mu(\ks_1)}-1, & \mbox{ if 
               }\mu(\ks_1)\leq 38,\\[0.3cm]
             \sqrt{29} \sqrt{\mu(\ks_1)}+\frac{11}{2}, & \mbox{
               if }\mu(\ks_1)\geq 39.
           \end{array}
         \right.$
     \theequationchangeback
   \end{equationlist}


   \subsection{Surfaces in $\PC^3$} 

   A smooth projective surface $\Sigma$ in $\PC^3$ is given by a single
   equation $f=0$ with $f\in\C[w,x,y,z]$ homogeneous, and by definition
   the degree of $\Sigma$, 
   say $n$, is just the degree of $f$. For $n=1$,
   $\Sigma\cong\PC^2$, for $n=2$, $\Sigma\cong\PC^1\times\PC^1$, and
   for $n=3$, $\Sigma$ is isomorphic to $\PC^2$ blown up in six points
   in general position. Thus the Picard number $\rho(\Sigma)$,
   i.~e.~the rank of the N\'eron-Severi group, in these cases is $1$, 
   $2$, and $7$ respectively.\tom{\footnote{See e.~g.~\cite{Har77}
       Example II.8.20.3 and Remark V.4.7.1.}} Note that these are
   also precisely the cases where $\Sigma$ is rational.

   In general the Picard number $\rho(\Sigma)$ of a 
   surface in $\PC^3$ may be arbitrarily large,\footnote{E.~g.~the $n$-th Fermat
     surface, given by $w^n+x^n+y^n+z^n=0$ has Picard number $\rho\geq
     3(n-1)(n-2)+1$, with equality if
     $\gcd(n,6)=1$. (cf.~\cite{Shi82} Theorem 7, see also \cite{AS83}
       pp.~1f.~and \cite{IS96} p.~146)}
   but the N\'eron-Severi group always contains a very special member,
   namely the class $H\in\NS(\Sigma)$ of a hyperplane
   section with $H^2=n$. And the class of the canonical divisor is then just
   $(n-4)H$. Moreover, if the degree of $\Sigma$ is at
   least four, that  
   is, if $\Sigma$ is not rational, then it is
   likely that  $\NS(\Sigma)=H\Z$. More precisely, if $n\geq 4$,
   Noether's Theorem says that 
   $\{\Sigma\;|\;\rho(\Sigma)=1, \deg(\Sigma)=n\}$ is a very general
   subset of the projective space of projective surfaces in $\PC^3$ of fixed
   degree $n$, i.~e.~it's complement is an at most countable union of
   lower dimensional subvarieties. (cf.~\cite{Har75} Corollary 3.5
   or \cite{IS96} p.~146)

   Since we consider the case of rational surfaces 
   separately the following considerations thus give a full answer for
   the ``general case'' of a surface in $\PC^3$.

   \renewcommand{\thesatz}{\ref{cor:vanishing}\alph{subsection}}
   \begin{corollary}\label{thm:vanishing-P3:1}
     Let $\Sigma\subset\PC^3$ be a surface in $\PC^3$ of degree $n$,
     $H\in\NS(\Sigma)$ be the algebraic class of a
     hyperplane section, and $d$ an integer satisfying
     \begin{equationlist}
       \theequationchange{\ref{eq:vanishing:1}\alph{subsection}}
       \item $n(d-n+4)^2 \geq 2 \sum_{i=1}^r(m_i+1)^2$, and
       \theequationchange{\ref{eq:vanishing:2}\alph{subsection}}
       \item $(d-n+4)\cdot H.B>\max\{m_i\;|\;i=1,\ldots,r\}$ for any irreducible curve
         $B$ with $B^2=0$ and $\dim|B|_a\geq 1$, and
       \theequationchange{\ref{eq:vanishing:3}\alph{subsection}}
       \item $d\geq n-4$,
       \theequationchangeback
     \end{equationlist}
     then for $z_1,\ldots,z_r\in\Sigma$ in very general position and $\nu>0$
     \begin{displaymath}
       H^\nu\left(\Bl_{\underline{z}}(\Sigma),d\pi^*H-\sum_{i=1}^r m_iE_i\right)=0.
     \end{displaymath}
   \end{corollary}
   \renewcommand{\thesatz}{\remembersatz}

   \begin{eremark}\leererpunkt
     \begin{enumerate}
       \item If $\NS(\Sigma)=H\Z$, then
         (\ref{eq:vanishing:2}\alph{subsection}) is redundant, since there 
         are no irreducible curves $B$ with $B^2=0$.
         Otherwise we would have $B\sim_a kH$ for some $k\in\Z$ and
         $k^2n=B^2=0$ would imply $k=0$, but then $H.B=0$ in contradiction
         to $H$ being ample.
       \item However, a quadric in $\PC^3$ or the K3-surface given by
         $w^4+x^4+y^4+z^4=0$ contain irreducible curves of
         self-intersection zero. 
       \item If $\sum_{i=1}^r(m_i+1)^2>\frac{n}{2}m_i^2$ for all
         $i=1,\ldots,r$ then again
         (\ref{eq:vanishing:2}\alph{subsection}) becomes obsolete in
         view of (\ref{eq:vanishing:1}\alph{subsection}), since
         $H.B>0$ anyway. 
         The above inequality is, for instance, fulfilled if the
         highest multiplicity occurs at least 
         $\frac{n}{2}$ times.
       \item In the existence theorems the condition
         depending on curves of self-intersection will vanish
         in any case. 
     \end{enumerate}
   \end{eremark}

   As for Corollary \ref{cor:existence-II} we claim that if
   $\NS(\Sigma)=H\Z$, then
   \begin{equationlist}
     \theequationchange{\ref{eq:existence-II:1}\alph{subsection}}
     \item
       $n(d-n+3)^2 \geq 2 \sum_{i=1}^r (m_i+1)^2$,
   \end{equationlist}
   ensures the existence of an irreducible curve $C\sim_a dH$
   with precisely $r$ ordinary singular points of multiplicities
   $m_1,\ldots, m_r$ and
   $h^1\big(\Sigma,\kj_{X(\underline{m};\underline{z})/\Sigma}(dH)\big)=0$. 

     The role of the very ample divisor $L$ is filled by a hyperplane
     section, and thus
     $g(L)=1+\frac{L^2+L.K_\Sigma}{2}=\binom{n-1}{2}$. Therefore,
     (\ref{eq:existence-II:1}\alph{subsection}) 
     obviously implies
     \eqref{eq:existence-II:1}, and \eqref{eq:existence-II:4} takes
     the form
     \begin{equation}\label{eq:ex}
       n\cdot(d-n+3)> m_i+2 \text{ for all } i=1,\ldots,r.
     \end{equation}
     However, from (\ref{eq:existence-II:1}\alph{subsection}) we
     deduce for any $i\in\{1,\ldots,r\}$
     \begin{displaymath}
       n\cdot (d-n+3)
       \geq 
       \sqrt{n}\cdot
       \sqrt{2}\cdot (m_i+1)
       \geq m_i+2,
     \end{displaymath}
     unless $n=r=m_1=1$, in which case we are done by the assumption
     $d\geq 3$. Thus \eqref{eq:existence-II:4} is redundant.

     Moreover, there are no curves of self-intersection zero on
     $\Sigma$, and
     it thus remains to verify \eqref{eq:existence-II:3},
     which in this situation takes the form 
     \begin{displaymath}
       d\geq n-3,
     \end{displaymath}
     and follows at once from \eqref{eq:ex}. 

   With the aid of this result the conditions of Corollary
   \ref{cor:existence-IV} for the existence of an irreducible curve
   $C\sim_a dH$ with prescribed singularities $\ks_i$ in this
   situation therefore reduce to
   \begin{equationlist}
     \theequationchange{\ref{eq:existence-IV:1}\alph{subsection}}
     \item $n(d-n+3)^2 
       \geq \frac{414}{5}\sum\limits_{\mu(\ks_i)\leq 38} \mu(\ks_i) +
       58\sum\limits_{\mu(\ks_i)\geq 39}
       \Big(\sqrt{\mu(\ks_i)}+\frac{13}{2\sqrt{29}}\Big)^2$, and
     \theequationchangeback
   \end{equationlist}


   \subsection{K3-Surfaces}

   We note that if $\Sigma$ is a K3-surface then the N\'eron-Severi
   group $\NS(\Sigma)$ and the Picard group $\Pic(\Sigma)$ of $\Sigma$ 
   coincide, i.~e.~$|D|_a=|D|_l$ for every divisor $D$ on
   $\Sigma$. Moreover, an irreducible curve $B$ has self-intersection
   $B^2=0$ if and only if the arithmetical genus of $B$ is one. In
   that case $|B|_l$ is a pencil of elliptic curves without base points
   endowing  $\Sigma$ with the structure of an elliptic fibration over
   $\PC^1$. (cf.~\cite{Mer85} or Proposition \ref{prop:condition}) We,
   therefore, distinguish two cases. 

   \subsubsection{Generic K3-Surfaces}

   Since a generic K3-surface does not possess an elliptic fibration
   the following version of Corollary \ref{cor:vanishing} applies for
   generic K3-surfaces. (cf.~\cite{FM94} I.1.3.7)

   \renewcommand{\thesatz}{\ref{cor:vanishing}\alph{subsection}.\roman{subsubsection}}
   \begin{corollary}\label{thm:vanishing-K3:1}
     Let $\Sigma$ be a K3-surface which is not elliptic, and let $D$ a
     divisor on $\Sigma$ satisfying
     \begin{equationlist}
       \theequationchange{\ref{eq:vanishing:1}\alph{subsection}}
       \item $D^2\geq 2\sum_{i=1}^r(m_i+1)^2$, and
       \theequationchange{\ref{eq:vanishing:3}\alph{subsection}}
       \item $D$ nef,
       \theequationchangeback
     \end{equationlist}
     then for $z_1,\ldots,z_r\in\Sigma$ in very general position and $\nu>0$
     \begin{displaymath}
       H^\nu\left(\Bl_{\underline{z}}(\Sigma),\pi^*D-\sum_{i=1}^r m_iE_i\right)=0.
     \end{displaymath}
   \end{corollary}
   \renewcommand{\thesatz}{\remembersatz}

   In view of equation
   (\ref{eq:existence-II:6})\tom{\footnote{
       $D^2+(2D-L-K_\Sigma)(L+K_\Sigma)+4\sum_{i=1}^rm_i+2r >
       (D-L)^2+2\big(D.L+(D-L).L\big)\geq 0$, by
       (\ref{eq:existence-II:1}\alph{subsection}),
       (\ref{eq:existence-II:3}\alph{subsection}), and since $|D|_l$
       is non-empty even without (\ref{eq:existence-II:5}).}}
   the conditions in
   Corollary \ref{cor:existence-II} reduce to
   \begin{equationlist}
     \theequationchange{\ref{eq:existence-II:1}\alph{subsection}}
     \item $(D-L)^2\geq 2\sum_{i=1}^r 
       (m_i+1)^2$, 
     \theequationchange{\ref{eq:existence-II:3}\alph{subsection}}
     \item $D-L$ nef, and
     \theequationchange{\ref{eq:existence-II:4}\alph{subsection}}
     \item $D.L-2g(L)\geq
       m_i+m_j$ for all $i,j$,
     \theequationchangeback
   \end{equationlist}         
   and, analogously, the conditions in Corollary
   \ref{cor:existence-IV} reduce to (\ref{eq:existence-IV:4}),
   \begin{equationlist}
      \theequationchange{\ref{eq:existence-IV:1}\alph{subsection}}
      \item $(D-L)^2
         \geq \frac{414}{5}\sum\limits_{\mu(\ks_i)\leq 38} \mu(\ks_i) +
         58\sum\limits_{\mu(\ks_i)\geq 39}
         \Big(\sqrt{\mu(\ks_i)}+\frac{13}{2\sqrt{29}}\Big)^2$, and
      \theequationchange{\ref{eq:existence-IV:3}\alph{subsection}}
       \item $D-L$ nef.
      \theequationchangeback
   \end{equationlist} 
   
   \subsubsection{K3-Surfaces with an Elliptic Structure}

   The hypersurface in $\PC^3$ given by the equation
   $x^4+y^4+z^4+u^4=0$ is an example of a K3-surface which is endowed with an
   elliptic fibration. Among the elliptic K3-surfaces the general one
   will possess a unique elliptic fibration while there are examples
   with infinitely many different such fibrations. (cf.~\cite{FM94} I.1.3.7)

   \renewcommand{\thesatz}{\ref{cor:vanishing}\alph{subsection}.\roman{subsubsection}}
   \begin{corollary}\label{thm:vanishing-K3:2}
     Let $\Sigma$ be a K3-surface which possesses an elliptic
     fibration, and let $D$ be a
     divisor on $\Sigma$ satisfying
     \begin{equationlist}
       \theequationchange{\ref{eq:vanishing:1}\alph{subsection}}
       \item $D^2\geq 2\sum_{i=1}^r(m_i+1)^2$, 
       \theequationchange{\ref{eq:vanishing:2}\alph{subsection}}
       \item $D.B>\max\{m_i\;|\;i=1,\ldots,r\}$ for any irreducible curve
         $B$ with $B^2=0$, and
       \theequationchange{\ref{eq:vanishing:3}\alph{subsection}}
       \item $D$ nef,
       \renewcommand{\thesatz}{\remembersatz}
     \end{equationlist}
     then for $z_1,\ldots,z_r\in\Sigma$ in very general position and $\nu>0$
     \begin{displaymath}
       H^\nu\left(\Bl_{\underline{z}}(\Sigma),\pi^*D-\sum_{i=1}^r m_iE_i\right)=0.
     \end{displaymath}
   \end{corollary}
   \renewcommand{\thesatz}{\remembersatz}

   \begin{remark}
     If $\Sigma$ is generic among the elliptic K3-surfaces, i.~e.~admits
     exactly one elliptic fibration, then Condition
     (\ref{eq:vanishing:2}\alph{subsection}) means that a curve in
     $|D|_l$ meets a general fibre in at least
     $k=\max\{m_i\;|\;i=1,\ldots,r\}$ distinct points.
   \end{remark}

   The conditions in Corollary \ref{cor:existence-II} then reduce to
   (\ref{eq:existence-II:1}\alph{subsection}), (\ref{eq:existence-II:3}\alph{subsection}),
   (\ref{eq:existence-II:4}\alph{subsection}), and
   \begin{equationlist}
     \theequationchange{\ref{eq:existence-II:2}\alph{subsection}}
     \item $(D-L).B>\max\{m_i\;|\;i=1,\ldots,r\}$ for any curve $B$ with $B^2=0$.
     \theequationchangeback
   \end{equationlist}         
   Similarly, the conditions in Corollary
   \ref{cor:existence-IV} reduce to
   (\ref{eq:existence-IV:1}\alph{subsection}),
   (\ref{eq:existence-IV:3}\alph{subsection}),
   (\ref{eq:existence-IV:4}), and 
   \begin{equationlist}
      \theequationchange{\ref{eq:existence-IV:2}\alph{subsection}}
      \item $(D-L).B>
         \left\{
           \begin{array}{ll}
             \sqrt{\frac{207}{5}}\sqrt{\mu(\ks_1)}-1, & \mbox{ if 
               }\mu(\ks_1)\leq 38,\\[0.3cm]
             \sqrt{29} \sqrt{\mu(\ks_1)}+\frac{11}{2}, & \mbox{
               if }\mu(\ks_1)\geq 39,
           \end{array}
         \right.$
         \\[0.3cm]
         for any irreducible curve $B$ with $B^2=0$.
      \theequationchangeback 
   \end{equationlist} 


   \newpage

\begin{appendix}
     
   \renewcommand{\thesatz}{\Alph{section}.\arabic{satz}}
   \renewcommand{\theequation}{\Alph{section}.\arabic{equation}}
   \renewcommand{\thesubsection}{\Alph{section}.\alph{subsection}}

   \section{Very General Position}\label{sec:verygeneralposition}
   \setcounter{equation}{0}





   It is our first aim to show that if there is a curve passing
   through points $z_1,\ldots,z_r\in\Sigma$ in very general position
   with multiplicities $n_1,\ldots,n_r$ then it can be equimultiply
   deformed in its algebraic system in a good way - i.~e.~suitable for
   Lemma \ref{lem:deformation}.

   \oldversion{
   \begin{lemma}\label{lem:alg-hilb}
     Let $\Sigma$ be embedded in $\PC^N$, and $h\in\Q[x]$ fixed. We denote by
     $\Hilb_\Sigma^h$ the Hilbert scheme of curves on $\Sigma$ with
     Hilbert polynomial $h$. 

     If $B\in\Hilb_\Sigma^h$,
     then $|B|_a$ is a connected component $\Hilb_\Sigma^h$.
   \end{lemma}
   \begin{proof}
     $\Hilb_\Sigma^h$ is a projective scheme (cf.~\cite{Mum66} Chapter 
     15) together with a universal family of curves. Restricting the
     family to the connected component $H\subseteq\Hilb_\Sigma^h$
     containing $B$ shows that the members of $H$ are indeed
     algebraically equivalent to $B$. Moreover, the universal property 
     of $\Hilb_\Sigma^h$ forces any connected algebraic family of curves
     containing $B$ to be induced by the family over $H$. Thus $|B|_a=H$.
   \end{proof}
   }

   \begin{lemma}\label{lem:closed}
     Let $B\subset\Sigma$ be a curve, and $\underline{n}\in\N_0^r$. Then
     \begin{displaymath}
       V_{B,\underline{n}} = \big\{\underline{z}\in\Sigma^r \;\big|\;
       \exists\: C\in |B|_a\: :\: \mult_{z_i}(C)\geq n_i
       \;\forall i=1,\ldots,r\big\}
     \end{displaymath}
     is a closed subset of $\Sigma^r$.
   \end{lemma}
   \begin{proof}
     \begin{varthm-roman}[Step 1]
       Show first that for $n\in\N_0$
       \begin{displaymath}
          X_{B,n}:= \big\{(C,z)\in H\times\Sigma \;|\; \mult_z(C)\geq n\big\}
       \end{displaymath}
       is a closed subset of $H\times\Sigma$, where $H:=|B|_a$.
     \end{varthm-roman}
     Being the reduction of a connected component of the Hilbert scheme $\Hilb_\Sigma$, 
     $H$ is a projective variety endowed
     with a universal family of curves, giving rise to the following
     diagram of morphisms
     \begin{displaymath}
       \xymatrix@C0.8cm{
         {\kc = \bigcup_{C\in H}\{C\}\times C\;} \ar@{^{(}->}[r]\ar[dr] &
         H\times\Sigma \ar@{->>}[r]_(0.58){\pr_\Sigma}\ar@{->>}[d]^{\pr_H}
         &{\Sigma}\\
         &H,&
         }
     \end{displaymath}
     \tom{\footnote{An open bracket had to be closed ).}}
     where $\kc$ is an effective Cartier divisor on $H\times \Sigma$ 
     with $\kc_{|\{C\}\times\Sigma}=C$.\tom{\footnote{For the
         definition of an algebraic family of curves see \cite{Har77}
         V.~ex.~1.7.}}  

     Let $s\in H^0\big(H\times\Sigma,\ko_{H\times\Sigma}(\kc)\big)$ be 
     a global section defining $\kc$. Then 
     \begin{displaymath}
       X_{B,n} = \big\{\eta=(C,z)\in H\times\Sigma \;\big|\; s_\eta \in
       (\m_{\Sigma,z}^n + \m_{H,C})\cdot\ko_{H\times\Sigma,\eta}\big\}.
     \end{displaymath}\tom{\footnote{$\m_{\Sigma,z}^n =
           \big(\m_{\Sigma,z}^n+\m_{H,C}\big)\cdot\ko_{H\times\Sigma,\eta}
           / \m_{H,C}\cdot\ko_{H\times\Sigma,\eta}$ and
           $C=\kc_{|\{C\}\times\Sigma}$ is locally in $z$ given by the 
           image of $s_\eta$ in
           $\ko_{\Sigma,z}=\ko_{H\times\Sigma,\eta}/\m_{H,C}\cdot\ko_{H\times\Sigma,\eta}$.}}

     We may consider a finite open affine covering of $H\times\Sigma$ of
     the form $\{H_i\times U_j\;|\; i\in I,j\in J\}$, $H_i\subset H$
     and $U_j\subset \Sigma$ open,  such that $\kc$
     is locally on $H_i\times U_j$ given by one polynomial equation,
     say 
     \begin{displaymath}
       s_{i,j}(\underline{a},\underline{b}) =0, \mbox{ for } \underline{a}\in H_i,\:
       \underline{b}\in U_j.
     \end{displaymath} \tom{\footnote{The $\underline{a}$ and
           $\underline{b}$ denote coordinates on the affine ambient
           spaces of $H_i\subseteq\A^{N_i}$ respectively $U_j\subseteq\A^{M_j}$.}}
     It suffices to show that $X_{B,n}\cap (H_i\times U_j)$ is
     closed in $H_i\times U_j$ for all $i,j$. 

     However, for $\eta=(C,z)=(\underline{a},\underline{b})\in H_i\times
     U_j$ we have
     \begin{displaymath}
       s_\eta \in \big(\m_{\Sigma,z}^n + \m_{H,C}\big)\cdot\ko_{H\times\Sigma,\eta}       
     \end{displaymath}
     if and only if
     \begin{displaymath}
       s_{i,j}(\underline{a},\underline{b}) =0
       \mbox{ and}
     \end{displaymath}
     \begin{displaymath}
       \frac{\partial^\alpha s_{i,j}}{\partial
         \underline{b}^\alpha}(\underline{a},\underline{b}) = 0,
       \mbox{ for all } |\alpha|\leq n-1.
     \end{displaymath}
     Thus,
     \begin{displaymath}
       X_{B,n}\cap (H_i\times U_j)
       =\Bigg\{(\underline{a},\underline{b})\in H_i\times U_j \;\Bigg|\;
      s_{i,j}(\underline{a},\underline{b}) = 0 = \frac{\partial^\alpha s_{i,j}}{\partial
         \underline{b}^\alpha}(\underline{a},\underline{b}),\: \forall \:
       |\alpha|\leq n-1\Bigg\}
     \end{displaymath}
     is a closed subvariety of $H_i\times U_j$\tom{, since the $\frac{\partial s_{i,j}}{\partial
         a_\nu}(\underline{a},\underline{b})$ and $\frac{\partial^\alpha s_{i,j}}{\partial
         \underline{b}^\alpha}(\underline{a},\underline{b})$ are
       polynomial expressions in $\underline{a}$ and $\underline{b}$}.

     \begin{varthm-roman}[Step 2]
       $V_{B,\underline{n}}$ is a closed subset of
       $\Sigma^r$.
     \end{varthm-roman}
     By Step 1 for $i=1,\ldots,r$ the set 
     \begin{displaymath}
       X_{B,\underline{n},i} := \big\{(\underline{z},C)\in
       \Sigma^r\times H \;\big|\; \mult_{z_i}(C)\geq n_i\big\} \cong
       \Sigma^{r-1}\times X_{B,n_i}
     \end{displaymath}
     is a closed subset of $\Sigma^r\times H\tom{\cong \Sigma^{r-1}\times
       H\times \Sigma}$.
     Considering now
     \begin{displaymath}
       \xymatrix@C0.8cm{
         X_{B,\underline{n}}:=\bigcap\limits_{i=1}^r X_{B,\underline{n},i}\;
         \ar@{^{(}->}[r]\ar[dr]_\rho &
         {\Sigma^r\times H} \ar@{->>}[d]\\
         &{\Sigma^r},
         }
     \end{displaymath}
     \tom{\footnote{An open bracket had to be closed ).}}
     we find that $V_{B,\underline{n}}=\rho(X_{B,\underline{n}})$, being the image of a closed
     subset under a morphism between projective varieties, is a closed subset of 
     $\Sigma^r$ (cf.~\cite{Har77} Ex.~II.4.4). 
   \end{proof}

   \begin{corollary}\label{cor:verygeneral}
     Then the complement of the set 
     \begin{displaymath}
       V = \bigcup_{B\in\Hilb_\Sigma}\;\bigcup_{\underline{n}\in\N_0^r}
       \big\{V_{B,\underline{n}} \;|\; V_{B,\underline{n}}\not=\Sigma^r\big\}
     \end{displaymath}
     is very general, where $\Hilb_\Sigma$ is the Hilbert scheme of
     curves on $\Sigma$.

     In particular, there is a very general subset $U\subseteq
     \Sigma^r$ such that if for some $\underline{z}\in U$ there is a
     curve $B\subset\Sigma$ with $\mult_{z_i}(B)=n_i$ for
     $i=1,\ldots,r$, then for any $\underline{z}'\in U$ there is 
     a curve $B'\in|B|_a$ with $\mult_{z_i'}(B')\geq n_i$.
   \end{corollary}
   \begin{proof}
     Fixing some embedding $\Sigma\subseteq\PC^n$ and $h\in\Q[x]$,
     $\Hilb_\Sigma^h$ is a projective variety and has thus only
     finitely many connected components. Thus the Hilbert scheme
     $\Hilb_\Sigma$ has only a countable number of connected components,
     and we have only a countable number of
     different $V_{B,\underline{n}}$, where $B$ runs through
     $\Hilb_\Sigma$ and $\underline{n}$ through $\N^r$.
     By Lemma \ref{lem:closed} the sets $V_{B,\underline{n}}$
     are closed, hence their complements
     $\Sigma^r\setminus
     V_{B,\underline{n}}$ are open.
     But then
     \begin{displaymath}
       U=\Sigma^r\setminus V
       =\bigcap_{B\in\Hilb_\Sigma}\;\bigcap_{\underline{n}\in\N_0^r}
       \big\{\Sigma^r \setminus V_{B,\underline{n}}\;|\; V_{B,\underline{n}}\not=\Sigma^r\big\}
     \end{displaymath}
     is an at most countable intersection of open dense subsets of
     $\Sigma^r$, and is hence very general. 
   \end{proof}

   \tom{
     If $\Sigma$ is regular, i.~e.~the irregularity
     $q(\Sigma)=h^1(\Sigma,\ko_\Sigma)=0$, algebraic and linear
     equivalence coincide, and thus $\Hilb_\Sigma^h=|D|_l$ if $D$ is
     any divisor with Hilbert polynomial $h$. This makes the proofs
     given above a bit simpler.

     \begin{lemma}\label{lem:closed-weak}
       Given $\underline{n}\in\N_0^r$ and $D\in\Pic(\Sigma)$, the set 
       \begin{displaymath}
         V_{D,\underline{n}} = \big\{\underline{z}\in\Sigma^r \;|\;
         \exists\: C\in|D|_l\: :\: \mult_{z_i}(C)\geq n_i\big\}
       \end{displaymath}
       is a closed subset of $\Sigma^r$.
     \end{lemma}
     \begin{proof}
       Fix an affine covering $\Sigma=U_1\cup\ldots\cup U_k$ of
       $\Sigma$, and a basis $s_0,\ldots,s_n$ of
       $H^0\big(\Sigma,\ko_\Sigma(D)\big)$. 
       
       It suffices to show that $V_{D,\underline{n}} \cap
       \big(U_{j_1}\times\ldots\times U_{j_r}\big)$ is closed
       in $U_{j_1}\times\ldots\times U_{j_r}$ for 
       all $\underline{j}=(j_1,\ldots,j_r)\in \{1,\ldots,k\}^r$\tom{,
         since those sets form an open covering of $V_{D,\underline{n}}$}.
       
       Consider the set
       \begin{displaymath}
         V_{\alpha,i,\underline{j}} =
         \big\{(\underline{z},\underline{a})\in
         U_{j_1}\times\ldots\times U_{j_r}\times\PC^n \;|\;
         a_0(D^\alpha s_0)(z_i) + \ldots + a_n(D^\alpha s_n)(z_i) =
         0\big\}, 
       \end{displaymath}
       where $\alpha$ is a multi index and $D^\alpha$ denotes the
       corresponding differential operator.\tom{\footnote{The dimension
           of the multi index depends on the embedding of the affine
           chart $U_{j_i}$ into some $\A^d$.}}
       This is a closed subvariety of $U_{j_1}\times\ldots\times
       U_{j_r}\times\PC^n$.\tom{\footnote{On $U_j$ the $s_i$ are given
           as quotients of polynomials where the denominator is not
           vanishing on $U_j$. Thus we may assume w.~l.~o.~g.~that $s_i$ 
           is represented by a polynomial, and hence the above equation
           is polynomial in $z_i$ and even linear in the $a_i$.}}
       Thus, also the set
       \begin{displaymath}
         V_{\underline{j}}=\bigcap_{i=1,\ldots,r}\;\bigcap_{|\alpha|\leq
           n_i-1} V_{\alpha,i,\underline{j}}
       \end{displaymath}
       is a closed subset of $U_{j_1}\times\ldots\times 
       U_{j_r}\times\PC^n$. Considering now the projection to $U_{j_1}\times\ldots\times
       U_{j_r}$,
       \begin{displaymath}
         \xymatrix@C0.6cm{
           V_{\underline{j}}\; \ar@{^{(}->}[r]\ar[dr]_\rho & 
           U_{j_1}\times\ldots\times U_{j_r}\times\PC^n \ar@{->>}[d] \\
           & U_{j_1}\times\ldots\times U_{j_r},
           }
       \end{displaymath}
       \tom{\footnote{An open bracket had to be closed ).}}
       we find that $\im(\rho)=V\cap \big(U_{j_1}\times\ldots\times       
       U_{j_r}\big)$ is closed in $U_{j_1}\times\ldots\times U_{j_r}$
       (cf.~\cite{Har92} Theorem 3.12).
     \end{proof}
     
     \begin{corollary}\label{cor:verygeneral-weak}
       Let $\Sigma$ be regular.
       Then the complement of the set 
       \begin{displaymath}
         V = \bigcup_{D\in\Pic(\Sigma)}\;\bigcup_{\underline{n}\in\N_0^r}
         \big\{V_{D,\underline{n}} \;|\; V_{D,\underline{n}}\not=\Sigma^r\big\}
       \end{displaymath}
       is very general.
       
       In particular, there is a very general subset $U\subseteq
       \Sigma^r$ such that if for some $\underline{z}\in U$ there is a
       curve $B\subset\Sigma$ with $\mult_{z_i}=n_i$ for
       $i=1,\ldots,r$, then for any $\underline{z}'\in U$ there is 
       a curve $B'\in|B|_l$ with $\mult_{z_i'}(B')\geq n_i$.
     \end{corollary}
     \begin{proof}
       By Lemma \ref{lem:closed-weak} the sets $V_{D,\underline{n}}$
       are closed, hence their complement
       $\Sigma^r\setminus
       V_{D,\underline{n}}$ is open.\tom{\footnote{See \cite{Har77}
           II.~ex.~3.18.}}
       Since $\Sigma$ is regular, $\Pic(\Sigma)\cong\NS(\Sigma)$ is a 
       finitely generated abelian group, hence
       countable.\tom{\footnote{See \cite{IS96} Chapters 3.2 and
           3.3. From the exponential sequence it follows that
           $\NS(\Sigma)\cong\Pic(\Sigma)/\Pic^0(\Sigma)$ with
           $\Pic^0(\Sigma)\cong \C^q/\Z^{2q}$, $q=q(\Sigma)$. Thus,
           $\Pic(\Sigma)$ is countable if and only if $\Sigma$ is
           regular.}}
       But then
       \begin{displaymath}
         U=\Sigma^r\setminus V
         =\bigcap_{D\in\Pic(\Sigma)}\;\bigcap_{\underline{n}\in\N_0^r}
         \big\{\Sigma^r \setminus V_{D,\underline{n}}\;|\; V_{D,\underline{n}}
         \not=\Sigma^r\big\}
       \end{displaymath}
       is an at most countable intersection of open dense subsets of
       $\Sigma^r$ and is hence very general.
     \end{proof}
     } 

   In the proof of Theorem \ref{thm:vanishing} we use at some place
   the result of Corollary \ref{cor:countable}. We could instead use
   Corollary \ref{cor:verygeneral}. However, since the results are quite
   nice and simple to prove we just give them.

   \oldversion{
   \begin{lemma}\label{lem:algsyst}
     Let $\Sigma$ be embedded in $\PC^N$, 
     and $h\in\Q[x]$. Then the number of curves $C$ 
     with Hilbert polynomial $h$ and $\dim|C|_a=0$ is finite.  
   \end{lemma}
   \begin{proof}
     If $\Hilb_\Sigma^h$ denotes the Hilbert scheme of curves 
     on $\Sigma$ with fixed Hilbert polynomial $h$, then 
     $\Hilb_\Sigma^h$ is a projective scheme (cf.~\cite{Mum66} Chapter
     15), and hence has only finitely many connected components.   
     If now $C\subseteq\Sigma$ is a curve with 
     Hilbert polynomial $h$, belonging to the connected component 
     $H=|C|_a$ of $\Hilb_\Sigma^h$
     and, if in 
     addition $\dim|C|_a=0$, then $H=\{C\}$, which proves the lemma. 
   \end{proof}
   }

   \begin{ecorollary}\label{cor:countable}\leererpunkt
     \begin{enumerate}
        \item The number of curves $B$ in $\Sigma$ with $\dim|B|_a=0$ is at most
          countable. 
        \item The number of exceptional curves in $\Sigma$ (i.~e.~curves with
          negative self intersection) is at most
          countable.
        \item There is a very general subset $U$ of $\Sigma^r$, $r\geq 
          1$, such
          that for $\underline{z}\in U$ no $z_i$ belongs to a curve $B\subset\Sigma$
          with $\dim|B|_a=0$, in particular to no exceptional curve.
     \end{enumerate}
   \end{ecorollary}
   \begin{proof}\leererpunkt
     \begin{enumerate}
        \item By definition $|B|_a$ is a connected
          component of $\Hilb_\Sigma$, whose number is at most
          countable\tom{ (see proof of Lemma \ref{cor:verygeneral})}.  
          If in addition $\dim|B|_a=0$, then $|B|_a=\{B\}$ which
          proves the claim.
        \item Curves of negative self-intersection are not
          algebraically equivalent to any other curve (cf.~\cite{IS96} 
          p.~153).
        \item Follows from (i).
     \end{enumerate}
   \end{proof}

   \begin{example}[Kodaira]
     Let $z_1,\ldots,z_9\in\PC^2$ be in very general position\footnote{To
       be precise, no three of the nine points should be collinear,
       and after any finite number of quadratic Cremona
       transformations centred at the $z_i$ (respectively the newly
       obtained centres) still no three should be collinear. Thus the
       admissible tuples in $(\PC^2)^9$ form a very general set,
       cf.~\cite{Har77} Ex.~V.4.15.} and let $\Sigma=\Bl_{\underline{z}}\big(\PC^2\big)$ be the blow up 
     of $\PC^2$ in $\underline{z}=(z_1,\ldots,z_9)$. Then $\Sigma$
     contains infinitely many irreducible smooth rational $-1$-curves, 
     i.~e.~exceptional curves of the first kind.
   \end{example}
   \begin{proof}
     It suffices to find an infinite number of irreducible curves $C$ in $\PC^2$
     such that 
     \begin{equation}\label{eq:ex-1}
       d^2-\sum_{i=1}^9 m_i^2 =-1,
     \end{equation}
     and 
     \begin{equation}\label{eq:ex-2}
       p_a(C)-\sum_{i=1}^9 \frac{m_i(m_i-1)}{2} =0,
     \end{equation}
     where $m_i=\mult_{z_i}(C)$ and $d=\deg(C)$, since the expression
     in (\ref{eq:ex-1}) is the self intersection of the strict
     transform $\widetilde{C}=\Bl_{\underline{z}}(C)$ of $C$ and (\ref{eq:ex-2}) gives its
     arithmetical genus.\tom{\footnote{
       \begin{displaymath}
         \begin{array}{lcl}
           p_a\big(\widetilde{C}\big)&=&1+\frac{K_\Sigma.\widetilde{C}+\widetilde{C}^2}{2} 
           =1+\frac{\big(K_{\PC^2}+\sum E_i\big).\big(C-\sum
             m_iE_i\big)+\big(C-\sum m_iE_i\big)^2}{2} \\
           &=& 1+\frac{K_{\PC^2}.C+C^2}{2} - \sum \frac{m_i(m_i-1)}{2}
           = p_a(C) - \sum \frac{m_i(m_i-1)}{2}.
         \end{array}
       \end{displaymath}}}
     In particular $\widetilde{C}$ cannot contain any singularities,
     since they would contribute to the arithmetical genus, and, being
     irreducible anyway, $\widetilde{C}$ is an exceptional curve of
     the first kind.

     We are going to deduce the existence of these curves with the aid 
     of quadratic Cremona transformations. 

     \begin{varthm-roman}[Claim]
       If for some $d>0$ and $ m_1,\ldots,m_9\geq 0$ with $3d-\sum_{i=1}^9 m_i=1$ there
       is an irreducible curve
       $C\in
       \big|\kj_{X(\underline{m};\underline{z})}(d)\big|_l$,
       then $T(C)\in
       \big|\kj_{X(\underline{m'};\underline{z'})}(d+a)\big|_l$ 
       is an irreducible curve, where
       \begin{itemize}
          \item $\{i,j,k\}\subset\{1,\ldots,9\}$ are such that
            $m_i+m_j+m_k<d$, 
          \item $T:\PC^2\dashrightarrow \PC^2$ is the quadratic
            Cremona transformation at $z_i,z_j,z_k$,
          \item $z_\nu'=\left\{\begin{array}{ll}z_\nu,&\mbox{if
                  }\nu\not=i,j,k,\\ T(\overline{z_\lambda
                  z_\mu}),&\mbox{if
                  }\{\nu,\lambda,\mu\}=\{i,j,k\},\end{array}\right.$
          \item $m_\nu'=\left\{\begin{array}{ll}m_\nu,&\mbox{if
                  }\nu\not=i,j,k,\\
                m_\nu+a,&\mbox{else, and}\end{array}\right.$
          \item $a=d-(m_i+m_j+m_k).$
       \end{itemize}
       Note that, $3(d+a)-\sum_{i=1}^9 m_i'=1$, i.~e.~we may iterate
       the process since  
       the hypothesis of the claim will be preserved.
     \end{varthm-roman}
     Since $3d > \sum_{i=1}^9 m_i$, there must be a triple $(i,j,k)$
     such that $d > m_i+m_j+m_k$. 

     Let us now consider the following diagram
     \begin{displaymath}
       \xymatrix@C0.6cm{
         &\Sigma=\Bl_{z_i,z_j,z_k}(\PC^2)=\Bl_{z_i',z_j',z_k'}(\PC^2)\ar[dl]_\pi\ar[dr]^{\pi'} &\\
         \PC^2\ar@{-->}[rr]^T&&\PC^2,
         }
     \end{displaymath}
     and let us denote the exceptional divisors of $\pi$ by $E_i$ and
     those of $\pi'$ by $E_i'$. Moreover, let
     $\widetilde{C}=\Bl_{z_i,z_j,z_k}(C)$ be the strict transform of 
     $C$ under $\pi$ and let
     $\widetilde{T(C)}=\Bl_{z_i',z_j',z_k'}\big(T(C)\big)$ be the strict
     transform of $T(C)$ under $\pi'$. Then of course
     $\widetilde{C}=\widetilde{T(C)}$, and
     $T(C)$, being the projection $\pi'\big(\widetilde{C}\big)$ 
     of the strict transform $\widetilde{C}$ of the
     irreducible curve $C$, is of
     course an irreducible curve. Note that the condition $d >
     m_i+m_j+m_k$ ensures that $\widetilde{C}$ is not one of the
     curves which are contracted. It thus suffices to verify 
     \begin{displaymath}
       \deg\big(T(C)\big) = d+a,
     \end{displaymath}
     and
     \begin{displaymath}
       m_i'=\mult_{z_i'}\big(T(C)\big)=
       \left\{\begin{array}{ll}
                m_\nu,&\mbox{if }\nu\not=i,j,k,\\
                m_\nu+a,&\mbox{else.}
              \end{array}\right.
     \end{displaymath}
     Since outside the lines $\overline{z_i z_j}$, $\overline{z_i
       z_k}$, and $\overline{z_k z_j}$ the transformation $T$ is an
     isomorphism and since by hypothesis none of the remaining $z_\nu$ 
     belongs to one of these lines we clearly have $m_\nu'=m_\nu$ for
     $\nu\not=i,j,k$.  Moreover, we have
     \begin{displaymath}
       \begin{array}{lcl}
         m_i'& = &\widetilde{T(C)}.E_i' =
       \widetilde{C}.\Bl_{z_i,z_j,z_k}(\overline{z_j z_k})\\
       &=& (\pi^*C-\sum m_lE_l).(\pi^*\overline{z_j z_k}-E_j-E_k)\\
       &=& C.\overline{z_j z_k} - m_j - m_k = d-m_j-m_k=m_i+a.
     \end{array}
     \end{displaymath}
     Analogously for $m_j'$ and $m_k'$.

     Finally we find
     \begin{displaymath}
       \begin{array}{lcl}
         \deg\big(T(C)\big) &=& T(C).\overline{z_i'
           z_j'}={\pi'}^*T(C).{\pi'}^*\overline{z_i' z_j'}\\
         &=&\big(\widetilde{T(C)}+\sum
         m_\nu'E_\nu'\big).\big(E_k+E_i'+E_j'\big)\\
         &=&\widetilde{C}.E_k+\sum m_\nu' E_\nu'.E_k \\
         &=& m_k+m_i'+m_j' =d+a.
       \end{array}
     \end{displaymath}
     This proves the claim.

     Let us now show by induction that for any $d>0$ there is an irreducible curve
     $C$ of degree 
     $d'\geq d$ satisfying (\ref{eq:ex-1}) and (\ref{eq:ex-2}). 
     For $d=1$ the line $C=\overline{z_1 z_2}$ through $z_1$ and $z_2$ 
     gives the induction start. Given some suitable curve of degree
     $d'\geq d$ the above claim then ensures that
     through points in very general position there is an irreducible curve of
     higher degree satisfying (\ref{eq:ex-1}) and (\ref{eq:ex-2}),
     since $a=d-(m_1+m_2+m_3)>0$. Thus the induction step is done. 
   \end{proof}

   The example shows that a smooth projective surface $\Sigma$ may
   indeed carry an infinite number of exceptional curves - even of the first
   kind. According to Nagata (\cite{Nag60} Theorem 4a, 
   p.~283) the example is due to Kodaira.  For further references
   on the example see 
   \cite{Har77} Ex.~V.4.15, \cite{BS95} Example 4.2.7, or \cite{Fra41}. \cite{IS96}
   p.~198 Example 3 shows that also $\PC^2$ blown up in the nine
   intersection points of two plane cubics carries infinitely 
   many exceptional curves of the first kind.


   \section{Condition (\ref{eq:vanishing:2})}\label{sec:condition}
   \setcounter{equation}{0}

   \begin{proposition}\label{prop:condition}
     Suppose that $B\subset\Sigma$ is an irreducible curve with
     $B^2=0$ and $\dim|B|_a\geq 1$, then
     \begin{equationlist}
        \item[eq:condition:1] $|B|_a$ is an irreducible projective curve, and 
        \item[eq:condition:2] there is a fibration
          $f:\Sigma\rightarrow H$ whose fibres are just the
          elements of $|B|_a$, and $H$
          is the normalisation of $|B|_a$.
     \end{equationlist}
   \end{proposition}

   We are proving the proposition in several steps. 

   \begin{proposition}\label{prop:isomorphism}
     Let $f:Y'\rightarrow Y$ be a finite flat morphism of noetherian schemes with
     $Y$ irreducible such that for some point $y_0\in Y$ the fibre
     $Y'_{y_0}=f^{-1}(y_0)=Y'\times_Y\Spec\big(k(y_0)\big)$ is a 
     single reduced point.\tom{\footnote{The assumption ``reduced'' is 
         necessary, since finite flat morphisms may very well have
         only non-reduced fibres. Consider the ring homomorphism
         $\varphi:A=\C[x]\rightarrow\C[x,y]/(y^2)=B : x\mapsto x$,
         making $B\cong A\oplus y\cdot A$ into a free and hence flat
         $A$-module with only non-reduced fibres (= double-points).}}

     Then the structure map $f^{\#}:\ko_Y\longrightarrow f_*\ko_{Y'}$
     is an isomorphism, and hence so is $f$. 
   \end{proposition}
   \begin{proof}
     Since there is at least one connected reduced fibre $Y'_{y_0}$,  
     semicontinuity of flat, proper\tom{\footnote{By \cite{Har77}
         Ex.~II.4.1 $f$ is proper since $f$ is finite.}} morphisms in the version
     \cite{EGA} IV.12.2.4 (vi) implies that there is an
     open dense subset $U\subseteq Y$ such that $Y'_y$ 
     is connected and reduced, hence a single reduced point,
     $\forall\;y\in U$. ($U$ dense in $Y$ is 
     due to the fact that $Y$ is irreducible.) 

     Thus the assumptions
     are stable under restriction to open subschemes of $Y$, and since
     the claim that we have to show is local on $Y$, we may assume that 
     $Y=\Spec(A)$ is affine. Moreover, $f$ being finite, thus affine,
     we have $Y'=\Spec(B)$ is also affine. 

     Since $f$ is flat it is open and
     hence dominates the irreducible affine variety $Y$ and, therefore, induces an
     inclusion of rings $A\hookrightarrow 
     B$. It now suffices to show:
     \begin{varthm-roman}[Claim]
       $A\hookrightarrow B$ is an isomorphism.
     \end{varthm-roman}
     By assumption there exists a $y=P\in\Spec(A)=Y$ such that
     $Y'_y=f^{-1}(y)=\Spec(B_P/PB_P)$ 
     is a single point with reduced structure. In particular we have
     for the multiplicity of $Y'_y=\Spec(B_P/PB_P)$ over
     $\{y\}=\Spec(A_P/PA_P)$ 
     \begin{displaymath}
       1=\mu(Y'_y)=\length_{A_P/PA_P}(B_P/PB_P),
     \end{displaymath}
     \tom{(cf.~\cite{Har77} p.~51 for the definition of the
       multiplicity)} which implies that  
     \begin{displaymath}
       A_P/PA_P\hookrightarrow B_P/PB_P
     \end{displaymath}
     is an isomorphism. Hence by Nakayama's
     Lemma\tom{ (cf.~Algebra II, SS 1996, Blatt 3, Aufgabe 3)} also 
     \begin{displaymath}
       A_P\hookrightarrow B_P
     \end{displaymath}
     is an isomorphism, that is, $B_P$ is free of rank 1 over
     $A_P$. $B$ being locally free\tom{\footnote{Since
         $A\hookrightarrow B$ is flat!}} over $A$, with $A/\sqrt{0}$ an
     integral domain, thus fulfils 
     \begin{displaymath}
       A_Q\hookrightarrow B_Q
     \end{displaymath}
     is an isomorphism for all $Q\in\Spec(A)$, and hence the claim
     follows.\tom{ (cf.~\cite{AM69} Proposition 3.9)}
   \end{proof}

   \tom{
     \begin{remark}
       Some comments from Bernd Kreu{\ss}ler.
       \begin{enumerate}
          \item Let $Y'$ be the standard example of a non-separated
            variety, the affine line over $k$ with the origin doubled, and let
            $f:Y'\rightarrow Y=\A^1_k$ be the projection to $Y$. Then
            $f$ is quasi-finite, but NOT finite, since the preimage of 
            the affine $Y$ is not affine. However, $f_*\ko_{Y'}=\ko_Y$ 
            (see \cite{Har77} II.2.3.5/6), but $f$ is certainly not an 
            isomorphism.
          \item By \cite{Har77} III.11.3 for $f:Y'\rightarrow Y$
            projective between noetherian schemes the condition
            $f_*\ko_{Y'}=\ko_Y$ implies that the fibres of $f$ are
            connected. The converse is also true, if the fibres are
            connected, the $f_*\ko_{Y'}=\ko_Y$.
          \item See (i) for an example, that the condition
            ``projective'' for the morphism is not obsolete.
          \item Let $f:Y'\rightarrow Y$ be the blow up of a point of
            $\P^2_k$, then the fibres are connected, and hence
            $f_*\ko_{Y'}=\ko_Y$. $f$ is generically finite, but not
            finite and certainly no isomorphism.
          \item If $f:X\rightarrow\Spec(A)$ is a morphism such that
            $f^\#:A\rightarrow\Gamma(X,\ko_X)$ is an isomorphism, then
            $f$ need not be an isomorphism - despite the fact that
            $\Hom\big(X,\Spec(A)\big)$ and $\Hom\big(A,\Gamma(X,\ko_X)\big)$ are
            naturally bijective. E.~g.~$\P^1_k$ and $f$ the
            constant map to $\Spec(k)$. Then $\Gamma(X,\ko_X)=k$ and
            thus $f^\#=\id$ is an isomorphism, but $f$ is certainly
            not! If $X$ is also affine, then however $f$ is an
            isomorphism if and only if $f^{\#}$ is so.
          \item If $A\rightarrow B$ is a ring homomorphism and $B$ is
            locally free of rank $1$ as an $A$-module, and if for some 
            $P\in\Spec(A)$ the localisation $A_P\rightarrow B_P$ is an 
            isomorphism, then the
            homomorphism itself is an isomorphism. If, however, $A\rightarrow 
            B$ is only an $A$-module homomorphism, this need not be
            true. E.~g.~let $X$ be a smooth curve, e.~g.~$\P^1$, and
            $x\in X$ a point; let $\kl:=\ko(-x)$ be the ideal sheaf of
            $x$, in the case of $\P^1$ it's just $\ko(-1)$. This is a
            locally free sheaf, and $x$ is a divisor on the smooth
            curve $X$. This gives an exact sequence $0\rightarrow L
            \rightarrow\ko_X \rightarrow \ko_x \rightarrow 0$. For the 
            stalks over $y\not=x$ the map from $L\rightarrow\ko_X$ is
            of course an isomorphism, but not for $y=x$.
          \item Let $T=\A^1_k$, $X_t=V(x^2 - txy) \subset \P^1_k$ with
            homogeneous coordinates $x,y$ on $\P^1_k$. Then for
            $t_0=0$ the fibre $X_{t_0}$ is connected, namely a
            non-reduced fat point, while for $t\not=0$ the fibre $X_t$ 
            consists of two single points, and is thus not
            connected. The condition $X_{t_0}$ reduced in the
            Corollary \ref{cor:principle-of-connectedness} is hence vital.
       \end{enumerate}
     \end{remark}
     }

   \begin{proposition}[Principle of Connectedness]\label{prop:principle-of-connectedness}
     Let $X$ and $Y$ be noetherian schemes, $Y$ connected, and let
     $\pi:X\rightarrow Y$ be a flat projective morphism such that for 
     some $y_0\in Y$ the fibre $X_{y_0}=\pi^{-1}(y_0)$ is reduced and connected.

     Then for all $y\in Y$ the fibre $X_y=\pi^{-1}(y)$ is connected.
   \end{proposition}
   \begin{proof}
     Considering points in the intersections of the finite number of
     irreducible components of $Y$ we can reduce to the case $Y$
     irreducible. 
 
     Stein Factorisation (cf.~\cite{EGA} III.4.3.3\tom{ or
       \cite{Har77} III.11.5}) gives a factorisation of $\pi$ of the form 
     \begin{displaymath}
       \xymatrix{
         {\pi:X}\ar[r]^(0.35){\pi'} & Y'=\Spec(\pi_*\ko_X)\ar[r]^(0.7)f & Y,
         }
     \end{displaymath}
     with 
     \begin{enumerate}
        \item[(1)] $\pi'$ connected (i.~e.~its fibres are connected),
        \item[(2)] $f$ finite,
        \item[(3)] $f_*\ko_{Y'}=\pi_*\ko_X$ locally free over $\ko_Y$,
          since $\pi$ is flat, and
        \item[(4)] $Y'_{y_0}=f^{-1}(y_0)$ is connected and reduced,
          i.~e.~a single reduced point.
     \end{enumerate}

     Because of (1) it suffices to show that $f$ is connected, and we
     claim that they are reduced as well. 
     Since $f$ is finite (3) is equivalent to saying that $f$ is
     flat\tom{ (cf.~\cite{AM69} Proposition 3.10)}. Hence $f$
     fulfils the assumptions of Proposition \ref{prop:isomorphism}, and we
     conclude that $\ko_Y=f_*\ko_{Y'}$ and the proposition follows from
     \cite{Har77} III.11.3.\tom{\footnote{Compare the result with
         \cite{EGA} III.4.3.10, which deals with the case of $X$ also being  
         integral, but weakens the hypothesis on $\pi$ to proper and 
         universally open. Also the assumption on the fibre is reduced 
         to $R(Y)$ being algebraically closed in $R(X)$ - in the case 
         of the characteristic of $R(Y)$ is zero. Compare also \cite{EGA} IV.12.2.4 (vi)
         and III.15.5.9, which both deal with an arbitrary number of
         connected components.}}

     Alternatively, from \cite{EGA} IV.15.5.9 (ii) it follows that there
     is an open dense subset $U\subseteq Y$ such that $X_y$ is
     connected for all $y\in U$. Since, moreover, by the same theorem
     the number of connected components of the fibres is a lower
     semi-continuous function on $Y$ the special fibres cannot have
     more connected components than the general ones, that is, all
     fibres are connected.
   \end{proof}

   \tom{
   \begin{corollary}[Principle of Connectedness]\label{cor:principle-of-connectedness}
     Let $\{X_t\}_{t\in T}$ be a flat family of closed subschemes 
     $X_t\subseteq\P_{k(t)}^n$, where $T$ is a connected noetherian scheme.
     Suppose that $X_{t_0}$ is connected and reduced for some point
     $t_0\in T$. Then $X_t$ is connected for all $t\in T$.
   \end{corollary}
   \begin{proof}
     The $X_t$ are the fibres of a flat projective morphism
     \begin{displaymath}
       \xymatrix@C0.6cm{
         {\kx}\ar@{^{(}->}[r]\ar[dr]^{\pi} & {\P_T^n}\ar[d]\\
         & T,
         }
     \end{displaymath}
     \tom{\footnote{An open bracket had to be closed ).}}
     where $\kx$ is a closed subscheme of $\P_T^n$. 
     The result thus follows from Proposition
     \ref{prop:principle-of-connectedness}.\tom{\footnote{Compare the result with
         \cite{Har77} Ex.~III.11.4.}}
   \end{proof}
   }

   \begin{lemma}\label{lem:condition:A}
     Under the hypotheses of Proposition \ref{prop:condition} let
     $C\in|B|_a$ then $C$ is connected.  
   \end{lemma}
   \begin{proof}
     Consider the universal family
     \begin{equation}\label{eq:universal-family}
       \xymatrix@C1cm{
         |B|_a\times\Sigma\;\ar[dr]_{\pr_{|B|_a}} & {\;\bigcup\limits_{C\in |B|_a} \{C\}\times
         C=:S}\ar@{_{(}->}[l]\ar[d]^{\pi \mbox{ \scriptsize
             flat}}\\
         &|B|_a
         }
     \end{equation}
     \tom{\footnote{An open bracket had to be closed ).}}
     over the connected projective scheme
     $|B|_a\subseteq\Hilb_\Sigma$. 
 
     Since the projection $\pi$ is a flat projective morphism, and 
     since the fibre $\pi^{-1}(B)=\{B\}\times B$ is connected and reduced, the 
     result follows from Proposition \ref{prop:principle-of-connectedness}.
   \end{proof}

   \begin{lemma}\label{lem:condition:B}
     Under the hypotheses of Proposition \ref{prop:condition} let
     $C\in|B|_a$ with $B\subseteq C$, then $C=B$. 
   \end{lemma}
   \begin{proof}
     Suppose $B\subsetneqq C$, then the Hilbert polynomials of $B$ and
     $C$ are different in contradiction to $B\sim_a C$.
   \end{proof}

   \begin{lemma}\label{lem:condition:C}
     Under the hypotheses of Proposition \ref{prop:condition} let
     $C\in |B|_a$ with $C\not= B$, then $C\cap B=\emptyset$.
   \end{lemma}
   \begin{proof}
     Since $B$ is irreducible by Lemma \ref{lem:condition:B} $B$ and
     $C$ do not have a common component. Suppose $B\cap
     C=\{p_1,\ldots,p_r\}$, then $B^2=B.C\geq r>0$ in contradiction to
     $B^2=0$. 
   \end{proof}

   \begin{proposition}[Zariski's Lemma]\label{prop:zariskis-lemma}
     Under the hypotheses of Proposition \ref{prop:condition} let
     $C=\sum_{i=1}^r n_iC_i\in|B|_a$, where the $C_i$ are pairwise
     different irreducible curves, $n_i>0$ for $i=1,\ldots,r$. 

     Then the intersection matrix $Q=(C_i.C_j)_{i,j=1,\ldots,r}$ is
     negative semi-definite, and, moreover, $C$, considered as an
     element of the vectorspace $\bigoplus_{i=1}^r \Q\cdot C_i$,
     generates the annihilator of $Q$. 

     In particular, $D^2\leq 0$ for all curves $D\subseteq C$, and,
     moreover, $D^2=0$ if and only if $D=C$.
   \end{proposition}
   \begin{proof}
     By Lemma \ref{lem:condition:A} $C$ is connected.
     We are going to apply \cite{BPV84} I.2.10, and thus we have to
     verify three conditions.
     \begin{enumerate}
        \item[(i')] $C.C_i=B.C_i=0$ for all $i=1,\ldots,r$ by Lemma
          \ref{lem:condition:C}. Thus $C$
          is an element of the annihilator of $Q$ 
          with $n_i>0$ for all $i=1,\ldots,r$.
        \item[(ii)] $C_i.C_j\geq 0$ for all $i\not=j$.
        \item[(iii)] Since $C$ is connected there is no non-trivial
          partition $I\cup J$ of $\{1,\ldots r\}$ such that
          $C_i.C_j=0$ for all $i\in I$ and $j\in J$.
     \end{enumerate}
     Thus \cite{BPV84} I.2.10 implies that $-Q$ is positive
     semi-definite. 
   \end{proof}

   \begin{lemma}\label{lem:condition:D}
     Under the hypotheses of Proposition \ref{prop:condition} let $C,
     C'\in|B|_a$ be two distinct curves, then $C\cap C'=\emptyset$.
   \end{lemma}
   \begin{proof}
     Suppose $C=A+D$ and $C'=A+D'$ such that $D$ and $D'$ have 
     no common component. 

     We have
     \begin{displaymath}
       0=B^2=(A+D)^2=(A+D')^2=(A+D).(A+D'),
     \end{displaymath}
     and thus
     \begin{displaymath}
       (A+D)^2+(A+D')^2=2 (A+D).(A+D'),
     \end{displaymath}
     which implies that 
     \begin{displaymath}
       D^2+{D'}^2=2D.D',
     \end{displaymath}
     where each summand on the left hand side is less than or equal to zero by
     Proposition \ref{prop:zariskis-lemma}, and the right hand side is
     greater than or equal to zero, since the curves $D$ and $D'$ 
     have no common component. We thus conclude 
     \begin{displaymath}
       D^2={D'}^2=D.D'=0.
     \end{displaymath}
     But then again Proposition \ref{prop:zariskis-lemma} implies that 
     $D=C$ and $D'=C'$, that is, $C$ and $C'$ have no common
     component.

     Suppose $C\cap C'=\{p_1,\ldots,p_r\}$, then $B^2=C.C'\geq r>0$
     would be a contradiction to $B^2=0$. Hence, $C\cap C'=\emptyset$.
   \end{proof}

   \begin{lemma}\label{lem:condition:E}
     Under the hypotheses of Proposition \ref{prop:condition} 
     consider once more the universal family
     (\ref{eq:universal-family}) together with its projection onto 
     $\Sigma$,
     \begin{equation}
       \label{eq:universal-family-diagram}
       \xymatrix{
         |B|_a\times\Sigma 
         \ar@/_2pc/[ddr]_{\pr_{|B|_a}}\ar@/^2pc/[drr]^{\pr_\Sigma}&&\\
         & S \ar[r]^{\pi'}\ar[d]^{\pi}\ar@{_{(}->}[lu] &{\Sigma}\\
         & |B|_a. &
         }
     \end{equation}
     \tom{\footnote{An open bracket had to be closed ).}}
     Then $S$ is an irreducible projective surface, $|B|_a$
     is an irreducible curve, and $\pi'$ is surjective.
   \end{lemma}
   \begin{proof}
     \begin{varthm-roman}[Step 1]
       $S$ is an irreducible projective surface and $\pi'$ is surjective.
     \end{varthm-roman}
     The universal property of $|B|_a$ implies that $S$ is an
     effective Cartier divisor
     of $|B|_a\times\Sigma$, and thus in particular projective of
     dimension at least $2\leq\dim|B|_a+\dim(\Sigma)-1$. Since $\pi'$ is
     projective, its image is closed in $\Sigma$ and of dimension 2,
     hence it is the whole of $\Sigma$, since $\Sigma$ is irreducible.

     By Lemma \ref{lem:condition:D} the fibres of $\pi'$ are all
     single points, and thus, by \cite{Har92} Theorem 11.14, $S$ is
     irreducible. 

     Moreover, 
     \begin{displaymath}
       \dim(S)=\dim(\Sigma)+\dim(\mbox{fibre})=2.
     \end{displaymath}

     \begin{varthm-roman}[Step 2] 
       $\dim|B|_a=\dim(|B|_a\times\Sigma)-\dim(\Sigma)=\dim(S)+1-2=1.$
     \end{varthm-roman}

     \begin{varthm-roman}[Step 3]
       $|B|_a$ is irreducible.
     \end{varthm-roman}
     Let $V$ be any irreducible component of $|B|_a$ of dimension one, 
     then we have a universal family over $V$ and the analogue
     of Step 1 for $V$ shows that the curves in $V$ cover
     $\Sigma$. But then by Lemma \ref{lem:condition:D} there can be no 
     further curve in $|B|_a$, since any further curve would
     necessarily have a non-empty intersection with one of the curves
     in $V$. 
   \end{proof}

   \begin{lemma}\label{lem:condition:F}\tom{\footnote{The proof uses
         that the characteristic of the ground field is zero, even
         though one might perhaps avoid this.}} 
     Let's consider the following commutative diagram of projective
     morphisms
     \begin{equation}\label{eq:reduced-diagram}
       \xymatrix@C0.8cm{
         S\ar[r]^{\pi'}\ar[d]_{\pi} &{\Sigma}\\
         |B|_a & S_{red}\ar@{_{(}->}[lu]\ar[u]_{\varphi'}\ar[l]^{\varphi}
         }
     \end{equation}
     \tom{\footnote{An open bracket had to be closed ).}}
     The map ${\varphi'}:S_{red}\longrightarrow \Sigma$ 
     is birational.
   \end{lemma}
   \begin{proof}
     Since $S_{red}$ and $\Sigma$ are irreducible and reduced, and since 
     ${\varphi'}$ is surjective, we may apply \cite{Har77}
     III.10.5, and thus there is an open dense subset
     $U\subseteq S_{red}$ such that ${\varphi'}_|:U\rightarrow \Sigma$ is
     smooth. Hence, in particular ${\varphi'}_|$ is flat and the fibres are
     single reduced points.\tom{\footnote{By definition ${\varphi'}_|$ is
         \'etale and hence the completed local rings of the fibres are
         isomorphic to the completed local rings of their base points
         and hence regular. But then the local rings themselves are
         regular and thus reduced.}}
     Since ${\varphi'}_|:U\rightarrow{\varphi'}(U)$ is
     projective and quasi-finite, it is finite 
     (cf.~\cite{Har77} Ex.~III.11.2), and it follows from
     Proposition \ref{prop:isomorphism} that ${\varphi'}_|$ is an
     isomorphism onto its image, i.~e.~${\varphi'}$ is birational.
   \end{proof}

   \begin{lemma}\label{lem:condition:G}
     If $\psi:\Sigma\dashrightarrow S_{red}$ denotes the 
     rational inverse of the map $\varphi'$ in
     (\ref{eq:reduced-diagram}), then $\psi$ is indeed a
     morphism, i.~e.~$\varphi'$ is an isomorphism.
   \end{lemma}
   \begin{proof}
     By Lemma \ref{lem:condition:D} the fibres of $\varphi'$ over the
     possible points of indeterminacy of $\varphi'$ are just points,
     and thus the result follows from \cite{Bea83} Lemma~II.9.
   \end{proof}

   \begin{lemma}\label{lem:condition:H}
     The map $g:\Sigma\rightarrow |B|_a$ assigning to each point
     $p\in\Sigma$ the unique curve $C\in|B|_a$ with $p\in C$ is a
     morphism, and is thus a fibration whose fibres are the
     curves in $|B|_a$. 
   \end{lemma}
   \begin{proof}
     We just have $g=\varphi\circ\psi$.
   \end{proof}

   \begin{proof}[Proof of Proposition \ref{prop:condition}]
     Let $\nu:H\rightarrow |B|_a$ be the normalisation of the
     irreducible curve $|B|_a$. Then $H$ is a smooth irreducible
     curve.

     Moreover, since $\Sigma$ is irreducible and smooth, and since
     $g:\Sigma\rightarrow|B|_a$ is surjective, $g$ factorises over
     $H$, i.~e.~we have the following commutative diagram
     \begin{displaymath}
       \xymatrix@C0.8cm{
         {\Sigma}\ar[r]^g\ar[d]_{\exists f} &|B|_a\\
         H\ar[ur]^\nu.
         }
     \end{displaymath}
     Then $f$ is the desired fibration.
   \end{proof}


   \section{Some Facts used in the Proofs of Section \ref{sec:Geng-Xu}}
   \setcounter{equation}{0}

   In this section we are, in particular, writing down some 
   identifications of certain sheaves respectively of their global
   sections. Doing this we try to be very formal. However, in a
   situation of the kind 
   $\xymatrix@C0.8cm{X\ar@{^{(}->}[r]^i & Y\ar[r]^\pi& Z}$ we usually do 
   not distinguish between $\ko_X$ and $i_*\ko_X$, or between $\pi$
   and any restriction of $\pi$ to $X$.

   \begin{lemma}\label{lem:powerseries}
     Let $\varphi(x,y,t)=\sum_{i=0}^\infty \varphi_i(x,y)\cdot t^i
     \in\C\{x,y,t\}$ with $\varphi(x,y,t)\in (x,y)^m$ for every fixed
     $t$ in some small disc $\Delta$ around $0$. Then
     $\varphi_i(x,y)\in (x,y)^m$ for every $i\in \N_0$.
   \end{lemma}
   \begin{proof}
     We write the power series as
     $\varphi=\sum_{\alpha+\beta=0}^\infty(\sum_{i=0}^\infty
     c_{\alpha,\beta,i}\cdot t^i)x^\alpha y^\beta$.\\
     $\varphi(x,y,t)\in (x,y)^m$ for every $t\in\Delta$
     implies 
     \begin{displaymath}
        \sum_{i=0}^\infty c_{\alpha,\beta,i}\cdot t^i =0 \;\;\forall\;
        \alpha+\beta<m \mbox{ and } t\in\Delta.
     \end{displaymath}
     The identity theorem for power series in $\C$ then implies that 
     \begin{displaymath}
        c_{\alpha,\beta,i} = 0 \;\;\forall\;
        \alpha+\beta<m \mbox{ and } i\geq 0.
     \end{displaymath}
   \end{proof}

   \begin{lemma}\label{lem:identifications}
     Let $X$ be a noetherian scheme, $i:C\hookrightarrow X$ a closed subscheme, 
     $\kf$ a sheaf of modules on $C$, and $\kg$ a 
     sheaf of modules on $X$. Then
     \begin{equationlist}
        \item[eq:identifications:1] $i_*\kf\cong i_*\kf\otimes_{\ko_X}\ko_C$,
        \item[eq:identifications:2] $H^0(C,\kf)=H^0(X,i_*\kf\otimes_{\ko_X}\ko_C)$, 
        \item[eq:identifications:3] $\kg\otimes_{\ko_X}\ko_C\cong
          i_*i^*(\kg\otimes_{\ko_X}\ko_C)$, and
        \item[eq:identifications:4] $H^0\big(C,i^*(\kg\otimes_{\ko_X}\ko_C)\big)=H^0(X,\kg\otimes_{\ko_X}\ko_C)$.
     \end{equationlist}
   \end{lemma}
   \begin{varthm-roman-break}[Proof:]
     \begin{wideitemize}
        \item[(\ref{eq:identifications:1})] For $U\subseteq X$ open, we define
          \begin{displaymath}
            \begin{array}{ccccc}
              \Gamma(U,i_*\kf)&\rightarrow&
              \Gamma(U,i_*\kf)\otimes_{\Gamma(U,\ko_X)}\Gamma(U,\ko_C)&
              \subseteq&\Gamma(U,i_*\kf\otimes_{\ko_X}\ko_C)\\
              s&\mapsto&s\otimes 1.&&
            \end{array}
          \end{displaymath}
          This morphism induces on the stalks the isomorphism
          $$i_*\kf_x=\left\{
            \begin{array}{ccc}
              \kf_x, \mbox{ (if $x\in C$) }& = & \kf_x\otimes_{\ko_{X,x}}\ko_{X,x}/I_{C,x}\\
              0,\mbox{ (else) } & =& 0\otimes_{\ko_{X,x}}\ko_{X,x}/I_{C,x}
            \end{array}
          \right\}
          \cong i_*\kf_x\otimes_{\ko_{X,x}}\ko_{C,x},$$ 
          where
          $I_{C,x}$ is the ideal defining $C$ in $X$ locally at $x$.
        \item[(\ref{eq:identifications:2})] The identification 
          (\ref{eq:identifications:1}) together with \cite{Har77} III.2.10 gives:
          \begin{displaymath}
            H^0(C,\kf)=H^0(X,i_*\kf)=H^0(X,i_*\kf\otimes_{\ko_X}\ko_C).
          \end{displaymath}
        \item[(\ref{eq:identifications:3})] The adjoint property of $i_*$ and $i^*$ together with
          $i^*i_*\cong \id$ gives rise to the 
          following isomorphisms:
          \begin{displaymath}
            \End\big(i_*i^*(\kg\otimes\ko_C)\big)\cong
            \Hom\big(i^*i_*i^*(\kg\otimes\ko_C),i^*(\kg\otimes\ko_C)\big)
          \end{displaymath}
          \begin{displaymath}
            \cong\End\big(i^*(\kg\otimes\ko_C)\big)
            \cong\Hom\big(\kg\otimes\ko_C,i_*i^*(\kg\otimes\ko_C)\big).
          \end{displaymath}
          That means, that the identity morphism on
          $i_*i^*(\kg\otimes\ko_C)$ must correspond to an isomorphism 
          from $\kg\otimes\ko_C$ to $i_*i^*(\kg\otimes\ko_C)$ via 
          these identifications.
        \item[(\ref{eq:identifications:4})] follows from
          (\ref{eq:identifications:3}) and once more \cite{Har77}
          III.2.10. 
     \end{wideitemize}
   \end{varthm-roman-break}

   \begin{corollary}\label{cor:identifications}
     In the situation of Lemma \ref{lem:lifting} we have:
     \begin{equationlist}
        \item[eq:identifications:5] $H^0\big(C,\pi_*\ko_{\widetilde{C}}(E)\otimes_{\ko_C}\ko_C(C)\big)
       =H^0\big(\Sigma,\pi_*\ko_{\widetilde{C}}(E)\otimes_{\ko_\Sigma}\ko_C(C)\big)$, and
     \item[eq:identifications:6] $H^0\big(C,\pi_*\ko_{\widetilde{\Sigma}}(E)\otimes_{\ko_\Sigma}\ko_C(C)\big)
       =H^0\big(\Sigma,\pi_*\ko_{\widetilde{\Sigma}}(E)\otimes_{\ko_\Sigma}\ko_C(C)\big)$.
     \end{equationlist}
   \end{corollary}
   \begin{proof}
     We denote by $j:\widetilde{C}\hookrightarrow \widetilde{\Sigma}$ and
     $i:C\hookrightarrow \Sigma$ respectively the given
     embeddings.
     \begin{wideitemize}
        \item[(\ref{eq:identifications:5})]  By 
          (\ref{eq:identifications:2}) in Lemma
          \ref{lem:identifications} we have:
          \begin{displaymath}
            H^0\big(C,\pi_*\ko_{\widetilde{C}}(E)\otimes_{\ko_C}\ko_C(C)\big)=
            H^0\Big(\Sigma,i_*\big(\pi_*\ko_{\widetilde{C}}(E)\otimes_{\ko_C}\ko_C(C)\big)
            \otimes_{\ko_\Sigma}\ko_C\Big).
          \end{displaymath}
          By the projection formula this is just equal to:
          \begin{displaymath}
            H^0\Big(\Sigma,\big(i_*\pi_*\ko_{\widetilde{C}}(E)\otimes_{\ko_\Sigma}\ko_\Sigma(C)\big)
            \otimes_{\ko_\Sigma}\ko_C\Big)
            =H^0\big(\Sigma,\pi_*j_*\ko_{\widetilde{C}}(E)\otimes_{\ko_\Sigma}\ko_C(C)\big)
          \end{displaymath}
          \begin{displaymath}
            =_{\mbox{\footnotesize def}}
            H^0\big(\Sigma,\pi_*\ko_{\widetilde{C}}(E)\otimes_{\ko_\Sigma}\ko_C(C)\big).
          \end{displaymath}
        \item[(\ref{eq:identifications:6})] Using
          (\ref{eq:identifications:4}) in Lemma \ref{lem:identifications} we get:
          \begin{displaymath}
            H^0\big(C,\pi_*\ko_{\widetilde{\Sigma}}(E)\otimes_{\ko_\Sigma}\ko_C(C)\big)=_{\mbox{\footnotesize def}}
            H^0\Big(C,i^*\big(\pi_*\ko_{\widetilde{\Sigma}}(E)\otimes_{\ko_\Sigma}\ko_C(C)\big)\Big)=
          \end{displaymath}
          \begin{displaymath}
            H^0\big(\Sigma,\pi_*\ko_{\widetilde{\Sigma}}(E)\otimes_{\ko_\Sigma}\ko_C(C)\big).
          \end{displaymath}
     \end{wideitemize}
   \end{proof}

   \tom{
   \begin{lemma}\label{lem:support:1}
     Let $\kf$ be any coherent sheaf on $\Sigma$. Then the kernel of
     the natural map
     \begin{displaymath}
       \xymatrix@C0.8cm{
         {\bigotimes_{i=1}^r \kj_{X(m_i;z_i)/\Sigma}\otimes\kf}\ar[r]^(0.58){\delta} &
         {\kj_{X(\underline{m};\underline{z})/\Sigma}\otimes\kf}
         }
     \end{displaymath}
     has support contained in $\{z_1,\ldots,z_r\}$.
   \end{lemma}
   \begin{proof}
     We have $\Ker(\delta)_z=\Ker(\delta_z)$ and for
     $z\not\in\{z_1,\ldots,z_r\}$ the map $\delta_z$ is given by
     \begin{displaymath}
       \xymatrix@C0.6cm@R0.1cm{
         {\bigotimes_{i=1}^r\ko_{\Sigma,z}\otimes\kf_z}\ar[r]&
         {\ko_{\Sigma,z}\otimes\kf_z}\\
         f_1\otimes\cdots\otimes f_r\otimes g\ar@{|->}[r] &
         f_1\cdots f_r\otimes g,
         }
     \end{displaymath}
     and is thus an isomorphism.
   \end{proof}
   }

   \begin{lemma}\label{lem:support:2}
     With the notation of Lemma \ref{lem:lifting} we show that
     $\supp\big(\Ker(\gamma)\big)\subseteq\{z_1,\ldots,z_r\}$.
   \end{lemma}
   \begin{proof}
     Since $\pi:\widetilde{\Sigma}\setminus\big(\bigcup_{i=1}^{r}E_i\big)\longrightarrow
     \Sigma\setminus\{z_1,\ldots,z_r\}$ is an isomorphism, we have for any
     sheaf $\kf$ of $\ko_{\widetilde{\Sigma}}$-modules and $y\in
     \widetilde{\Sigma}\setminus\big(\bigcup_{i=1}^{r}E_i\big)$: 
     \begin{displaymath}
       (\pi_*\kf)_{\pi(y)}=\lim_{\pi(y)\in
         V}\kf(\pi^{-1}(V))=\lim_{y\in
         U}\kf(U)=\kf_y.
     \end{displaymath}
     In particular, 
     \begin{displaymath}
       \big(\pi_*\ko_{\widetilde{\Sigma}}(E)\otimes_{\ko_\Sigma}\ko_C(C)\big)_{\pi(y)}\cong
       \ko_{\widetilde{\Sigma},y}\otimes_{\ko_{\Sigma,\pi(y)}}\ko_{C,\pi(y)}\cong
       \ko_{\Sigma,\pi(y)}\otimes_{\ko_{\Sigma,\pi(y)}}\ko_{C,\pi(y)},
     \end{displaymath}
     and
     \begin{displaymath}
       \big(\pi_*\ko_{\widetilde{C}}(E)\otimes_{\ko_\Sigma}\ko_C(C)\big)_{\pi(y)}\cong
       \ko_{\widetilde{C},y}\otimes_{\ko_{\Sigma,\pi(y)}}\ko_{C,\pi(y)}\cong
       \ko_{C,\pi(y)}\otimes_{\ko_{\Sigma,\pi(y)}}\ko_{C,\pi(y)}.
     \end{displaymath}
     Moreover, the morphism $\gamma_{\pi(y)}$ becomes under these
     identifications just the morphism given by
     $a\otimes\overline{b}=1\otimes\overline{ab}\mapsto\overline{a}\otimes\overline{b}=1\otimes\overline{ab}$, 
     which is injective. Thus,
     $0=\Ker(\gamma_{\pi(y)})=\Ker(\gamma)_{\pi(y)}$, and
     $\pi(y)\not\in\supp(\Ker(\gamma))$.     
   \end{proof}

   \begin{lemma}\label{lem:torsion}
     Let $X$ be an irreducible noetherian scheme, 
     $\kf$ a coherent sheaf on $X$, and $s\in
     H^0(X,\kf)$ such that $\dim\big(\supp(s)\big)<\dim(X)$.
     Then $s\in H^0\big(X,\Tor(\kf)\big)$.
   \end{lemma}
   \begin{proof}
     The multiplication by $s$ gives rise to the following exact
     sequence:
     \begin{displaymath}
       \xymatrix@C0.6cm{
         0\ar[r] & {\Ker(\cdot s)} \ar[r] & {\ko_X}\ar[r]^{\cdot
         s} & {\kf}.
       }
     \end{displaymath}  
     Since $\ko_X$ and $\kf$ are coherent, so is $\Ker(\cdot s)$, and
     hence $\supp\big(\Ker(\cdot s)\big)$ is closed in $X$.
     Now,  
     \begin{displaymath}
       \supp(\Ker(\cdot s))=\{z\in X\;|\;\exists\;
       0\not=r_z\in\ko_{X,z}:r_z\cdot s_z=0\}
     \end{displaymath}
     \begin{displaymath}
       =\{z\in X\;|\;s_z\in \Tor(\kf_z)\}.       
     \end{displaymath}
     But then the complement $\{z\in X|s_z\not\in \Tor(\kf_z)\}$ is
     open and is contained in $\supp(s)$ (since $s_z=0$ implies that
     $s_z\in\Tor(\kf_z)$), and is thus empty since $X$ is irreducible 
     and $\supp(s)$ of lower dimension. 
   \end{proof}

   \section{The Degree of a Line Bundle on a Curve}
   \setcounter{equation}{0}

   \begin{remark}
     Let $C=C_1\cup\ldots\cup C_k$ be a reduced curve on a
     smooth projective surface $\Sigma$ over $\C$, where the $C_i$ are irreducible, 
     and let $\kl$ be a line bundle on $C$. Then we define the \emph{degree}
     of $\kl$ with the aid of the normalisation
     $\nu:C'\rightarrow C$. We have
     $H^2(C,\Z)\cong\bigoplus_{i=1}^kH^2(C'_i,\Z)=\Z^k$, and
     thus the image of $\kl$ in $H^2(C,\Z)$, which is the first Chern
     class of $\kl$, can be viewed as a vector $(l_1,\ldots,l_k)$ of
     integers, and we may define the degree of $\kl$ by
     \begin{displaymath}
       \deg(\kl):=l_1+\cdots+l_k.
     \end{displaymath}
     In particular, if $C$ is irreducible, we get:
     \begin{displaymath}
       \deg(\kl)=\deg(\nu^*\kl)=c_1(\nu^*\kl).
     \end{displaymath}
     Since $H^0(C,\kl)\not=0$ implies that
     $H^0(C',\nu^*\kl)\not=0$, and since the existence of a
     non-vanishing global section of $\nu^*\kl$ on the smooth curve
     $C'$ implies that the corresponding divisor is effective,
     we get the following lemma. (cf.~\cite{BPV84} Section~II.2)
   \end{remark}

   \begin{lemma}\label{lem:degree}
     Let $C$ be an irreducible reduced curve on a
     smooth projective surface $\Sigma$, 
     and let $\kl$ be a line bundle on $C$.
     If $H^0(C,\kl)\not=0$, then $\deg(\kl)\geq 0$.
   \end{lemma}


   \section{Two Results used in the Proof of Theorem \ref{thm:existence-I}}
   \setcounter{equation}{0}

   \begin{lemma}\label{lem:smoothcurves}
     Let $L$ be very ample over $\C$ on the smooth projective
     surface $\Sigma$, and let $z,z'\in\Sigma$ be two distinct
     points. 
     Then there is a smooth curve through $z$ and $z'$ in $|L|_l$.
   \end{lemma}
   \begin{proof}\tom{\footnote{See also \cite{Har77} Proof of
         Corollary III.9.13.}}
     Considering the embedding into $\PC^n$ defined by
     $L$ there is an $n-2$-di\-men\-sio\-nal family of hyperplane sections going 
     through two fixed points of $\Sigma$, which in local
     coordinates w.~l.~o.~g.~is given by the family of equations
     $\kf=\big\{a_1x_1+\ldots+a_{n-1}x_{n-1}=0\;\big|\;
     (a_1:\ldots:a_{n-1})\in\PC^{n-2}\big\}$. Since the local
     analytic rings of $\Sigma$ in every point are smooth,
     hence, in particular complete intersections, they are given
     as $\C\{x_1,\ldots,x_n\}$ modulo some ideal generated by
     $n-2$  power series $f_1,\ldots,f_{n-2}$ forming a regular
     sequence. Thus, 
     having $n-2$ free indeterminates in our family $\kf$ of equations,
     a generic equation $g$ will lead to a regular
     sequence $f_1,\ldots,f_{n-2},g$, i.~e.~the hyperplane section 
     defined by $g$ is smooth in each of the two points, and thus
     everywhere.
   \end{proof}

   \begin{lemma}\label{lem:inductionstep}
     Let $L\subset \Sigma$ be a smooth curve and $X\subset \Sigma$  
     a zero-dimensional scheme. If $D\in\Div(\Sigma)$ such that
     \begin{equationlist}
        \item[eq:inductionstep:1] $h^1\big(\Sigma,\kj_{X:L/\Sigma}(D-L)\big)=0$, and
        \item[eq:inductionstep:2] $\deg(X\cap L)\leq D.L+1-2 g(L)$,
     \end{equationlist}
     then
     \begin{displaymath}
       h^1\big(\Sigma,\kj_{X/\Sigma}(D)\big)=0.
     \end{displaymath}
   \end{lemma}

   \begin{proof}
     Condition (\ref{eq:inductionstep:2})) implies
     \begin{displaymath}
       2 g(L)-2 < D.L-\deg(X\cap L) = \deg\big(\ko_L(D)\big) +
       \deg\big(\kj_{X\cap L/L}\big) 
     \end{displaymath}
     \begin{displaymath}
       =\deg\big(\kj_{X\cap
         L/L}(D)\big),
     \end{displaymath}
     and thus by Riemann-Roch (cf.~\cite{Har77} IV.1.3.4)
     \begin{displaymath}
       h^1\big(\kj_{X\cap L/L}(D)\big)=0.
     \end{displaymath}
     Consider now the exact sequence
     \begin{displaymath}
       \xymatrix@C0.6cm{
         0\ar[r] & {\kj_{X:L/\Sigma}(D-L)}\ar[r]^(0.6){\cdot L} & {\kj_{X/\Sigma}(D)}\ar[r]
         &{\kj_{X\cap L/L}(D)}\ar[r] &0.
         }
     \end{displaymath}
     The result then follows from the corresponding long exact cohomology 
     sequence
     \begin{displaymath}
       \xymatrix@C0.6cm{
         0 = H^1\big(\kj_{X:L/\Sigma}(D-L)\big)\ar[r] &
         H^1\big(\kj_{X/\Sigma}(D)\big)\ar[r] &
         H^1\big(\kj_{X\cap L/L}(D)\big)= 0.
         }
     \end{displaymath}
   \end{proof}


   \section{Product Surfaces}\label{app:product-surfaces}
   \setcounter{equation}{0}


   Througout this section we stick to the notation of Section
   \ref{subsec:product-curves} and \ref{subsec:elliptic-curves}.
   Let $C_1$ and $C_2$ be two smooth projective curves of genus
   $g_1\geq 0$
   and $g_2\geq 0$ respectively, and let $\Sigma=C_1\times C_2$. 

   Supposed that one of the curves is rational, the surface is geometrically
   ruled and the Picard number of $\Sigma$ is two. Whereas, if both $C_1$ and $C_2$
   are of strictly positive genus, this need no longer be the
   case as we have seen in Remark \ref{rem:isogenous}. 
   Thus the following proposition is the best we may expect.

   \begin{proposition}\label{prop:product-surfaces}
     For a generic choice of smooth projective curves $C_1$ and $C_2$
     the Neron-Severi group of $\Sigma=C_1\times C_2$ is $\NS(\Sigma)\cong
     C_1\Z\oplus C_2\Z$.

     More precisely, fixing $g_1$ and $g_2$ there is a very general
     subset $U\subseteq M_{g_1}\times M_{g_2}$ such that for any
     $(C_1,C_2)\in U$ the Picard number of $C_1\times C_2$ is two,
     where $M_{g_i}$ denotes the moduli space of smooth projective
     curves of genus $g_i$, $i=1,2$.
   \end{proposition}
   \begin{proof}
     As already mentioned, if either $g_1$ or $g_2$ is zero, then we may take
     $U=M_{g_1}\times M_{g_2}$. 

     Suppose that $g_1=g_2=1$. Given an elliptic curve $C_1$ there is
     a countable union $V$ of proper subvarieties of $M_1$ such that for
     any $C_2\in M_1\setminus V$ the Picard number of $C_1\times C_2$
     is two - namely, if $\tau_1$ and $\tau_2$ denote the periods as
     in Section \ref{subsec:elliptic-curves}, then we have to require
     that there exists no invertible integer matrix $\left(
     \begin{smallmatrix}
       z_1&z_2\\z_3&z_4
     \end{smallmatrix}\right)$
     such that $\tau_2 = \frac{z_4-z_3\tau_1}{z_2-z_1\tau_1}$. (Compare
     also \cite{GH94} p.~286.)

     We, therefore, may assume that $g_1\geq 2$ and $g_2\geq 1$.
     The claim then follows from Lemma \ref{lem:product-surfaces},
     which is due to Denis Gaitsgory. 
   \end{proof}

   \begin{lemma}[Denis Gaitsgory]\label{lem:product-surfaces}
     Let $C_2$ be any smooth projective curve of genus $g_2\geq 1$. Then
     for any $g_1\geq 2$ there is a very general subset $U$ of the moduli space
     $M_{g_1}$ of smooth projective curves of genus $g_1$ such that the
     Picard number of $C_1\times C_2$ is two for any $C_1\in U$. 
   \end{lemma}

   \begin{proof}
     We note that a curve $B\subset \Sigma=C_1\times C_2$ with
     $C_1\not\sim_a B\not\sim_a C_2$ induces a non-trivial
     morphism
     $\mu_B:C_1\rightarrow\Pic(C_2):p\mapsto\pr_{2*}\big(\pr_1^*(p)\big)$, 
     where $\pr_i:\Sigma\rightarrow C_i$, $i=1,2$,
     denote the canonical projections. It thus makes sense to study
     the moduli problem of (non-trivial) maps from curves of genus $g_1$
     into $\Pic(C_2)$. 

     More precisely, let $k\in\N$ and let
     $0\not=\beta\in H_2\big(\Pic_k(C_2),\Z\big)=\Z^{2g_2}$ be
     given, where $\Pic_k(C_2)$ is the Picard variety of divisors of
     degree $k$ on $C_2$. Following the
     notation of \cite{FP97} we 
     denote by $M_{g_1,0}\big(\Pic_k(C_2),\beta\big)$ the moduli space
     of pairs $(C_1,\mu)$, where $C_1$ is a smooth projective curve of
     genus $g_1$ and $\mu:C_1\rightarrow\Pic_k(C_2)$ a morphism with
     $\mu_*\big([C_1]\big)=\beta$.  We then have the
     canonical morphism 
     \begin{displaymath}
       F_{k,\beta}:M_{g_1,0}\big(\Pic_k(C_2),\beta\big)\rightarrow
       M_{g_1}:(C_1,\mu)\mapsto C_1,
     \end{displaymath}
     just forgetting the map $\mu$, and the proposition reduces to
     the following claim:

     \begin{claim}
       For no choice of $k\in\N$ and
       $0\not=\beta\in H_2\big(\Pic_k(C_2),\Z\big)$ the morphism
       $F_{k,\beta}$ is dominant.
     \end{claim}

     Let $\mu:C_1\rightarrow \Pic_k(C_2)$ be any morphism with
     $\mu_*\big([C_1]\big)=\beta$. Then $\mu$ is not a contraction and
     the image of $C_1$ is a projective
     curve in the abelian variety $\Pic_k(C_2)$. Moreover, we have the
     following exact sequence of sheaves
     \begin{equation}\label{eq:curveproduct}
       \xymatrix{
         0\ar[r]& \kt_{C_1}\ar[r]^(0.3){d\mu}&
         \mu^*\kt_{\Pic_k(C_2)}=\ko_{C_1}^{g_2}\ar[r]& 
         \kn_\mu:=\coker(d\mu)\ar[r]&0.
         }
     \end{equation}
     Since $d\mu$ is a non-zero inclusion, its dual
     $d\mu^\vee:(\mu^*\kt_{\Pic_k(C_2)})^\vee=\ko_{C_1}^{g_2} \rightarrow
     \Omega_{C_1}=\omega_{C_1}$ is not zero on global sections, that is
     \begin{displaymath}
       H^0\big(d\mu^\vee\big) :
       H^0\big(C_1,\ko_{C_1}^{g_2}\big)=\Hom_{\ko_{C_1}}\big(\ko_{C_1}^{g_2},\ko_{C_1}\big)
       \rightarrow
       H^0(C_1,\omega_{C_1})=\Hom_{\ko_{C_1}}\big(\kt_{C_1},\ko_{C_1}\big)
     \end{displaymath}
     is not the zero map. Since $g_1\geq 2$ we have
     $h^0(C_1,\omega_{C_1})=2g_1-2>0$, and thus $\omega_{C_1}$ has global
     sections. Therefore, the induced map
     $H^0\big(C_1,\omega_{C_1}\otimes\ko_{C_1}^{g_2}\big)
     \rightarrow H^0(C_1,\omega_{C_1}\otimes\omega_{C_1})$ is not the zero map,
     which by Serre duality gives that the map 
     \begin{displaymath}
       H^1(d\mu):H^1(C_1,\kt_{C_1})\rightarrow H^1\big(C_1,\mu^*\kt_{\Pic_k(C_2)}\big)
     \end{displaymath}
     from the long exact cohomology sequence of
     (\ref{eq:curveproduct}) is not zero. Hence the coboundary map
     \begin{displaymath}
       \delta:H^0(C_1,\kn_\mu)\rightarrow H^1(C_1,\kt_{C_1})
     \end{displaymath}
     cannot be surjective. According to \cite{Har98} p.~96 
     we have
     \begin{displaymath}
       \delta=dF_{k,\beta}:\kt_{M_{g_1,0}(\Pic_k(C_2),\beta)}
       =H^0(C_1,\kn_\mu)\longrightarrow \kt_{M_{g_1}}=H^1(C_1,\kt_{C_1}).
     \end{displaymath}
     But if the differential of $F_{k,\beta}$ is not surjective, then
     $F_{k,\beta}$ itself cannot be dominant.
   \end{proof}

\end{appendix}

   \bibliographystyle{amsalpha}
   \bibliography{bibliographie}
   \nocite{AM69}
   \nocite{GLS97}
   \nocite{GLS98b}
   \nocite{GLS00}
   \nocite{GS99}
   \nocite{Los99}
   \nocite{Shu91b}
   \nocite{Shu96b}
   \nocite{ST96}

\end{document}